\documentclass[10pt]{article}
\usepackage[mathscr]{eucal}
\usepackage{mathrsfs}
\usepackage{amsfonts}
\usepackage{amsmath}
\usepackage{amsthm}
\usepackage{amssymb}
\usepackage{latexsym}
\usepackage[all]{xy}
\usepackage{pstricks,pst-node}

\theoremstyle{plain}
\newtheorem{definition}[equation]{Definition}
\newtheorem{corollary}[equation]{Corollary}
\newtheorem{lemma}[equation]{Lemma}
\newtheorem{proposition}[equation]{Proposition}
\newtheorem{theorem}[equation]{Theorem}

\theoremstyle{remark}
\newtheorem{remark}[equation]{Remark}

\numberwithin{equation}{subsection}

\makeatletter
\renewcommand{\subsection}{\@startsection{subsection}{2}{0pt}{-3ex
plus -1ex minus -0.2ex}{-2mm plus -0pt minus
-2pt}{\normalfont\bfseries}} \makeatother

\newcommand{\Lmod}[1]{#1\text{-}\mathsf{mod}}

\DeclareMathOperator{\Res}{\mathsf{Res}}
\DeclareMathOperator{\Ker}{\mathrm{Ker}}
\DeclareMathOperator{\gr}{\mathrm{gr}}
\DeclareMathOperator{\Lie}{\mathrm{Lie}}
\DeclareMathOperator{\Tr}{\mathrm{Tr}}
\DeclareMathOperator{\End}{\mathrm{End}}
\DeclareMathOperator{\Hom}{\mathrm{Hom}}
\DeclareMathOperator{\Hilb}{\mathrm{Hilb}}

\DeclareMathOperator{\diag}{\mathrm{diag}}
\newcommand{\opp}{\operatorname}
\DeclareMathOperator{\im}{\mathrm{Im}}
\newcommand{\step}[1]{\noindent\vskip 2pt{{\text{\sc{Step }}#1}\en}}
\newcommand{\Mod}[1]{#1\text{-}\mathsf{mod}}
\newcommand{\vi}{${\sf {(i)}}\;$}
\newcommand{\vii}{${\sf {(ii)}}\;$}
\newcommand{\viii}{${\sf {(iii)}}\;$}

\newcommand{\beq}{\begin{equation}\label}
\newcommand{\eeq}{\end{equation}}
\newcommand{\en}{{\enspace}}
\newcommand{\CC}{\mathbb C}
\newcommand{\ZZ}{\mathbb Z}
\newcommand{\solu}[1]{\begin{sol}{\bf (\ref{#1})}}
\newcommand{\dd}{\partial}

\def\ccirc{{{}_{\,{}^{^\circ}}}}
\newcommand{\sset}{\subset}

\newcommand{\Ann}{{\mathtt{Ann}^{\,}}}
\newcommand{\into}{{}^{\,}\hookrightarrow^{\,}}
\newcommand{\too}{\,\longrightarrow\,}
\newcommand{\mto}{{\mapsto}}
\newcommand{\onto}{\twoheadrightarrow}

\newcommand{\iso}{{\;\stackrel{_\sim}{\to}\;}}
\newcommand{\cd}{\!\cdot\!}
\newcommand{\id}{{{\mathtt{Id}}}}
\newcommand{\inv}{^{-1}}
\newcommand{\ad}{{\mathtt{{ad}}}}
\newcommand{\Ad}{{\mathtt{{Ad}}}}

\newcommand{\GL}{\operatorname{GL}}

\newcommand{\Sp}{\operatorname{Sp}}
\newcommand{\hh}{{\mathsf{H}}}
\newcommand{\e}{{\mathbf{e}}}
\newcommand{\bone}{{\boldsymbol{1}}}
\newcommand{\eu}{{\mathsf{eu}}}
\newcommand{\la}{\lambda}
\newcommand{\op}{{\mathsf{op}}}
\newcommand{\QQ}{{\overline{Q}}}
\newcommand{\Qc}{{Q_{\mathsf{CM}}}}
\newcommand{\QQc}{{\overline{Q_{\mathsf{CM}}}}}
\newcommand{\al}{{\alpha}}
\newcommand{\bbe}{{\boldsymbol{\epsilon}}}
\newcommand{\bxi}{{\boldsymbol{\xi}}}
\newcommand{\bmu}{{\boldsymbol{\mu}}}
\newcommand{\bkap}{{\boldsymbol{\kappa}}}
\newcommand{\eps}{{\epsilon}}
\newcommand{\bept}{{\boldsymbol{\varepsilon}_{_\TT}}}
\newcommand{\ga}{{\gamma}}
\newcommand{\om}{\omega}
\newcommand{\rt}{\rtimes }
\newcommand{\half}{{\mbox{$\frac{1}{2}$}}}
\newcommand{\ttimes}{{\rtimes_{_{\mathbf{S}}}}}
\newcommand{\beps}{{{\mathbf{\varepsilon}_{_{\mathbf{\Gamma}}}}}}
\newcommand{\BH}{{\mathbb{H}}}

\newcommand{\at}{{{\mathsf{a}}}}
\renewcommand{\ss}{{{\mathbf{S}}}}
\newcommand{\GG}{{{\mathbf{\Gamma}}_n}}
\newcommand{\GGG}{{{\mathbf{\Gamma}}_{n-1}}}
\newcommand{\bis}{_{(1)}}
\newcommand{\omg}{\omega(\gamma}
\newcommand{\ttt}{{\hh_{(1)}}}
\newcommand{\hol}{{\mathscr{H}\!\text{\it ol}}}
\newcommand{\ft}{{\mathfrak{t}}}
\newcommand{\ls}{{L^{\times}}}
\newcommand{\eel}{{{\mathsf e}_T^n}}

\renewcommand{\k}{\Bbbk}
\newcommand{\BQ}{{\mathbb Q}}
\newcommand{\BF}{{\mathbb F}}
\newcommand{\FA}{{\mathscr A}}
\newcommand{\FM}{{\mathsf{M}}}
\newcommand{\N}{{\mathsf{N}}}
\newcommand{\A}{{\mathfrak{A}}}
\newcommand{\UU}{{\mathsf{U}}}

\newcommand{\ot}{{\otimes}}
\newcommand{\otot}{{\otimes\cdots\otimes}}
\newcommand{\wt}{\widetilde}
\newcommand{\Th}{{\Theta}}
\newcommand{\LL}{{L_{\text{reg}}}}
\newcommand{\FF}{{F^{\ot n}}}
\newcommand{\M}{{\mathsf{B}}}
\newcommand{\bbM}{{\mathbb{M}}}
\renewcommand{\O}{{\mathcal{O}}}
\newcommand{\T}{{\mathfrak{T}}}
\newcommand{\D}{{\mathscr{D}}}
\newcommand{\TD}{{\mathscr{D}}}
\renewcommand{\P}{{\mathbb{P}}}
\newcommand{\pr}{{\mathrm{pr}}}

\newcommand{\F}{{\mathsf{N}}}
\newcommand{\Z}{{\mathcal{Z}}}
\newcommand{\Zc}{{\mathcal{Z}_{\mathsf{CM}}}}

\renewcommand{\mod}{{\mathrm{mod}}}

\newcommand{\sfA}{{\mathsf{A}}}
\newcommand{\sfR}{{\mathsf{R}}}
\newcommand{\sfE}{{\mathsf{E}}}

\def\C{{\mathbb{C}}}
\def\Rep{{\mathsf{Rep}}}
\def\g{{\mathfrak{g}}}

\def\gl{{\mathfrak{g}\mathfrak{l}}}
\def\sln{{\mathfrak{s}\mathfrak{l}}_n}

\def\C{{\mathbb{C}}}
\def\Rep{{\mathsf{Rep}}}
\def\g{{\mathfrak{g}}}

\def\gl{{\mathfrak{g}\mathfrak{l}}}
\def\sln{{\mathfrak{s}\mathfrak{l}}_n}

\def\JJ{{\mathbb{J}}}

\def\pa{{\partial}}
\def\G{{\Gamma}}
\def\Ic{{I_{\mathsf{CM}}}}

\def\de{\delta}
\def\dq{\overline{Q}}
\def\al{\alpha}
\def\be{\beta}
\def\th{\theta}
\def\Th{\Theta}

\def\hp{\hphantom{x}}

\def\bbt{{h}}
\def\btt{{\mathbf{t}}}
\def\TT{{\mathbb{T}}}
\def\X{{\mathfrak{X}}}

\def\sminus{{\smallsetminus}}
\newcommand{\dis}{\displaystyle}
\newcommand{\lreg}{{L^n_{\text{reg}}}}
\newcommand{\reg}{_{\text{reg}}}

\begin{document}

\centerline{\Large {\textbf{Harish-Chandra
homomorphisms and}}}
\vskip 3pt
\centerline{\Large {\textbf{symplectic reflection algebras for
wreath-products }}}
\bigskip
\centerline{\sc {\large{Pavel Etingof, Wee Liang Gan, Victor Ginzburg, Alexei Oblomkov}}}
\vskip 9mm

\hfill{\em To Joseph Bernstein on the occasion of his 60th
Birthday}\break
\bigskip
\begin{abstract}
The main result of the paper is  a natural
 construction of the spherical subalgebra in a symplectic reflection algebra
associated with a wreath-product in terms of quantum hamiltonian
reduction of an algebra of differential operators
on a representation space of an extended Dynkin
quiver. The existence of such a construction has been conjectured in
\cite{EG}.

We also present a new approach to reflection functors
and shift functors for generalized preprojective algebras and 
 symplectic reflection algebras
associated with wreath-products.
\end{abstract}
\bigskip

\centerline{\sf Table of Contents}
\vskip -1mm

$\hspace{30mm}$ {\footnotesize \parbox[t]{115mm}{
\hp${}_{}$\!
\hp\!1.{ $\;\,$} {\tt Introduction} \newline
\hp2.{ $\;\,$} {\tt Calogero-Moser quiver} \newline
\hp3.{ $\;\,$} {\tt Radial part map} \newline
\hp4.{ $\;\,$} {\tt Dunkl representation} \newline
\hp5.{ $\;\,$} {\tt Harish-Chandra homomorphism} \newline
\hp6.{ $\;\,$} {\tt Reflection isomorphisms} \newline
\hp7.{ $\;\,$} {\tt Appendix A: Extended Dynkin quiver} \newline
\hp8.{ $\;\,$} {\tt Appendix B: Proof of Proposition \ref{radialpart}} \newline
\hp9.{ $\;\,$}  {\tt Appendix C: Proof of Theorem \ref{mess}}}
}


\section{Introduction}
The main result of the paper is   the proof of
\cite[Conjecture 11.22]{EG} that provides a natural
 construction of the spherical subalgebra in a symplectic reflection algebra
associated with a wreath-product in terms of quantum hamiltonian
reduction of an algebra of differential operators.

To state the main result  we briefly recall
a few basic definitions.

\subsection{Quantum Hamiltonian reduction.}
Throughout, we work with associative unital
$\C$-algebras. We write $\Hom=\Hom_\C,\,\ot=\ot_\C,$ etc.

Let
 $A$ be an associative
algebra, that may also be viewed as a Lie algebra
with respect to the commutator
Lie bracket. Given a  Lie algebra
 $\g$  and a Lie algebra  homomorphism
 $\rho: \g\to A$,
one has an {\em adjoint} $\g$-action on $A$ given by
$\ad x: a\mapsto \rho(x)\cdot a - a\cdot \rho(x),\, x\in\g,$
$a\in A.$
The left ideal $A\cdot \rho(\g)$ is stable under the adjoint action.
Furthermore, one shows that multiplication in $A$ induces
a well defined associative algebra structure on
$$
\A(A,\g,\rho):=\bigl(A/A\cd\rho(\g)\bigr)^{\ad\g},
$$
 the space
of $\ad\g$-invariants in $ A/ A\cd \rho(\g)$.
The resulting algebra $\A(A,\g,\rho)$ is called
 the {\em quantum  Hamiltonian reduction} of $A$ at $\rho$.

Observe  that, if $a\in A$ is such that the element
$a\,\text{mod}\, A\cd\rho(\g) \,\in A/ A\cd \rho(\g)$
is $\ad\g$-invariant, then
the operator of right multiplication by $a$
descends to a well-defined map
$R_a:  A/ A\cdot \rho(\g)\to A/ A\cdot \rho(\g)$.
Moreover, the assignment $a\mapsto R_a$
induces an algebra isomorphism
$\dis \A(A,\g,\rho)=
( A/ A\cd \rho(\g)\bigr)^{\ad\g}\iso \bigl(\End_ A( A/ A\cd \rho(\g))\bigr)^{\op}.
$

If $A$, viewed as an $\ad\g$-module, is semisimple, i.e.,
splits into a
(possibly infinite) direct sum of
irreducible finite dimensional $\g$-representations, then
the operations of taking $\g$-invariants and taking the quotient
commute, and we may write
\beq{BA}
\A(A,\g,\rho)=\bigl(A/A\cd\rho(\g)\bigr)^{\ad\g}=
A^{\ad\g}\big/(A\cd\rho(\g))^{\ad\g}.
\eeq
Observe that, in this formula, $(A\cdot\rho(\g))^{\ad\g}$ is
a {\em two-sided} ideal of the algebra $A^{\ad\g}$.

Any $A$-module $M$ may be viewed also as a $\g$-module, via the 
homomorphism $\rho$, and
 we write
$M^\g:=\{m\in M\mid \rho(x)m=0,\;\forall x\in \g\}$
for the corresponding space of $\g$-invariants.
Let $\Lmod{(A,\g)}$ be the full subcategory of the abelian
 category of  left $A$-modules whose objects
are semisimple as
 $\g$-modules. Let $\Lmod{\A(A,\g,\rho)}$ be the abelian category
of left  $\A(A,\g,\rho)$-modules.

One defines an exact functor,
called {\em Hamiltonian reduction functor},
as follows
\beq{BH}
\BH:\ \Lmod{(A,\g)}\to\Lmod{\A(A,\g,\rho)},\quad
M\longmapsto \BH(M):=\Hom_{A}\bigl(A/ A\cd \rho(\g), M\bigr)=
M^{\g},
\eeq
where the  action of $\A(A,\g,\rho)$ on $\BH(M)$
comes from the tautological
{\em right} action
of $\End_ A( A/ A\cd \rho(\g))$ on
$A/ A\cd \rho(\g)$ and the above mentioned
isomorphism $\dis \A(A,\g,\rho)=\bigl(\End_ A( A/ A\cd \rho(\g))\bigr)^{\op}.
$ 

\subsection{Symplectic reflection algebras for wreath-products.}\label{sra}
Let $n$ be a positive integer. Let $S_n$  be the permutation
group of $[1,n] := \{1,\ldots, n\}$, and write $s_{\ell m}\in S_n$
for the transposition $\ell\leftrightarrow m$.
Let $L$ be a 2-dimensional complex vector space,
and $\om $ a symplectic form on $L$.

Let $\Gamma$ be a finite subgroup of $\Sp(L)$, and
let $\GG:=S_n\ltimes\Gamma^n$ be a wreath product group
acting naturally in $L^n$.
Given $\ell\in [1,n]$ and $\ga\in\Gamma$, resp. $v\in L$,
we will write
$\ga_{(\ell)}\in\GG$ for $\ga$ placed in the $\ell$-th factor $\Gamma$,
resp. $v_{(\ell)}\in L^n$ for $v$ placed in the $\ell$-th factor $L.$

According to  \cite{EG}, there is a family of associative algebras,
 called  symplectic reflection
algebras, attached to the pair $(L^n,\GG)$ as above.
To define these algebras, write $Z\Gamma$ for the center of the group
algebra $\C[\Gamma]$ and let $Z_o\Gamma\sset Z\Gamma$ be a codimension
1 subspace formed by the elements
\beq{c}
c = \sum\nolimits_{\ga \in \Gamma\smallsetminus
\{1\}}\, c_{\ga}\cdot\ga \in Z\Gamma, \quad\forall
c_\ga\in\C.
\eeq

Given $t,k\in\C$ and $c\in Z_o\Gamma$,  the
corresponding \emph{symplectic reflection algebra} $\hh_{t,k,c}(\GG)$,
with parameters $t,k,c$, may be defined, cf. \cite[Lemma 3.1.1]{GG},
as a quotient of the smash product algebra
$T(L^n) \rtimes \C[\GG]$ by the following relations:
\begin{align}
[x_{(\ell)}, y_{(\ell)}]
=& t\cdot 1+ \frac{k}{2} \sum_{m\neq \ell}\sum_{\ga\in\Gamma}
s_{\ell m}\ga_{(\ell)}\ga_{(m)}^{-1}
+ \sum_{\ga\in\Gamma\smallsetminus\{1\}} c_{\ga}\ga_{(\ell)} ,
\qquad \forall  \ell\in[1,n];  \label{relation1}\\
[u_{(\ell)},v_{(m)}]=& -\frac{k}{2} \sum_{\ga\in\Gamma} \omega(\ga u,v)
s_{\ell m}\ga_{(\ell)}\ga_{(m)}^{-1},
\qquad \forall  u,v\in L,\ \ell,m\in [1,n],\ \ell\neq m,\label{relation2}
\end{align}
where  $\{x,y\}$ is a fixed basis  for $L$ with $\om (x,y)=1$.

\subsection{Quivers.}
Let $Q$ be  an extended Dynkin quiver with  vertex set $I$,
and let $o\in I$ be  an extending vertex of $Q$.

\begin{definition}\label{CMQ} The quiver  $\Qc$  obtained from $Q$ by adjoining an
additional vertex $s$
and an arrow $b: s\to o$ is called
 the \emph{Calogero-Moser quiver} for $Q$.
Thus, $\Ic=I\sqcup\{s\}$ is the vertex set for $\Qc$, and
the vertex $s$ is called the \emph{special} vertex.
\end{definition}

Given $\al=\{\al_i\}_{i\in \Ic}\in \ZZ^{\Ic}$,  a dimension vector for $\Qc$,
write
\beq{rep_CM}
 \Rep_\al(\Qc) := \bigoplus_{\{a:i\to j\;\mid\; a\in \Qc\}}
\Hom(\C^{\al_i},\C^{\al_j})=\bigoplus_{\{a:i\to j\;\mid\; a\in \Qc\}}
\mathrm{Mat}(\al_j \times \al_i, \C)
\eeq
for
the  space of representations of  $\Qc$ of dimension $\al.$
Let $\D(\Qc,\al)$ be the algebra of polynomial
differential operators on the vector space
$\Rep_\al(\Qc)$.

The group $\GL(\al):=\prod_{i\in \Ic}\GL(\C^{\al_i})$ acts naturally on
 $\Rep_\al(\Qc),$
by conjugation.
Hence, each element $h$
of the Lie algebra  $\gl(\al) :=\Lie \GL(\al)$ gives rise to a vector
field $\bxi_h$ on $\Rep_\al(\Qc)$. This yields a Lie algebra map
$\bxi: \gl(\al)\to \D(\Qc,\al)$.

The center of the reductive Lie algebra
 $\gl(\al)=\oplus_{i\in I}\gl(\al_i)$
is clearly isomorphic to $\C^I$.
Therefore, associated with any $\chi=\{\chi_i\}_{i\in I}\in\C^I$,
one has a  Lie algebra homomorphism ${\chi:\gl(\al)\to\C,}$
$x=\oplus_{i\in I} x_i\, \mto \, \sum_{i\in I} \chi_i\cdot\Tr^{} x_i.$
We will use additive notation for such homomorphisms
and write $\bxi -\chi: \gl(\al)\to \D(\Qc,\al)$
(rather than $\bxi\ot (-\chi)$)
for the Lie algebra map
$h\mto\bxi_h-\chi(h)\cd 1_\D$.
Let $\im(\bxi- \chi)$  denote the image of the latter map.

We may apply  Hamiltonian reduction 
 \eqref{BA} to the algebra
$\D(\Qc,\al)$ and to the Lie algebra map
$\bxi -\chi$. This way, we get the  algebra
\begin{align}\label{Achi}
 \A\big(\D(\Qc,\al),\,\gl(\al),\,\bxi- \chi\big)&\,=
\D(\Qc,\al)^{\GL(\al)}/\JJ_\chi,\\
\text{where}\quad
\JJ_\chi&:=
\big( \D(\Qc,\al)\cd
\im(\bxi-\chi)\big)^{\GL(\al)}.\nonumber
\end{align}

Let $T^*\Rep_\al(\Qc)$ be the cotangent bundle on $\Rep_\al(\Qc).$
The  total space of the
cotangent bundle comes equipped with the canonical symplectic
structure and with a moment map
\beq{moment}
\bmu:\ T^*\Rep_\al(\Qc)\too \gl(\al)^*\cong\gl(\al).
\eeq
We may apply the {\em classical} Hamiltonian reduction 
to  $\C[T^*\Rep_\al(\Qc)],$ the Poisson algebra of polynomial functions
on $T^*\Rep_\al(\Qc)$.
This way, we get the Poisson algebra
$\C[\bmu\inv(0)]^{\GL(\al)}$
of $\GL(\al)$-invariant polynomial functions on the zero fiber
of the moment map.
The  algebra in \eqref{Achi}
 may be viewed as a quantization of the Poisson algebra
 $\C[\bmu\inv(0)]^{\GL(\al)}$.

\subsection{Main result.}\label{main_sec}
From now on, we fix  $n\in{\mathbb{N}}$, a 2-dimensional
symplectic vector space $L$ and  $\G\sset Sp(L)$,
a  finite subgroup as in \S\ref{sra}.
To $(n,L,\G),$ we will associate a quiver $Q$, a dimension vector $\alpha$,
 and a character $\chi$ as follows.

We let $Q$ be an affine Dynkin quiver associated to
$\Gamma$ via the McKay correspondence.
Thus, the set $I$ of vertices of $Q$ is identified with the
set of isomorphism classes of irreducible representations
of $\Gamma$. Let $N_i$ be the irreducible representation of
$\Gamma$ corresponding to the vertex $i\in I$, and
let $\delta_i = \dim N_i$.
The extending vertex $o\in I$ corresponds to the
trivial representation of $\Gamma$, so $\de_o=1.$
The vector $\de=\{\de_i\}_{i\in I}\in\ZZ^I$ is
 the minimal positive imaginary
root of the  affine root system associated to $Q$.
Motivated by M. Holland \cite{Ho}, we put
\beq{data}
\pa=\{\pa_i\}_{i\in I}\in\ZZ^I,\quad
\pa_i:=n\bigl(-\de_i + \sum\nolimits_{\{a\in \Qc\;\mid\; t(a)=i\}}
\de_{h(a)}\bigr),
\quad\forall i\in I.
\eeq

Given
 a central element  $c\in Z\G$,
 write $\Tr(c; N_i)$ for the trace
of $c$ in the simple $\G$-module $N_i,\, i\in I$. Thus,
for any $c\in Z_o\G,$ see  \eqref{c},
we have $\sum_{i\in I}\de_i\cdot\Tr(c; N_i)=0$.
Associated with any data $n\in{\mathbb{N}},k\in\C,$ and $c\in Z_o\G$,
 we introduce 
 three vectors
$$
 \chi=\{\chi_i\}_{i\in \Ic},\en
  \chi'=\{\chi'_i\}_{i\in \Ic}\in\C^{\Ic},\en
\text{and}\en \la(c)=\{\la(c)_i\}_{i\in I}\in\C^I,\en
\text{such that}\en\de\cdot\la(c)=1,
$$
where we have used standard
notation $\de\cdot\la=\sum_i \de_i\cdot\la_i.$
These vectors are defined as follows
\begin{align}\label{vectors}
&\la(c)_i:=\Tr(c; N_i)+\de_i/|\G|,\quad\forall i\in I;\nonumber\\
\chi_s&:=n
 (k\cd|\G|/2-1)+1,\quad\chi_o=\la(c)_o-\pa_o-k\cd|\G|/2,\quad
\chi_i=\la(c)_i-\pa_i,\quad\forall i\in I\sminus\{o\};\nonumber\\
\chi'_s&:=\chi_s-1=n(k\cd|\G|/2-1),\quad
\chi'_i=\chi_i,\quad\forall i\in I.
\end{align}

We are going to consider representations of the quiver
$\Qc$ with dimention vector
\beq{al}
\al=\{\al_i\}_{i\in \Ic}\in\ZZ^{\Ic}_{\geq 0},\quad
\text{where}\quad\al_s:=1,
\en\text{and}\en\al_i:=n\cd\de_i,\quad\forall i\in I.
\eeq

Let  $\chi'\in\C^\Ic$ be as in   \eqref{vectors}, and let
$\JJ_{\chi'}=\big(\D(\Qc,\al)\cdot
\im(\bxi-\chi')\big)^{\GL(\al)}$
be the  corresponding  two-sided ideal in 
$\D(\Qc,\al)$, cf. \eqref{Achi}. 
Write $\e:=\frac{1}{|\GG|}\sum_{g\in\GG} g$ for the `symmetrizer'
idempotent
viewed as an element of the symplectic reflection algebra  $\hh_{t,k,c}(\GG)$.

We are now in a position to state our main result
about  deformed Harish-Chandra
homomorphisms for symplectic reflection algebras
associated with a wreath-product.
According to \cite{EG},
the importance of
the  deformed Harish-Chandra
homomorphism is due to the fact that this homomorphism
 provides a description of the  {\em spherical subalgebra}
$\dis\e\hh_{t,k,c}(\GG)\e\sset\hh_{t,k,c}(\GG)$
 in terms of quantum Hamiltonian reduction of the ring
of polynomial differential operators on the vector space
$\Rep_\al(\Qc)$.
In the special case of a {\em cyclic} group $\G\sset SL_2(\C)$,
that is, for
quivers $Q$ of type $\widetilde {\mathbf{A}}_m$ (equipped
with the cyclic orientation),
the deformed Harish-Chandra
homomorphism has been already constructed in \cite{Ob}, see also
 \cite{Go}. 
In all other cases, a construction of
the  deformed Harish-Chandra
homomorphism $\Phi_{k,c}$ will be given in the present paper.

Our  main result reads

\begin{theorem}\label{main} Assume that $\G\sset SL_2(\C)$ is not 
a cyclic group of {\em odd} order (i.e. $Q$ is not
of type $\widetilde {\mathbf{A}}_{2m}$), and put
$t:=1/|\G|$.  Then, for any $n\in{\mathbb{N}},k\in\C,c\in Z_o\G$,
there is an algebra isomorphism
$$\Phi_{k,c}:\
\A\big(\D(\Qc,\al),\gl(\al),\bxi-
\chi'\big)=\D(\Qc,\al)^{\GL(\al)}/\JJ_{\chi'}
\,\iso \,
\e\hh_{t,k,c}(\GG)\e.$$

Furthermore, the map $\Phi_{k,c}$ is compatible with natural increasing
filtrations
on the algebras involved and the corresponding associated graded map 
gives rise to a graded Poisson
algebra isomorphism, cf \eqref{moment}:
$$\gr\Phi_{k,c}:\
\C[\bmu\inv(0)]^{\GL(\al)}\iso \gr\bigl(\e\hh_{t,k,c}(\GG)\e\bigr).
$$
\end{theorem}

This theorem is a slightly modified and corrected version
of \cite[Conjecture 11.22]{EG}
(in \cite{EG}, as well as in the main body of the
present paper, everything is stated in terms of
the quiver $Q$ rather than in terms of
the Calogero-Moser quiver $\Qc$, see Definition \ref{hc_def} and
Theorem \ref{maintheorem} in \S\ref{hc_sec} below;
however, the two approaches are easily seen to be
equivalent).
 Theorem \ref{main} is a common generalization of
two earlier results. The first one is \cite[Theorem 6.2.3]{GG2}, cf. also
\cite[Corollary 7.4]{EG};
it corresponds to the (somewhat degenerate)
case of $\G=\{1\}$. The second result, due to
M.~Holland \cite{Ho}, is a special case of Theorem \ref{main}
for $n=1$, where the symplectic reflection algebra
is Morita equivalent to a deformed preprojective algebra of \cite{CBH}.
Also, in the special case of a cyclic group $\G=\ZZ/m\ZZ$ the isomorphism
of  Theorem \ref{main} has been recently constructed
in \cite{Go} using the results from \cite{Ob}.

A `classical' counterpart of  Theorem \ref{main} involving
classical Hamiltonian reduction (at {\em
generic} values of the moment map \eqref{moment})  has been proved
in \cite[Theorem 11.16]{EG}.

Combining Theorem \ref{main} with \eqref{BH}, and using  
the same argument as in the proof of 
\cite[Proposition 6.8.1]{GG2}, we deduce
\begin{corollary}\label{BHcor}
There exists an exact functor of Hamiltonian reduction
$$\BH: \Lmod{(\D(\Qc,\al),\gl(\al))}\too
\Lmod{{\e\hh_{t,k,c}(\GG)\e}}.
$$
This functor induces an equivalence
$\Lmod{(\D(\Qc,\al),\gl(\al))}/\Ker\BH
\iso\Lmod{{\e\hh_{t,k,c}(\GG)\e}}.$\qed
\end{corollary}

We expect that the Hamiltonian reduction functor
induces an equivalence between the subcategory of 
$\Lmod{(\D(\Qc,\al),\gl(\al))}$ formed by $\D$-modules
whose characteristic variety is contained in the
{\em Nilpotent Lagrangian}, see \cite[\S12]{Lu1},
and the category of finite dimensional $\e\hh_{t,k,c}(\GG)\e$-modules.

\subsection{Four homomorphisms.}\label{strat}
Our construction of the isomorphism $\Phi_{k,c}$ in
Theorem \ref{main} is rather indirect. It involves
four additional algebras and four homomorphisms between those algebras,
which are important in their own right.

The first
algebra, to be denoted $\Pi'(\Qc)$, is a slightly renormalized
version of the deformed preprojective algebra, with
appropriate parameters, cf. \cite{CBH}, associated to the Calogero-Moser
quiver $\Qc$.
The second algebra, to be denoted $\M$,
contains the spherical algebra $\e\hh_{t,k,c}(\GG)\e$ as a subalgebra.
The algebra $\M$ is a `Calogero-Moser cousin' of  {\em generalized
preprojective algebras}
introduced by two of us in \cite[(1.2.3)]{GG},  see
also Definition
\ref{dea} below.

 The third algebra,  $\T_\chi,$ is a
`matrix-valued' counterpart
of the algebra introduced in \eqref{Achi}.
To define this algebra, we
introduce the following vector spaces
\beq{N}
\N=\oplus_{i\in\Ic} \N_i,\en\text{where}\en
\N_s:=N^*_o \cong \C,\en\text{and}\en
\N_i:=N_i^*\ot\C^n,\en
\forall i\in I.
\eeq
Thus, we have $\N_i\cong\C^{\al_i},$
so the group $\GL(\al)$ acts  on $\N$ in an obvious way, and this gives
the tautological
representation $\tau: \gl(\al)\to \End \N$.
Following  M. Holland
\cite{Ho}, we apply the
 quantum Hamiltonian reduction to the algebra
$\D(\Qc,\al) \ot \End\N$ and to the Lie algebra
homomorphism
$$\bxi- (\chi-\tau): \gl(\al) \to\D(\Qc,\al) \ot \End\N,\quad
h\mto\bxi_h\ot\id_\N-1_\D\ot \bigl(\chi(h)\id_\N-\tau(h)\bigr),$$
where $\chi: \gl(\al) \to\C$ is as
in \eqref{vectors}. This way, we get an algebra
\beq{Tchi}
\T_\chi :=
\frac{\big(\D(\Qc,\al) \ot
\End\N \big)^{\GL(\al)}}{\Big( \big( \D(\Qc,\al) \ot \End\N \big)\cd
\im(\bxi-(\chi-\tau)) \Big)^{\GL(\al)}}
\eeq

Now, let ${\mathbb P}^1=(L\sminus\{0\})/\C^\times$
be  the projective line.
We will consider  an appropriate $\GG$-equivariant vector bundle
of rank $\dim\N$ on $\X$, where
$\X\sset ({\mathbb P}^1)^n$ is
 a $\GG$-stable Zariski open dense
 subset
in the cartesian product of $n$ copies of ${\mathbb P}^1$.
Further, we will define  a certain
algebra
$\D(\X,p,\varrho)$  of {\em twisted}
 differential operators acting
in that vector bundle, see
\S\ref{tdo} for the notation and also
\eqref{thrad}.

One has the following diagram of four algebra homomorphisms, 
all denoted by various $\Th$'s,
involving the four algebras introduced above

\beq{cccc}
\xymatrix{
&&\Pi'(\Qc)\ar[dll]_<>(0.5){\Th^\text{Holland}}
\ar[drr]^<>(0.5){\Th^\text{Quiver}}&&\\
{\T_\chi\en}
\ar[drr]_<>(0.5){\Th^\text{Radial}}
&&&&
{\en\M}
\ar[dll]^<>(0.5){\Theta^\text{Dunkl}}\\
&&\D^{}(\X,p,\varrho)^\GG &&
}
\eeq

In this diagram,  the map $\Th^\text{Holland}$ is (a slightly renormalized
version of) an algebra homomorphism introduced by M. Holland in \cite{Ho}.
 The map
$\Theta^\text{Dunkl}$ is a $\G$-analog of the
Dunkl representation for rational Cherednik algebras,
cf. \cite{EG}.
The map $\Th^\text{Radial}$ is obtained by a `radial part' type
construction
with respect to an appropriate transverse slice
to generic $\GL(\al)$-orbits in $\Rep_{\al}(\Qc)$.
We produce such a slice using a
map
$L^{\oplus n}\to \Rep_{\al}(\Qc)$,
which is  generically injective and is such
that its image is  generically transverse to
 $\GL(\al)$-orbits in $\Rep_{\al}(\Qc)$.
Our radial
part construction associates
to a  $\GL(\al)$-invariant  differential
operator $u\in\big(\D(\Qc,\al) \ot
\End\N \big)^{\GL(\al)}$
a
$\GG$-invariant twisted differential operator
 $\Th^\text{Radial}(u)\in \D(\X,p,\varrho)^\GG.$

The fourth  map, $\Th^\text{Quiver}$, is new.
The main idea behind the construction of this map,
as well as the definition of the algebra $\M$,
will be outlined in \S\ref{pp}
 below and a more rigorous  treatment will be given later,
in \S\ref{quivermap}.

\begin{remark}  In the special case of a cyclic group
$\G=\ZZ/m\ZZ$, the Dunkl operators that we consider
 are {\em not} the same as those introduced
earlier by Dunkl-Opdam in \cite{DO}.
\end{remark}
\subsection{Strategy of the proof of Theorem \ref{main}.}\label{strategy_pf}
The proof of the main theorem is based
on the following
 key result

\begin{theorem}\label{bulk} Diagram \eqref{cccc} commutes, i.e.,
we have:
$\dis\;\Th^{\opp{Radial}}\ccirc\Th^{\opp{Holland}}=
\Th^{\opp{Dunkl}}\ccirc\Th^{\opp{Quiver}}.
$
\end{theorem}

The proof of this Theorem is long and messy;
it occupies about one half of the paper.
In the proof, we explicitly compute
both sides of the equation
$\Th^{\text{Radial}}\ccirc\Th^{\text{Holland}}(x)=
\Th^{\text{Dunkl}}\ccirc\Th^{\text{Quiver}}(x),$
for an appropriate set $\{x,\,x\in\Pi'(\Qc)\}$ of generators
 of the algebra~$\Pi'(\Qc)$.

To deduce Theorem \ref{main} from Theorem
\ref{bulk}, one has to be able to replace in diagram \ref{cccc}
the algebra $\T_\chi$, of `matrix valued' twisted differential operators,
by  a `smaller' algebra of   {\em scalar-valued}
twisted differential operators of
the form $\A\big(\D(\Qc,\al),\,\gl(\al),\,\bxi-
\chi\big)$,  that appears in
 Theorem~\ref{main}.

To this end,  let   $ {\mathsf p}_s\in\End\N$
denote the idempotent
 corresponding to the projection $\N=\bigoplus_{j\in\Ic}
\N_j\onto \N_s.$ For $\chi,\chi'$  as in \eqref{vectors},
one proves
\begin{equation} \label{achi}
 {\mathsf p}_s\T_{\chi} {\mathsf p}_s\cong\D(\Qc,\al)^{\GL(\al)}/\JJ_{\chi'}=
\A\big(\D(\Qc,\al),\,\gl(\al),\,\bxi-
\chi'\big)=:\A_{\chi'}.
\end{equation}

Write $e_i$  for
the idempotent in the algebra $\Pi'(\Qc)$
 corresponding to the trivial
path at $i$.
It is easy to see that the map $\Th^{\text{Quiver}}$
sends the subalgebra $e_s\Pi'(\Qc)e_s\sset\Pi'(\Qc)$,
spanned by paths beginning and ending at the special vertex $s$, into
 $\e\hh_{t,k,c}(\GG)\e,$ a  subalgebra in $\M$.
Furthermore,
restricting diagram \eqref{cccc} to the subalgebra
$e_s\Pi'(\Qc)e_s$, one obtains   four
 algebra homomorphisms along the perimeter of the
following diagram
\beq{ccc}
\xymatrix{
&&e_s\Pi'(\Qc)e_s\ar@{->>}[dll]_<>(0.5){\Th^\text{Holland}}
\ar[drr]^<>(0.5){\Th^\text{Quiver}}&&\\
{{{}^{^{}}}\D(\Qc,\al)^{\GL(\al)^{\!\!{}^{}}}/^{\!}\JJ_{\chi'}{}^{^{}}
\en}\ar@{.>}[rrrr]^<>(0.5){\Phi_{k,c}}
\ar[drr]_<>(0.5){\Th^\text{Radial}}
&&&&
{\en\e\hh_{t,k,c}\e}
\ar@{^{(}->}[dll]^<>(0.5){\Theta^\text{Dunkl}}\\
&&\D^{}(\X,p,\varrho_s)^\GG&&
}
\eeq
Here, $\D(\X,p,\varrho_s)^\GG$ stands for
an appropriate ring of {\em scalar-valued}
$\GG$-invariant twisted differential operators on $\X$.

The perimeter of  diagram \eqref{ccc} commutes
by  Theorem \ref{bulk}. In addition, one proves
\begin{lemma} In diagram \eqref{ccc},
the map $\Th^{\text{Holland}}$ is {\em surjective}
and the map $\Th^{\text{Dunkl}}$
is {\em injective}.
\end{lemma}
It is clear that the Lemma yields
$$
\Ker\Th^{\text{Holland}}\sset
\Ker(\Th ^{\text{Radial}}\ccirc\Th^{\text{Holland}})=
\Ker(\Th^{\text{Dunkl}}\ccirc\Th^{\text{Quiver}})
=\Ker\Th^{\text{Quiver}}.
$$
The resulting inclusion
$
\Ker\Th^{\text{Holland}}\sset\Ker\Th^{\text{Quiver}}$ implies
that we  may (and will) define the dashed arrow $\Phi_{k,c}$
in diagram \eqref{ccc}
to be the composite

$$
{\Large\xymatrix{
{\frac{\D(\Qc,\al)^{\GL(\al)}}{\JJ_{\chi'}}\;}
\ar@{=}[rr]^<>(0.5){(\Th^{\text{Holland}})\inv}&&
{\;\frac{e_s\Pi'(\Qc)e_s}{\Ker\Th^{\text{Holland}}}\;}
\ar@{->>}[r]^<>(0.5){_\text{proj}}&
{\;\frac{e_s\Pi'(\Qc)e_s}{\Ker\Th^{\text{Quiver}}}\;}
\ar@{^{(}->}[rr]^<>(0.5){\Th^{\text{Quiver}}}
&&}}
\e\hh_{t,k,c}\e.
$$

To complete the proof of Theorem \ref{main}, one observes that
all the objects appearing in diagram \eqref{ccc} come equipped with
natural filtrations, and all the maps in the diagram are
filtration preserving. Therefore, to prove that the map
$\Phi_{k,c}$ is bijective, it suffices to  show a similar statement for
$\gr\Phi_{k,c}$, the associated graded map. The latter
statement follows readily from the results of \cite{CB} and \cite{GG2}
concerning  the geometry
of moment maps arising from representations of affine Dynkin quivers.

\subsection{The algebra $\M$ and the map $\Th^{\text{Quiver}}$.}\label{pp}
To define the algebra $\M$ that appears in diagram
\eqref{Tchi}, we will first introduce in \eqref{e}
certain idempotents $e_{i,n-1}\in \C[\GG],\,i\in I$.
Then, we let
\beq{M}
\bbM := \hh_{t,k,c}(\GG) \e \bigoplus \big( \oplus_{i\in I}
\hh_{t,k,c}(\GG) e_{i,n-1} \big).
\eeq
Thus, $\bbM$ is a left
$\hh_{t,k,c}(\GG)$-module, and we put
 $\M:=(\End_{\hh_{t,k,c}(\GG)} \bbM)^{\op}.$
This endomorphism algebra
 is built out of $\Hom$-spaces
 between  various $\hh_{t,k,c}(\GG)$-modules
which appear as
direct summands in \eqref{M}. The  $\Hom$-spaces are easily computed,
and we find
\begin{align}\label{MM}
&\M=\bigoplus_{i,j\in\Ic}\M_{i,j},\quad\text{where}\quad
\M_{s,s}
= \e\hh\e,\quad\text{and}\\
&\M_{s,j}=\e \hh e_{j,n-1} ,\quad
\M_{i,s}
=  e_{i,n-1}\hh \e,\quad
\M_{i,j}=e_{i,n-1}\hh e_{j,n-1},\quad\forall i,j\in
I.\nonumber
\end{align}

Each direct summand $\M_{i,j}$ here is a subspace
of the algebra $\hh_{t,k,c}(\GG)$, and  multiplication in the algebra
$\M$ is given by   `matrix multiplication'
$\M_{i,j}\times \M_{j,k}\to\M_{i,k}$ where, for each $i,j,k\in\Ic,$
the corresponding
pairing is induced by  the multiplication in  $\hh_{t,k,c}(\GG)$.

Our construction of
the map $\Th^{\text{Quiver}}$ is based on  an  exact functor
\beq{funct}
\Mod{\hh_{t,k,c}(\GG)}\too \Mod{\Pi(\Qc)},\quad
M\mto \wt M.
\eeq

To define this functor, let  $L\bis$, resp. $\G\bis$,
 be a copy (inside the algebra $\hh_{t,k,c}(\GG)$)
of our 2-dimensional vector space $L$,
resp. copy of the group $\G$, corresponding
to the first direct summand in $L^{\oplus n}$.
Further, let
 $S_{n-1}$ be the subgroup of $S_n$
which permutes $[2,n]$, and
let $\GGG=S_{n-1}\rtimes \G^{n-1}\sset
\GG$
be the wreath-product subgroup corresponding to the
last $n-1$ factors in $\G^n$.
It is clear from
the commutation relations in $T(L^{\oplus n})\rtimes \C[\GG]$
that
any element of the subalgebra  $\ttt\sset\hh_{t,k,c}(\GG),$
generated by $L\bis$ and $\G\bis$,
commutes with $\GGG$.

Now, let $M$ be  an arbitrary left $\hh_{t,k,c}(\GG)$-module.
We deduce that
the space $M^{\GGG}\sset M$, of $\GGG$-invariants, is stable
under the action of the subalgebra $\ttt$.
Thus,
to each vertex $i\in Q$ we may attach
the vector space
$\dis M_i:=\Hom_{\G\bis}(N_i, M^{\GGG})$,
the corresponding $\G\bis$-isotypic component.
Further, following the strategy of \cite{CBH} and using the McKay correspondence,
we see that the action map
$L\bis\ot M^{\GGG}\to M^{\GGG}$ induces
linear maps between various isotypic components $M_i$.
This way,
the collection $\{M_i\}_{i\in I}$ acquires
the structure of
a representation of the quiver $\dq$.
In addition, the subspace $M_s:=M^\GG\sset M$
is clearly contained in $M_o=\Hom_{\G\bis}(N_o, M^\GGG)=M^\GGG$
as a {\em canonical} direct summand. Therefore
the imbedding $b: M_s\to M_o$ and the
projection $b^*: M_o\to M_s$ provide additional
maps, making the collection  $\{M_i\}_{i\in \Ic}$
a  representation of the quiver $\QQc$.
One can check that this representation
descends to a representation of
the algebra $\Pi(\Qc)$, which is a quotient
of the path algebra of  $\QQc$.
Thus,  to any  $\hh_{t,k,c}(\GG)$-module
$M$ we have assigned a
$\Pi(\Qc)$-module $\wt M=\oplus_{i\in \Ic}M_i$.
This gives the desired functor \eqref{funct}, cf. \S\ref{uu} below for a generalization.

Finally, we
apply the functor $M\mto \wt M$ to $M:=\bbM$,
the  $\hh_{t,k,c}(\GG)$-module in \eqref{M}.
It is immediate from \eqref{MM}
that one has  a natural bijection
$\M\cong{\widetilde{{\mathbb M}}}$. The bijection
gives $\M$ the structure of a left $\Pi(\Qc)$-module,
 moreover, the action of  $\Pi(\Qc)$ on $\M$
commutes with right multiplication (with respect to the algebra
structure)
by the elements of $\M$. It follows that the  $\Pi(\Qc)$-module
 structure on $\M$ comes, via left multiplication,
from  an algebra homomorphism
$\Pi(\Qc) \to \M.$
The latter homomorphism clearly restricts to
a  homomorphism
$e_s \Pi(\Qc) e_s \to \M_{s,s}=\e\hh_{t,k,c}(\GG)\e.
$

There is a modification of the above construction,
to be explained in \S\ref{quivermap}, in which
the algebra $\Pi(\Qc)$ is replaced
by the renormalized algebra $\Pi'(\Qc)$.
This way, one obtains similar algebra homomorphisms
\beq{Thquiv}
\Theta^{\text{Quiver}}:\Pi'(\Qc) \to \M,
\quad\text{and}\quad \Theta^{\text{Quiver}}:
e_s \Pi'(\Qc) e_s \to  \M_{s,s}=\e\hh_{t,k,c}(\GG)\e.
\eeq

\subsection{Applications to Reflection functors and Shift functors.}
\label{uu} 
In
\S\ref{refl_sec}, we study reflection functors and shift
functors for  generalized preprojective
algebras and symplectic reflection algebras
associated with  wreath-products, cf.~\cite{GG}.

More  generally,
let $Q$ be an arbitrary (not necessarily extended Dynkin) quiver,
with vertex set $I$. Write  $\mathfrak C=(\mathfrak C_{ij})$ for the
generalized Cartan matrix of $Q$ and $W$ for
the Weyl group $W$, defined as the group generated by the
simple reflections $r_i$ for $i\in I$. The group $W$ acts on $\C^I$
as $r_i: \la=\sum_{j\in I} \la_je_j\mto
 \la-\sum_{j\in I} \mathfrak C_{ij}\la_ie_j$.

For any  $\la\in\C^I$,
one has an algebra
 $\Pi'_\la(Q),$ a renormalized version
of the corresponding {\em deformed preprojective
algebra}  studied in \cite{CBH}.
Further,
for any integer $n\geq 1$, and complex parameters $\nu\in\C$
and $\la\in\C^I$, we have associated in  \cite[(1.2.3)]{GG}, see
also Definition
\ref{dea} below,
a  {\em generalized preprojective
algebra} $\sfA_{n,\la,\nu}(Q).$

For each $i\in I$, there are reflection functors $F'_i$
for the corresponding
deformed preprojective algebras $\Pi'_\la(Q)$, introduced in \cite{CBH},
and also their analogues for generalized  preprojective algebras,
introduced  in \cite{Ga}:
\beq{gan}
 F_i:\
 \Mod{\mathsf A_{n,\la,\nu}(Q)}\too
\Mod{\mathsf A_{n,r_i(\la),\nu}(Q)}.
\eeq

We will show in \S\ref{high} that these functors satisfy standard
Coxeter relations:

\begin{proposition} \label{weylrel} 
For the  reflection functors $F_i$ 
for generalized  preprojective algebras, one has:

\vi If $\la_i\pm p\nu\neq 0$ for $p=0,1,...,n-1$, then
$F_i^2=\id$.

\vii Suppose $\mathfrak{C}_{ij}=0$.
If $\la_i\pm p\nu\neq 0$ and $\la_j\pm p\nu\neq 0$
for $p=0,1,...,n-1$, then $F_iF_j=F_jF_i$.

\viii Suppose $\mathfrak{C}_{ij}=-1$.
If $\la_i\pm p\nu\neq 0$, $\la_j\pm p\nu\neq 0$ and
$\la_i+\la_j\pm p\nu\neq 0$ for $p=0,1,...,n-1$, then
$F_iF_jF_i = F_jF_iF_j$.
\end{proposition}
Part \vi of the Proposition has been already proved
in \cite[Theorem 5.1]{Ga};
parts \vii and \viii are new. 
In the special
case $n=1$, the Proposition  is due to 
\cite{CBH},
\cite{Na}, 
\cite{Lu2}, and \cite{Maf}.
However, we believe that, even in that
special case,  our proof appears to be simpler.

Next, given $c$ as  in \eqref{c}, we put
$$
c':=\sum_{\ga\in\G\sminus\{1\}}\, (2t-c_\ga)\cd\ga\inv
\quad\text{and}\quad
 \e_- := \frac{1}{n!}\sum_{\sigma\in S_n} (-1)^\sigma \sigma
(e_o\otot e_o).
$$
Using  our main   Theorem \ref{main} 
and reflection functors,  we will deduce
\begin{corollary}  \label{maincorollary} 
For $t=1/|\G|$ and any $c$ as in \eqref{c},
there are algebra isomorphisms
$$\e\hh_{t,k,c}\e  \simeq \e_-\hh_{t, k-2t,c'}
\e_-\simeq \e_-\hh_{t,k-2t,c}\e_-.
$$
\end{corollary}
We will prove the first isomorphism above in \S\ref{fisc}
and the second in \S\ref{shiftsubsec}.
Using the composite isomorphism in  Corollary \ref{maincorollary}, we define the
{\it shift functor} to be the functor
\beq{SS}
 \mathbb S : \Mod{\hh_{t,k,c}} \too \Mod{\hh_{t,k-2t, c}},
\quad V \mapsto \hh_{t,k-2t,c}\e_-\ot_{\e\hh_{t,k,c}\e} \e V. 
\eeq

Finally, we can extend the construction
exploited in the definition of the map $\Th^{\text{Quiver}}$
to an appropriate, more general, context as follows.

Let $T$ be any nonempty subset of $I$.
Generalizing the definition of Calogero-Moser quiver,
let $Q_T$ be a quiver obtained from $Q$ by adjoining a vertex $s$,
called the \emph{special} vertex,
and arrows $b_i : s \to i$, one  for each $i\in T$.
Recall that $e_i$ denotes the idempotent in the path algebra
corresponding to a vertex $i$. Thus, given  $\la\in\C^I$, we write
$\la=\sum \la_ie_i,$ and we also put
$e_T := \sum_{i\in T} e_i$.

In \S\ref{qf},
for any $n\geq 1,\,\la\in\C^I,\,\nu\in\C,$
we introduce  an exact functor
\beq{G}G':\
 \Mod{\sfA_{n,\la,\nu}(Q)}\too
\Mod{\Pi'_{\la-\nu e_T +n\nu e_s}(Q_T)}.
\eeq

The construction of reflection functors 
for  generalized  preprojective algebras, see \eqref{gan},
 implies readily that, for any $i\in I$,
one has the following commutative diagram
\beq{diag_commutes}
\xymatrix{
\Mod{\sfA_{n,\la,\nu}(Q)}
\ar[rr]^{F_i}
\ar[d]_{G'} && \Mod{\sfA_{n,r_i(\la),\nu}(Q)}
\ar[d]^{G'} \\
\Mod{\Pi'_{\la- \nu e_T + n\nu e_s}(Q_T)}
\ar[rr]^{F'_i}
&& \Mod{\Pi'_{r_i(\la)- \nu e_T + n\nu e_s}(Q_T)}.
}
\eeq

The functor \eqref{G} is a generalization of the functor $M\mto\wt M$ considered in
\S\ref{pp} in the following sense. Let
$Q$ be the extended Dynkin quiver associated to
a finite subgroup $\G\sset SL_2(\C).$
Given a data $(n,k,c)$, as in \eqref{data},
put $t=1/|\G|$ and $\nu=k\cdot |\G|/2.$
The generalized preprojective algebra $\sfA_{n,\la,\nu}(Q)$
is {\em Morita equivalent}, according to \cite{GG},
to the symplectic
reflection algebra $\hh_{t,k,c}(\GG),$
so one has a  category
equivalence
$\Mod{\hh_{t,k,c}(\GG)}\iso\Mod{\sfA_{n,\la,\nu}(Q)}.$
Therefore, composing this
equivalence
with \eqref{G}, yields  a functor
$\dis\Mod{\hh_{t,k,c}(\GG)}\too\Mod{\Pi'_{\la-\nu e_T +n\nu e_s}(Q_T)}.$
The latter functor reduces, in the special case
of the one point set $T=\{o\},$ to the functor $M\mto \wt M$  considered in
\S\ref{pp}.

\subsection{Quantization of the Hilbert scheme of points on the resolution
of Kleinian singularity.}
The shift functor \eqref{SS} is the $\G$-analogue of the 
shift functor introduced in \cite{BEG} 
in the case of the trivial group $\Gamma$. The latter functor has been used by
 Gordon-Stafford \cite{GS}  to construct 
quantization of the
Hilbert scheme of $n$ points of the plane $\C^2$.

Now, let $X\to L/\G$ be the minimal resolution
of the Kleinian singularity $L/\G$ and
let $\Hilb^nX$ be the
Hilbert scheme of $n$ points in $X$.
It should be possible to use the shift functor \eqref{SS} and
  Theorem \ref{main}   to construct 
quantizations of  $\Hilb^nX$.
This would provide a common
 generalization to the case of wreath-products $\GG=S_n\ltimes \G^n$ of the results
 of Gordon-Stafford \cite{GS} in the special case $\Gamma=1$ and $n\geq 1$,
 and also of the results  of
Boyarchenko \cite{Bo} in the special case of arbitrary $\G\sset SL_2(\C)$  and $n=1$,
cf. also~\cite{Mu} for the case of cyclic group $\G$ (and  $n=1$).

In a different direction,  the construction of the  algebra
 $\e\hh_{t,k,c}(\GG)\e$ in terms of   Hamiltonian reduction
provided by Theorem \ref{main} gives way to  applying the
machinery of \cite{BFG} to symplectic reflection
algebras over $\k$,
an algebraic closure of the finite field $\BF_p$.

In more detail, fix a finite group $\G\sset SL_2(\C)$ and
a positive integer $n$. Then, a  routine argument
shows that, for all large
enough primes $p> n$,
each of the schemes $X,\,\Hilb^nX$,
and  $\bmu\inv(0)$, cf. \eqref{moment}, has a  well defined reduction
to a reduced scheme over $\k$.
Further, let
${\mathscr M}_0$ be the irreducible component 
of $\bmu\inv(0)$, cf. \eqref{moment},
 as defined in \cite[Theorem 3.3.3(ii)]{GG2}.
Then, the action of the  group $\GL(\al)/G_m$ 
on ${\mathscr M}_0$ generically free.
Moreover, according to H. Nakajima, 
there exists a $\GL(\al)$-stable
Zariski open dense subset
$\FM\sset {\mathscr M}_0,$ of {\em stable} points, such
that one has  a smooth universal geometric
quotient morphism $\FM\to \Hilb^nX.$ Furthermore, in this case all
the {\em Basic assumptions}   of \cite[4.1.1]{BFG} hold.

Next, let  $\BQ[\G]$ be the group algebra of $\G$ with
{\em rational} coefficients.
Write $Z(\G,\BQ)$ for the center
of  $\BQ[\G]$,
and  $Z_o(\G,\BQ)$ for the corresponding codimension
1 subspace, cf. \eqref{c}. 
Fix $k\in \BQ$
and $c\in Z_o(\G,\BQ)$ and let
 $\e\hh_{t,k,c}(\GG,\BQ)\e$ be the  $\BQ$-rational
version of the $\C$-algebra
 $\e\hh_{t,k,c}(\GG)\e$.
Then,  there exists a  large
enough  constant
$N(k,c)>\opp{max}(n, |\G|)$
such that for all primes $p> N(k,c)$ the $\BQ$-algebra
 $\e\hh_{t,k,c}(\GG,\BQ)\e$ has a well defined
reduction to a $\k$-algebra  $\e\hh_{t,k,c}(\GG,\k)\e$.

On the other hand, one can 
apply  a characteristic $p$ version of quantum
Hamiltonian reduction, as explained in \cite[\S3]{BFG},
in our present situation. This way,  for all  large
enough primes $p$, Theorem 4.1.4 from \cite{BFG}
provides a construction of a sheaf of Azumaya algebras $\FA_{k,c}$ on  $(\Hilb^nX)^{(1)},$
the Frobenius twist of the scheme $\Hilb^nX$.

Mimicing the proof of \cite[Theorem 7.2.4(i)-(ii)]{BFG}, and using
our Theorem \ref{main}, one obtains the following result

\begin{theorem}\label{bfg} Fix $k\in \BQ$
and $c\in Z_o(\G,\BQ)$. Then, there exists a  constant
$d(k,c)>\opp{max}(n, |\G|)$, such that for all primes $p> d(k,c)$
and $t=1/|\G|\in\k,$
we have a natural algebra isomorphism
$$H^0\big((\Hilb^nX)^{(1)},\,
\FA_{k,c}\big)\cong \e\hh_{t,k,c}(\GG,\k)\e,
\quad\text{moreover},\quad
H^i\big((\Hilb^nX)^{(1)},\,
\FA_{k,c}\big)=0,\;\en\forall i>0.
$$
\end{theorem}

\subsection{Acknowledgments.}{\small
We are grateful to Iain Gordon for a careful reading of a preliminary
draft of the paper.
The work of P.E., W.L.G., and V.G.
 was  partially supported by the NSF grants
 DMS-0504847, DMS-0401509, and DMS-0303465,
 respectively.
The work of P.E., V.G., and A.O. was partially supported by
the CRDF grant RM1-2545-MO-03.}

\section{Calogero-Moser quiver}

\subsection{Intertwiners.}\label{intertwiners}
Let $\QQ$ be the double quiver of $Q$, obtained from
$Q$ by adding a reverse edge $a^*:j\to i$
for each edge $a:i\to j$ in $Q$. For any edge $a:i\to j$
in $\QQ$, we write its tail $t(a) := i$ and its head $h(a) := j$.

We have an identification $L\stackrel{\sim}{\to}
L^*: u \mapsto \om (u,\cdot)$.
Let $a\in \QQ$ be an edge, then for each intertwiner
$ \phi_a: L\ot N^*_{t(a)} \to N^*_{h(a)} $,
we have a corresponding intertwiner
$ \phi'_a: N^*_{t(a)} \to L\ot N^*_{h(a)} $.

Suppose $Q$ is not of type $\wt A_1$.
Following \cite{CBH} (cf. also \cite{Me}), 
we normalize the intertwiners
such that for each edge $a\in Q$, we have
$\phi_{a^*}\phi'_a = \delta_{h(a)}\id_{N^*_{t(a)}}$, and
so $\phi_a \phi'_{a^*} = - \delta_{t(a)}\id_{N^*_{h(a)}}$.
Thus, $\phi'_a\phi_{a^*}$ is $\delta_{h(a)}$ times the projection
of $L\ot N^*_{h(a)}$ to $N^*_{t(a)}$, and
$\phi'_{a^*}\phi_a$ is $-\delta_{t(a)}$ times the projection
of $L\ot N^*_{t(a)}$ to $N^*_{h(a)}$.
Hence, for any vertex $i$,
\begin{equation} \label{rp}
\sum_{a\in Q; h(a)=i} \phi'_a\phi_{a^*}
- \sum_{a\in Q; t(a)=i} \phi'_{a^*}\phi_a
= \delta_i \id_{L\otimes N^*_i}.
\end{equation}

Suppose now that $Q$ is of type $\wt A_1$.
Then $\Gamma$ is the group with $2$ elements $1,\zeta$.
Moreover, $\zeta x= -x$ and $\zeta y=-y$.
Write the vertices of $Q$ as $o$ and $i$, where
$N_o$ is the trivial representation of $\Gamma$ and
$N_i$ is the sign representation of $\Gamma$.
We have a decomposition
$L=N^x_i \oplus N^y_i$ where $N^x_i$ is spanned by $x$ and
$N^y_i$ is spanned by $y$.
Let $\pr^x_i: L \ot N_o \to N_i$ be the projection map
to $N^x_i\ot N_o = N_i$,
and $\pr^y_i: L \ot N_o \to N_i$ be the projection
map to $N^y_i\ot N_o=N_i$.
Let $\pr^x_o: L\ot N_i \to N_o$ be the projection map
to $N^x_i\ot N_i = N_o$,
and $\pr^y_o: L \ot N_i \to N_o$ be the projection
map to $N^y_i\ot N_i = N_o$.
Denote the edges of $Q$ by $a_1$ and $a_2$.
If $a_1: o\to i$, then let
$\phi_{a_1} = \pr^y_i$ and $\phi_{a^*_1} = \pr^x_o$.
If $a_1: i\to o$, then let
$\phi_{a_1} = \pr^x_o$ and $\phi_{a^*_1} = -\pr^y_i$.
If $a_2: o \to i$, then let
$\phi_{a_2} = \pr^x_i$ and $\phi_{a^*_2} = -\pr^y_o$.
If $a_2: i\to o$, then let
$\phi_{a_2} = \pr^y_o$ and $\phi_{a^*_2} = \pr^x_i$.
It is easy to see that with these choices, we again
have (\ref{rp}).

\subsection{Quiver map.}   \label{quivermap}
For convenience, we shall fix
an isomorphism $N_i=\C^{\delta_i}$, where $\delta_i = \dim N_i$.
We have $\C\Gamma=\bigoplus_{i\in I}  \End N_i
= \bigoplus_{i\in I} \mathrm{Mat}_{\delta_i}(\C)$.
Let $e^i_{p,q}$ ($1\leq p,q \leq \delta_i$)  be the element of
$\C\Gamma$ with $1$ in the $(p,q)$-entry of the matrix for the
$i$-th summand and zero elsewhere.
Let $e_i$ be the idempotent $e^i_{1,1}$.
In particular, $e_o = \sum_{\ga\in\Gamma} \ga /|\Gamma|$,
where $o$ is the extending vertex of the affine Dynkin
quiver $Q$.
Note that $N_i = \C[\Gamma]e_i$ and
$\phi_a\in e_{h(a)} (L\ot \C[\Gamma] e_{t(a)})$.
Here, the left action of $\Gamma$ on $L\ot \C[\Gamma]$
is the diagonal one.
When $Q$ is not of type $\wt A_1$,
$\phi_a$ spans $e_{h(a)} (L\ot \C[\Gamma] e_{t(a)})$.
When $Q$ is of type $\wt A_1$ with vertices $o$ and $i$,
$e_o (L\ot \C[\Gamma] e_i)$ and $e_i (L\ot \C[\Gamma] e_o)$
are both 2 dimensional and spanned by the intertwiners $\phi_a$
which they contain.

We put $e_n :=\frac{1}{n!} \sum_{\sigma\in S_n} \sigma  \in \C[S_n]$.
For any $i\in I$, let
\beq{e} e_{i,n-1} := e_{n-1} (e_i \ot e_o \otot e_o)
\in \C[\GG],\quad\text{and}\quad
\e := e_n (e_o \otot e_o) \in \C[\GG].
\eeq
The idempotent $\e$ is same as the one that appears in 
Theorem \ref{main} of the Introduction.

For each vertex $i$ of the Calogero-Moser quiver $\Qc$,
cf. Definition \ref{CMQ}, the idempotent $e_i$
is the trivial path at the vertex $i$.
Let $\la_i$ be the trace of $t\cdot 1+ c$ on $N_i$,
let $\la = \sum_{i\in I} \la_i e_i$,
and let $$\nu := k|\Gamma|/2.$$
Let $\Pi = \Pi_{\la-\nu e_o + n\nu e_s}(\Qc)$, the
deformed preprojective algebra of $\Qc$ with parameter
$\la-\nu e_o + n\nu e_s$ as defined in \cite{CBH}.
So $\Pi$ is the quotient of the path algebra
$\C\QQc$ of $\QQc$ by the following relations:
$$ \sum_{a\in Q} [a,a^*] + bb^* = \la -\nu e_o, \qquad
b^*b = -n\nu e_s. $$

We shall define a functor from $\hh$-modules
to $\Pi$-modules. Let $M$ be a $\hh$-module.
We want to define a $\Pi$-module $\wt M$.
For each $i\in I$, let $\wt M_i := e_{i,n-1}M$. Also, let $\wt M_s
:= \e M$. If $a$ is an edge in $\QQ$, then define $a : \wt
M_{t(a)} \too \wt M_{h(a)}$ to be the map given by the element $
\phi_a \ot e_o \otot e_o \in \hh$. Define $b:\wt M_s \too \wt M_o$
to be  the inclusion map,
and define $b^*: \wt M_o \too \wt M_s$
to be $-\nu\cdot(1+s_{12} + \cdots + s_{1n})$.

\begin{lemma} \label{pimodule}
With the above action, $\wt M$ is a $\Pi$-module.
\end{lemma}
\begin{proof}
It is clear that
$(1+s_{12}+\cdots+s_{1n})e_{n-1}=ne_n$.
On $\wt M$, we have $b^*b = -n\nu$, and
$bb^* = -n\nu e_n = -\nu\cdot(1+s_{12} + \cdots + s_{1n})$.

By \cite[(3.5.2)]{GG}, we have an isomorphism
$f^{\ot n} \hh f^{\ot n} = \mathsf A_{n,\la,\nu}(Q)$
where $f = \sum_{i\in I}e_i$, cf.  Definition \ref{dea} below.
Now $f^{\ot n} M$ is a $\mathsf A_{n,\la,\nu}(Q)$-module, and
$e_{i,n-1}M = e_{i,n-1}f^{\ot n}M$, $\e M = \e f^{\ot n}M$.
The action of the edge $a: \wt
M_{t(a)} \too \wt M_{h(a)}$ is the action given by the element
$a\ot e_o^{\ot (n-1)}\in \mathsf A_{n,\la,\nu}(Q)$.

Now on $\wt M$,
at a vertex $i\neq o,s$, we have
$$ \sum_{a\in Q; h(a)=i} aa^* - \sum_{a\in Q; t(a)=i} a^*a
= \la_i $$
by the relation (i) in Definition \ref{dea}.
At the vertex $o$, we have
$$
\sum_{a\in Q; h(a)=o} aa^* - \sum_{a\in Q; t(a)=o} a^*a
=\la_o + \nu\cdot(s_{12} + \cdots + s_{1n}) =\la_o - \nu -  bb^*,
$$
using again the relation (i) in Definition \ref{dea}.
\end{proof}

It is clear that the assignment $G: M\mapsto \wt M$ is functorial.
We will give a more general construction in \S\ref{qf}.

Applying the functor  $M\mapsto \wt M$
to the $\hh$-module $\bbM$ introduced in \eqref{M} we
construct, as has been explained in
\S\ref{pp}, the algebra homomorphism $\theta^{\text{Quiver}}:
\Pi\to \mathsf B.$

Observe that $\theta^{\text{Quiver}}(b^*)$ is $0$ when $\nu=0$.
For this reason,
we shall need a slight modification of the above constructions.

Define $\Pi'$ to be the
quotient of the path algebra $\C\QQc$ by the following
relations:
$$
\sum_{a\in Q} [a,a^*] + \nu bb^* = \lambda - \nu e_o,
\qquad b^*b = -n e_s.
$$
We have a morphism of algebras
$\Pi\too\Pi'$ defined on the edges by
$$a \mapsto a  \text{ for }a\neq b^*,
\qquad b^* \mapsto \nu b^*.$$
This is an isomorphism only when $\nu\neq 0$.

Given a $\hh$-module $M$, we define
a $\Pi'$-module structure on $\wt M$
as above, except that now,
we let $b^*: \wt M_o \too \wt M_s$ be
$-(1+s_{12}+\cdots+s_{1n})$.
Hence, as above, we obtain a morphism
$\Th^{\text{Quiver}}: \Pi' \too \M, $ cf. \eqref{Thquiv}.

\subsection{Holland's map.}\label{holland}
In this subsection, we recall a construction of Holland
from \cite{Ho}.

Let $\eps_i \in \ZZ^I$ denote the coordinate vector corresponding
to the vertex $i\in I$. Let $\delta= \sum_{i\in I} \delta_i \eps_i$,
the minimal positive imaginary root of $Q$.
Let $\al := n\delta+\eps_s$, a dimension vector for $\Qc$.
Thus, $\al_i = n\delta_i$ for $i\in I$, and $\al_s = 1$.
\emph{We shall assume that $\lambda \cdot \delta=1$, that is,
$t = 1/|\Gamma|$.}

Let $e^a_{p,q}$ and $t^a_{p,q}$
($a\in \Qc$, $1\leq p \leq \al_{h(a)}$,
$1\leq q \leq \al_{t(a)}$) be, respectively, the coordinate
vectors and the coordinate functions on $\Rep_\al(\Qc)$.
We write $e^a_{q,p}$ for the transpose of $e^a_{p,q}$.
Now define a representation of $\QQc$ on $\O(\Rep_\al(\Qc))\ot N$,
the space of $N$-valued regular functions on $\Rep_\al(\Qc)$,
as follows.
For any $a\in \Qc$, we define the following
$\End\N$-valued differential operators of order
$0$ and $1$, respectively
$$ \hat{a}:= -\sum\nolimits_{p,q}\, e^a_{p,q} \ot t^a_{p,q},
\quad\text{resp.},\quad
\hat{a}^* :=
\sum\nolimits_{p,q}\, e^a_{q,p} \ot \frac{\partial}{\partial t^a_{p,q}}.$$

The assignment $a\mto \hat{a}, a^*\mto\hat{a}^*$ extends by
multiplicativity to an algebra homomorphism
\beq{wth} \wt \theta^{\text{Holland}}:
\C\QQc \too \big(\D(\Qc,\al) \ot \End\N \big)^{\GL(\al)},
\eeq
where $\C\QQc$ denotes the path algebra of the double quiver
$\QQc$. By \cite[Theorem 3.14]{Ho} and
\cite[Lemma 3.16]{Ho}, $\wt \theta^{\text{Holland}}$ descends to
homomorphisms, cf. notation in \eqref{vectors}:
$$ \theta^{\text{Holland}} : \Pi \too \T_\chi  \qquad\text{ and }\qquad
 \theta^{\text{Holland}} : e_s \Pi e_s \too \A_{\chi'}.$$
We remind that $\A_{\chi'}$ is the algebra in (\ref{achi}).

Later, we will define a homomorphism
$ \Theta^{\text{Holland}} : e_s \Pi' e_s \too \A_{\chi'} $.

\section{Radial part map}\label{Rad_sec}
From now on, we 
assume that $Q$ is not of type $\widetilde A_n$
where $n$ is even.
Equivalently, that means that $\Gamma$ has the (necessarily unique)
 {\em central
element of order $2$}, to be denoted $\zeta\in \Gamma$.

\subsection{Twisted differential operators.}\label{tdo}
Let $\TT\cong (\C^\times)^m$ be a torus
with Lie algebra $\ft:=\Lie\TT,$
and $p: X\to \X$ a principal $\TT$-bundle.
For any
$\bbt\in \Lie \TT$, the infinitesimal $\bbt$-action yields
a vector field $\bxi_\bbt$ on $X$.
Let $\D_X$ be the sheaf of algebraic differential
operators of $X$. The action of $\TT$ on $X$ makes
$\D_X$ a $\TT$-equivariant sheaf of algebras,
and we write
$\G(X,\D_X)^\TT$ for the algebra of $\TT$-invariant
global differential operators on $X$.
The assignment $\bbt\mto\bxi_\bbt$
gives a Lie algebra homomorphism
$\ft\too\G(X,\D_X)^\TT$.

Let $\rho: \ft\to \End F$ be a finite dimensional
representation.
We form $\D_{X,F}:= \D_X\ot\End F,$ a  sheaf of
associative algebras on $X$.
Let $\bxi-\rho:\ft\to\D_{X,F}= \D_X\ot\End F,$
$h\mto \bxi_h\ot\id_F-\rho(h)$
be the diagonal  Lie algebra homomorphism.
We write $\im(\bxi-\rho)$ for the image
of this homomorphism, and
$(p_*\D_{X,F})^{\ad\ft}$ for
the subsheaf of those
sections of the push-forward sheaf
$p_*\D_{X,F}$, on $\X$,  which commute
with $\im(\bxi-\rho)$.
Thus, $\im(\bxi-\rho)$ is a {\em central} subspace
of
$(p_*\D_{X,F})^{\ad\ft}$, a
 sheaf of
associative algebras on $\X$.
We write $(p_*\D_{X,F})^{\ad\ft}\cdot\im(\bxi-\rho)$
for the (automatically two-sided)
ideal in $(p_*\D_{X,F})^{\ad\ft}$
generated by the image of $\bxi-\rho$.
Thus, the quotient  $(p_*\D_{X,F})^{\ad\ft}/(p_*
\D_{X,F})^{\ad\ft}\cdot\im(\bxi-\rho)$
is a well-defined sheaf of associative algebras on $\X$.
Let
\beq{dtw}
\D(\X,p,\rho):
=\G\big(\X,(p_*\D_{X,F})^{\ad\ft}/
(p_*\D_{X,F})^{\ad\ft}\cd\im(\bxi-\rho)\big).
\eeq
be the algebra of its global sections,
to be referred to as  the algebra of {\em twisted differential
operators on $\X$} associated with the principal
$\TT$-bundle $p: X \to\X$ and representation $\rho$.

For any open set $U$ (in the ordinary, Hausdorff topology),
we write $\hol(U,F)$ for the vector
space of all holomorphic maps $U\to F$.
Given such an open subset $U\sset X$, put
\beq{hol}\hol_\rho(U):=\{f\in\hol(U, F)
\en\mid\en \bxi_h f=\rho(h) f,\en\forall h\in\ft\}.
\eeq

There is a natural action
of the algebra $\D(\X,p,\rho)$ on the vector
space $\hol_\rho(U)$ given, in local coordinates,
by  differential
operators with $\End F$-valued coefficients.

Given  a 
decomposition $F=F_1 \oplus \cdots \oplus F_r$ into a 
direct sum of $\ft$-subrepresentations, we have
$\End F=\bigoplus_{1\leq m,l\leq n} \Hom(F_m,F_l).$
This gives the induced direct sum decomposion
$$
\D(\X,p,\rho) = \bigoplus_{1\leq m,l\leq n}
\D(\X, p, \rho, F_m\rightarrow F_l).
$$
Thus, for each $(m,l)$, the direct summand
$\D(\X, p, \rho, F_m\rightarrow F_l)$ has
a natural structure of
left $\D(\X,p,\rho|_{F_l})$-module and of
right $\D(\X,p,\rho|_{F_m})$-module.

\subsection{The radial part construction.}\label{rad_general}
Let $G$ be a linear algebraic group and
 $Y$ a smooth $G$-variety. Assume in addition that
we have a smooth subvariety  $X\sset Y$ which is stable under
the action of a torus $\TT\sset G$, and we also have
a smooth morphism  $p:X\to \X$, which is
a principal $\TT$-bundle with respect to the induced $\TT$-action on $X$.
Thus, we have the following diagram
\beq{act}
\xymatrix{
\X\;&&X\;\ar@{->>}[ll]^<>(0.5){p}
\ar@{^{(}->}[rr]^<>(0.5){x\mto 1\times x}_<>(0.5){\jmath}&&
\;G\times_\TT X\; \ar[rr]^<>(0.5){g\times x\mto g(x)}_<>(0.5){\mathsf{act}}
&&\;Y.
}
\eeq

Let $\g:=\Lie G$ and let $\rho: \g\to\End F$ be a finite
dimensional representation.
For any open subset $U_Y\sset Y$, we
may consider the vector space
$\hol_\rho(U_Y)$ defined similarly to \eqref{hol}
but with respect to the $\g$-representation $\rho$.
Write $\rho_\ft=\rho|_\ft$ for the restriction
of $\rho$ to the Lie algebra $\ft=\Lie \TT.$
Restriction of functions gives the map
$$
\Res:\ \hol_\rho(U_Y)\too \hol_{\rho_\ft}(X\cap U_Y),
\quad f\mto \Res f:=\jmath^*f.
$$

Let $\D(Y,F)=\G(Y,\D_{Y,F})=\G(Y,\D_Y)\ot\End F$ be the algebra of
$\End F$-valued differential operators on $Y$. As above,
we have the Lie algebra map $\bxi-\rho: \g\to\D(Y,F)$
and the
subalgebra
 $\D(Y, F)^{\ad\g}$,   of the operators commuting with the image
of that map.

 Let $K\sset G$ be a finite
subgroup that preserves $X$ and normalizes the torus $\TT$.
The action of $K$ on $X$, resp. on $\D_Y$, induces a natural
$K$-action on $\X$, resp. on the algebra
$\D(\X,p,\rho_\ft),$ of twisted differential operators on $\X$.
We write  $\D(\X,p,\rho_\ft)^K$
for the subalgebra of $K$-invariants.

One has the following standard result.
\begin{proposition}[\textbf{Radial part map}]\label{rad} Assume that
$G$ is connected and the differential of the map $\mathsf{act},$ in \eqref{act},
 is an isomorphism at every point $\jmath(x), x\in X$.
Then,

\vi There is a natural radial part homomorphism
$\th^{\opp{Radial}}: \D(Y,F)^{\ad\g}\too \D(\X,p,\rho_\ft)^K$ such that,
for any open (in the  Hausdorff topology) subset $U_Y\sset Y$,
we have
$$\th^{\opp{Radial}}(D)\cd(\Res f)=\Res(Df),\quad\forall D\in\D(Y, F)^{\ad\g},\,
f\in \hol_\rho(U_Y).
$$

\vii The two-sided ideal
$\big(\D(Y,F)\im(\bxi-\rho)\big)^{\ad\g}$ is contained
in the kernel of the radial part map~$\th^{\opp{Radial}}.$

\viii Assume, in addition, that  $\X$ is affine
 and  the restriction of $\rho$ to $\ft$ is diagonalizable. Then,
there are canonical  algebra isomorphisms, cf. \eqref{BA},
$$
\A\big(\D(X,F),\ft,\bxi-\rho\big)\cong
\D(X,F)^{\ad\ft}/
\big(\D(X,F)\im(\bxi-\rho)\big)^{\ad\ft}
\;\iso\;\D(\X,p,\rho_\ft).
\qquad\Box
$$
\end{proposition}

\subsection{A slice in $\Rep_\al(\Qc)$.}\label{radialpartconstr}
We choose the following orientation on $Q$:
the vertex $o$ is a sink, and any other vertex is
a source or a sink.
Thus, the second order element
$\zeta$ acts by $1$ at sinks and by $-1$ at sources.
Note also that, see \eqref{vectors}
$$
\partial_i = - n \Tr|_{N_i}(\zeta), \quad i\in I.
$$

The collection of  intertwiners $\phi= (\phi_a)_{a\in Q}$
introduced in \S\ref{intertwiners} gives rise
to a linear map
$$\phi: L \to\Rep_\delta(Q),
\en\text{where} \en
\phi_a : L \to \Hom(N^*_{t(a)}, N^*_{h(a)}),\en u \mapsto \phi_a(u). $$

We  also define a linear map
$\jmath: L^n \to \Rep_\al(\Qc)$  by
$$
\jmath(u_1, \ldots, u_n)_b = (1, 1, \ldots, 1),\quad\text{and}\quad
\jmath(u_1, \ldots, u_n)_a =
(\phi_a(u_1), \ldots, \phi_a(u_n)),
\en\forall a\in Q.
$$

\begin{lemma}\label{lll}
Let $u,w\in L$. Suppose there are
$\beta_i \in \End(N_i^*)$ for $i\in I$ such that
$\phi_a(u)\beta_{t(a)} = \beta_{h(a)}\phi_a(w)$
for all edges $a\in Q$. If the $\beta_i$ are not
all zero, then
$u$ is proportional to $\ga w$ for some $\ga\in \Gamma$.
\end{lemma}

The lemma will be proved later, at the end of \S\ref{lllpf}.
From this lemma, we deduce

\begin{corollary}\label{generic}
There exists a
$\GG$-stable Zariski-open dense subset
$L^n\reg\sset L^n$ such that
$\jmath(L^n\reg)$ is contained in the set of
generic   representations of $Q$ and, moreover,
$\jmath(L^n\reg)$ meets generic $\GL(\al)$-orbits  in $\Rep_\al(\Qc)$ transversely.
\end{corollary}
\begin{proof}
First, we claim that the affine space $\jmath(L^n)$ meets
generic $\GL(\al)$-orbits in $\Rep_\al(\Qc)$.

Recall that
 generic representations of $Q$ with dimension vector $n\delta$
are direct sums of $n$ representations with dimension vector $\delta$,
it suffices to show that the subspace consisting of the representations
$\{\phi_a(u)\}_{a\in Q}$ for all $u\in L$ intersects generic
$\GL(\delta)$-orbits in $\Rep_\delta(Q)$.
By the preceding lemma, the (rational)
map $\Rep_\delta(Q)\to \mathbb P^1$ defined in \cite[Theorem 6.2]{Ri}
(which parametrizes generic orbits)
is nonconstant on $L$.

The Corollary  now follows from
the standard Bertini-Sard theorem,
cf. eg. \cite{Ve}. \end{proof}

\subsection{Discriminant function.}\label{d_fn} 
Put $\ls:=L\sminus \{0\}$.
The multiplicative group
$\C^\times$ is imbedded in $\GL(L)$ as scalar matrices,
and we have the standard
$\C^\times$-bundle
$\ls\to {\mathbb{P}}_L:=\ls/\C^\times\cong{\mathbb{P}}^1$
(projective line).
The group $S:=\C^\times \cap \G$ consists of
two elements, $S=\{1,\zeta\}$, where $\zeta\in\G$ is 
our second order element. 

Given  $\ell\in {\mathbb{P}}_L$,
write  $\G^\ell\sset\G$ for
the isotropy group of the line $\ell\sset L$.
Clearly, one has $S\sset \G^\ell.$
Therefore, $|\G^\ell|/2=|\G^\ell/S|$ is a positive integer,
and we put $\bkap_\ell:=|\G^\ell/S|-1$.
Thus, we have $\bkap_\ell>0$ if and only if
the group $\G^\ell\sset\G,$
is {\em strictly larger} than
$S$. The lines $\ell$ with this property form a
{\em finite} set ${\mathbb{P}}_L^{\text{sing}}\subset{\mathbb{P}}_L $.
For each $\ell\in
{\mathbb{P}}_L^{\text{sing}}$, we choose and fix a
vector representative
$v_\ell\in\ell\sminus\{0\}\sset L.$

We introduce
the following function
$$\Delta:=\prod\nolimits_{\ell\in {\mathbb{P}}_L^{\text{sing}}}
\om (v_\ell,-)^{\bkap_\ell}\in\C[L],$$
which is uniquely
defined up to a nonzero constant factor depending on the
choice of representatives $v_\ell,\ell\in{\mathbb{P}}_L^{\text{sing}}.$
Further, we introduce a {\em discriminant}
function on $\lreg$, defined by
\beq{Delta}
\Delta_n(u_1,\ldots,u_n):=\prod_{k=1}^n \frac{1}{\Delta(u_k)}
\prod_{m\ne l}\prod_{\gamma\in\Gamma}\frac{1}{\om (u_m,\gamma
u_l)},
\quad\text{for}\en u_1,\ldots,u_n\in L.
\eeq

Let  $\lreg$ be a Zariski open set
as in  Corollary \ref{generic}.
Shrinking this set if necessary, from now on we assume
in addition that $\lreg$ is
an {\em affine} $\TT$-stable  subset
contained in $(\ls)^n$
and, moreover, that  the denominator of the
function $\Delta_n$ does not vanish 
on $\lreg$.
Thus, the set $\lreg$ is
  $\GG\ltimes\TT$-stable, and we have $\Delta_n\in\C[\lreg].$

The natural action of the torus $\TT$ on $\lreg$ induces
an action of the Lie algebra $\ft=\Lie\TT$ on 
the coordinate ring $\C[\lreg]$.
Given $h=(h_1,\ldots,h_n)\in\ft=\C^n,$
we write the action map 
as $h:\C[\lreg]\to\C[\lreg],\, f\mto \bxi_h f$,
and also put $\Tr(h)=h_1+\ldots+ h_n.$
\begin{lemma}\label{discr} The discriminant $\Delta_n\in \C[\lreg]$
is a weight vector for the $\ft$-action, specifically,
we have
$$\bxi_h\Delta_n=(n|\G|-2)\Tr(h)\cd\Delta_n,
\quad\forall h\in\ft.
$$
\end{lemma}
\begin{proof}
Note that $\sum_{\ell\in {\mathbb{P}}_L^{\text{sing}}} \bkap_\ell
=|\Gamma|-2$. Hence,
   $\Delta$ is a homogeneous polynomial of
degree $|\Gamma|-2.$
We see that any vector  $u_m$  appears on the RHS of \eqref{Delta}
 with degree $-(|\Gamma|-2),$ in the factor $\prod\Delta(u_k)^{-1}$,
and with degree  $-(n-1)|\Gamma|$, in the factor
 $1/\prod\prod\om(u_m,\gamma u_l)$.
\end{proof}

\subsection{Compatibility with group actions.}
Let $\TT:=(\C^\times)^n$ be the torus,
and form the  wreath product $\GG\ltimes \TT=
S_n\ltimes (\G\times\C^\times)^n$.
We are going to define a group imbedding
\beq{imbed}
{\jmath_{_{\text{Lie}}}}:\
 \GG\ltimes \TT\to \GL(\al).
\eeq

To this end,
we recall the direct sum decomposition
\eqref{N}, so $\dim\N_i=\al_i$
and one may identify
the group $\GL(\al_i)$ with $\GL(\N_i)$, for any $i\in\Ic$.
Now, by the structure of group algebras,
we have the canonical algebra isomorphism
$\C[\Gamma] \iso \oplus_{i\in I} \End(N_i^*).$
Thus, we have a group  imbedding $\Gamma \into \GL(\delta)$
and, therefore, an imbedding $\Gamma^n \into
\GL(\delta) \times \cdots \times \GL(\delta) \into
\prod_{i\in I} \GL(\al_i).$
Further, we define a homomorphism
 $S_n\to \GL(\N_i)$ by $\sigma\mto \id_{N^*_i}\ot \sigma_{_{\C^n}},$
where $\sigma_{_{\C^n}}\in\GL(\C^n)$ stands for the permutation
matrix corresponding to $\sigma\in S_n$.
This way, combining together the above defined imbeddings of
 $\G^n$ and  $S_n$,
we obtain a  group imbedding
${\jmath_{_{\text{Lie}}}}:\GG\to\GL(\al)$, such that its component
$\GG\to\GL(\al_s)$,
at  the special vertex $s$, sends every element of $\GG$ to 1.

It remains to  define
the torus imbedding ${\jmath_{_{\text{Lie}}}}:\TT\to \GL(\al),$
 ${\mathbf{t}}\mapsto g({\mathbf{t}})=
\{g_i({\mathbf{t}})\in\GL(\al_i)\}_{i\in \Ic}$.
The latter is given
as  follows.
We put  $g_s({\mathbf{t}})=1,\,\forall{\mathbf{t}}\in\TT,$ and, for any $i\in I$, let
$$
g_i({\mathbf{t}}) := {\mathbf{t}}\inv\ot\id_{N^*_i}
\en \text{ if $i$ is
a source in $Q$ },\en\text{and}\en g_i({\mathbf{t}}) = \id_{\N_i} \text{ if $i$ is
a sink in $Q$,}
$$
where, for ${\mathbf{t}}=(t_1,\ldots,t_n)$
 we let ${\mathbf{t}}\inv\in\GL(\C^n)$ denote the  diagonal
matrix with diagonal entries $t_1\inv,\ldots,t_n\inv$.
We note that the image of $\TT$ under the above imbedding
is {\em not} contained in the center of
the group $\GL(\al).$

The torus  $\TT:=(\C^\times)^n$ acts naturally
on $L^n$; the element ${\mathbf{t}}=(t_1,\ldots,t_n) \in \TT$ sends
$(u_1, \ldots, u_n)\in L^n$ to $(t_1u_1, \ldots, t_nu_n)$.
This action of $\TT$ commutes with that of the group $\G^n$.
Thus, we get an action of the group
$\GG\ltimes \TT$ on~$L^n$.
It is easy to show that the group imbedding ${\jmath_{_{\text{Lie}}}}
: \GG\ltimes \TT\hookrightarrow \GL(\al)$ agrees with the slice imbedding $\jmath:
L^n\into\Rep_\al(\Qc).$ Specifically, one has
\beq{agree}\jmath(g(u))={\jmath_{_{\text{Lie}}}}(g)(\jmath(u)),
\quad\forall u\in  L^n,\,
g\in\GG\ltimes \TT.
\eeq

\subsection{The homomorphisms $\th^{\text{Radial}}$ and
$\Th^{\text{Radial}}$.}\label{bundles}

The element $\chi\in\C^\Ic$, in \eqref{vectors}, gives
a  Lie
 algebra homomorphism $\chi:\gl(\al)\to \C$.
We also have the
 tautological  representation $\tau:\gl(\al)\to\End\N$,
see \eqref{N},
and we let
 $\chi-\tau:\gl(\al)\to \End\N,$
$h\mto\chi(h)\cdot\id_\N -\tau(h),\,\forall h\in\gl(\al)$.

The group imbedding \eqref{imbed} induces
 the corresponding  Lie algebra imbedding
$\jmath_{_{\text{Lie}}}: \ft=\Lie\TT\hookrightarrow \gl(\al)$.
We let $\rho:=(\chi-\tau)\ccirc \jmath_{_{\text{Lie}}}$
be
 the pull-back of the representation $\chi-\tau$ to the
Lie algebra $\ft$
via the  imbedding $\ft\hookrightarrow
 \gl(\al)$.

 We are now in a position to apply the general radial
part construction of \S\ref{rad_general}
in our special case.
Specifically,
the $n$-th cartesian power of the
 $\C^\times$-bundle
$\ls\to {\mathbb{P}}_L$ gives
 a principal $\TT$-bundle
$(\ls)^n\to ({\mathbb{P}}_L)^n.$
We  set $X:=\lreg\sset (\ls)^n,$ and let  $\X\sset ({\mathbb{P}}_L)^n$
be the image of $X.$
Write $p: X\to\X$ for the restriction of the  bundle projection
to $X$. Thus, $\X$ is a $\GG$-stable Zariski open dense subset
of
$ ({\mathbb{P}}_L)^n$, and $p: X\to\X$ is  a principal $\TT$-bundle.

We apply Proposition \ref{rad} to
$$G=\GL(\al),\quad \TT=(\C^\times)^n,\quad
K=\GG,\quad Y=\Rep_\al(\Qc),\quad
p:X=\lreg\to\X=\lreg/\TT.
$$
Thus, we obtain an algebra homomorphism,
cf. \eqref{Tchi}:
\beq{thrad}
\th^{\text{Radial}}: \T_\chi=
\frac{\big(\D(\Qc,\al) \ot
\End\N \big)^{\GL(\al)}}{\Big( \big( \D(\Qc,\al) \ot \End\N \big)
\im(\bxi- ( \chi-\tau)) \Big)^{\GL(\al)}}
\;\too\;\D(\X,p,\rho)^\GG.
\eeq

Further, we introduce another representation   $\varrho: \ft\to\End \N,
\, h\mto\varrho(h)$
by the formula $\varrho(h):=\rho(h) + \frac{1}{2}(n|\G|-2)\Tr(h)\id_\N.$

It is easy to see that each of the direct summands
in the decomposition $\N=\oplus_{i\in\Qc}\N_i,$ cf. \eqref{N},
is stable under the  $\ft$-action
via either representation 
 $\rho$ or $\varrho$.
Thus we can write $\rho=\oplus_{i\in\Ic} \rho_i,$
 and $\varrho=\oplus_{i\in\Ic} \rho_i.$
To describe these representations
more explicitly, 
let $c_\zeta$ be 
the  coefficient in  \eqref{c} corresponding to our
second order element $\zeta\in\G$, and put
\beq{mu}
\mu:=-(\mbox{$\frac{c_\zeta\cdot|\Gamma|}{2}$}+1),\quad
\text{and}\quad\psi := \sum_{\{i\in I\;\mid\; i\text{ is a source in }
Q\}}\delta_i\cd \chi_i\in \C. 
\eeq

One finds that
 the representations $\rho_i$ and $\varrho_i$ are given by the
following explicit formulas
\begin{align*}
&\rho_i(h) = \psi\cd \Tr(h)\cd\id_{\N_i},&\varrho_i(h) =
(\mu+\mbox{$\frac{1}{2}$})\cd\Tr(h)\cd\id_{\N_i}
\qquad\text{if $i=s$ or $i$ is a sink in $Q$}, \\
&\rho_i(h)= \psi\cd \Tr(h)\cd\id_{\N_i}-
h\ot\id_{N^*_i}, &\varrho_i(h) =(\mu+\mbox{$\frac{1}{2}$})\cd \Tr(h)\cd\id_{\N_i}-
h\ot\id_{N^*_i}
\en\text{if $i$ is a source in $Q$},
\end{align*}
 where, for any $h=(h_1,\ldots,h_n)\in\ft=\C^n,$
we write $\Tr(h):=h_1+\ldots+h_n,$ and
where the tensor factor $h$, in  $h\ot\id_{N^*_i}$,
stands for the map $\C^n\to\C^n$ given by
the diagonal matrix with diagonal entries
$h_1,\ldots,h_n.$

Next, according to Lemma~\ref{params} below, we have
\beq{mupsi}
2(\mu-\psi)=-2(\mbox{$\frac{c_\zeta\cdot|\Gamma|}{2}$}+1)
+ 1 - 2\psi=(|\Gamma|-2) + (n-1)|\Gamma|=
n|\G|-2.
\eeq
Hence,  Lemma \ref{discr} shows that
$\varrho -\rho=\frac{1}{2}(n|\G|-2)\cdot\Tr(-)$ is
 nothing but  the weight of $\sqrt{\Delta_n}$, the square root of the
discriminant function. Thus,
given
 $D\in \D(\X,p,\rho)$, we may
 conjugate $D$ by the operator of multiplication by 
the function  $\sqrt{\Delta_n}$ to obtain
 a twisted differential
operator
$\frac{1}{\sqrt{\Delta_n}}\ccirc 
D\ccirc \sqrt{\Delta_n}
\in \D(\X,p,\varrho),$
such that for any open set $U\sset \lreg$ the induced action
on functions is given by the map
$$
\G(U, \varrho)\to\G(U, \varrho),
\quad
f\mto (1/\sqrt{\Delta_n})\cd\th^{\text{Radial}}(D)(\sqrt{\Delta_n}\cd f).
$$
 
We note that although the square root of the discriminant
function $\Delta_n$ is ill defined as a function,
conjugation 
by the operator of multiplication by such a function
is an unambiguously defined operation
on twisted differential operators. 
Furthermore, the result of conjugation by $\sqrt{\Delta_n}$
is clearly independent of the choice of
 nonzero constant factor involved in the definition of $\Delta_n$,
cf. \S\ref{d_fn}.
Thus,  we have a canonically defined
 algebra homomorphism
 $$\Th^{\opp{Radial}}: \
\T_\chi\too\D(\X,p,\varrho)^\GG,\quad
u\mto \Theta^{\opp{Radial}}(u):=\frac{1}{\sqrt{\Delta_n}}\ccirc 
\th^{\text{Radial}}(u)\ccirc \sqrt{\Delta_n}.
$$

\subsection{Formulas for the map
$\th^{\text{Radial}}\ccirc\theta^{\text{Holland}}$.}
\label{formulas} 
We introduce some notation. Given a map
$f: L\to U$ and any $1\le l \le n$,
we write $f_l$ for the
map $L^n\to U,\,
(u_1,\ldots,u_n)\mto f(u_l).$
Thus,  given
 $\gamma\in \G,$ we have the composite $L\stackrel{\ga}\too
 L\stackrel{f}\too U,$ and the corresponding map
$(f\ccirc\gamma)_l:L^n\to U$.

Let $\om$ denote the symplectic form on $L$.
For any vector
$v\in L$, we have the linear function
$v^\vee: u\mto \om(v,u)$.
Thus, for any
 $v\in L$ and $\gamma\in \G$,
we may consider the following
 functions 
\begin{align}\label{notation}
(\gamma v)^\vee_l:=(v^\vee\ccirc \gamma\inv)_l:&\
L^n\to\C,\quad
(u_1,\ldots,u_n)\mto \om(\gamma v,u_l)=\om(v,\gamma\inv u_l)\quad\text{and}\nonumber\\
\omg;m,l):&\
L^n\to\C,\quad (u_1,\ldots,u_n)\mto \om(u_m,\gamma u_l),
\qquad\forall 1\le m\ne l \le n.\end{align}
The definition
of the open subset $\lreg\sset L^n$ insures,
see \S\ref{d_fn}, that none of the
functions $\omg;m,l)$ vanishes on $\lreg$.
Hence,  we have $1/\omg;m,l)\in\C[\lreg].$

In \S\ref{intertwiners}, for each edge $a\in \QQ,$
we have introduced  intertwiners $
\phi_a: L\ot N^*_{t(a)} \to N^*_{h(a)}$
and $\phi'_a: N^*_{t(a)} \to L\ot N^*_{h(a)}$.
Below, we shall view $\phi_a$ as
a $\Hom(N^*_{t(a)},N^*_{h(a)})$-valued linear
function on $L$, written  as $u\mto \phi_a(u)$.
The $n$-th cartesian power of this function
gives a $\GG$-equivariant linear map
$$\phi_a^n: L^n\to
\Hom(\N_{t(a)},\N_{h(a)})=\Hom(N^*_{t(a)},N^*_{h(a)})\ot\End\C^n,\quad
(u_1,\ldots,u_n)\mto\sum_{1\leq l\leq n}\,
\phi_a(u_l)\ot E_{ll},$$
where $E_{ll}$ stands for the $n\times n$-matrix unit
with the only nonzero entry at the place $(l,l)$.

Similarly,  we shall view $\phi'_a$ as a  $\Hom(N^*_{t(a)},N^*_{h(a)})$-valued 
constant vector field on $L$ whose value at each point $u\in L$ is equal
to $\phi'_a$.
Thus, given
$m\in [1,n],$ we shall write  $\frac{\partial}{(\partial\phi_{a^*})_m}$
for the $\Hom(N^*_{h(a)},N^*_{t(a)})$-valued  first order
differential operator on $L^n$ corresponding to the
constant
 vector field
$\phi'_{a^*}\in \Hom(N^*_{h(a)}, N^*_{t(a)}) \ot L$ 
along the $m$-th direct factor $L$ in $L^n$.

Next, recall the map $\theta^{\text{Holland}}$, introduced in
\S\ref{holland}.
The composite $\th^{\text{Radial}}\ccirc\theta^{\text{Holland}}$
associates to every edge $a\in\QQ$
a  twisted differential
operator from the algebra $\D(\X,p,\rho)$.
By definition, such an operator is a coset modulo the ideal
$\D(\lreg,\N)^{\ad\ft}\cdot\im(\bxi-\rho)$,
see Proposition \ref{rad}(iii), of an element
$$\th^{\text{Radial}}\ccirc\theta^{\text{Holland}}(a)\in
\D(\lreg)\ot\Hom(\F_{t(a)},\F_{h(a)})
\sset \D(\lreg,\N)$$
We will write such an element $D$ as an
$n\times n$-matrix with entries in $\Hom(N^*_{t(a)},N^*_{h(a)})$,
and write $D_{ml}$ for $(m,l)$-th entry of
that matrix.

\begin{proposition}  \label{radialpart}
Let $a\in Q$. Then

\vi $\,\th^{\opp{Radial}}\ccirc\theta^{\opp{Holland}}(a)$ is a zero-order
differential
operator on $\lreg$ given by multiplication by the 
function $\phi_a^n$.

\vii 
For $ l \neq m$, the $(m, l)$-entry of
$\th^{\opp{Radial}}\ccirc\theta^{\opp{Holland}}(a^*)$ is
is a zero-order
differential
operator on $\lreg$ given by multiplication by the 
function
\begin{equation}  \label{radialml}
-k/2 \sum\nolimits_\ga\, \frac{(\phi_{a^*}\ccirc\ga)_ l }
{\omg;m,l)} \ga,
\end{equation}
and the $(m,m)$-entry of $\th^{\opp{Radial}}\ccirc\theta^{\opp{Holland}}(a^*)$ is
a first order differential operator 
\begin{gather}
\frac{1}{|\Gamma|}\left(\frac{\partial}{\partial (\phi_{a^*})_m} +
\sum_{\ga \neq 1,\zeta}
\frac{(\phi_{a^*}\ccirc(\ga^{-1}+\id))_{mm}}{\om(\ga;m,m)}
\left( -1 + |\Gamma|c_\ga \ga^{-1} \right)
+ \frac{1}{|\Gamma|}\sum_{\ell\neq m}  \sum_\ga
\frac{(\phi_{a^*}\ccirc\ga)_\ell}{\omega(\ga;m,\ell)}\right). \label{radialmm}
\end{gather}

\viii
For the edge $b:s\to o$ we have
\begin{align}
\th^{\text{Radial}}\ccirc\theta^{\text{Holland}}(b)
=& -\sum\nolimits_p\, e^b_{p,1}, \nonumber  \\  \label{radialbb}
\th^{\text{Radial}}\ccirc\theta^{\text{Holland}}(b^*)
=& (1- \sum\nolimits_{j\in I}\, \delta_j\chi_j) \sum\nolimits_p\, e^b_{1,p}
= \nu \sum\nolimits_p \,e^b_{1,p}.
\end{align}
\end{proposition}

The proof of  Proposition \ref{radialpart} will be given later,
in \S\ref{pf_radialpart}.  We will use

\begin{lemma}\label{params}
We have $c_\zeta +n =  (1-2\psi)/|\Gamma|,$ and
$k = 2(1 - \sum_j \delta_j \chi_j)/|\Gamma|.$ Furthermore,
$$
c_\ga = \Big( 1- \sum\nolimits_j\, \chi_j
(\delta_j- \Tr|_{N^*_j}(\ga)) \Big)/|\Gamma|
\quad\text{ for } \ga \neq\zeta .
$$
\end{lemma}
\begin{proof}
Since $\la_i = \Tr|_{N_i}(t\cdot 1+c)$, we obtain by
orthogonality relations that
$c_\ga =$\newline
$\dis 1/|\Gamma| \sum\nolimits_{j\in I}\, \la_j \Tr|_{N^*_j}(\ga). $
Hence, for $\pa$ as in \eqref{vectors}, we compute
\begin{align*}
& \sum\nolimits_j \chi_j (\delta_j - \Tr|_{N^*_j}(\ga))  \\
=& \la\cdot \delta - \nu - \partial \cdot \delta
- \sum\nolimits_j \la_j \Tr|_{N^*_j}(\ga)  + \nu
- n \sum\nolimits_j \Tr|_{N_j}(\zeta) \Tr|_{N^*_j}(\ga)  \\
=& 1 - |\Gamma|c_\ga - n \sum\nolimits_j \Tr|_{N_j}(\zeta) \Tr|_{N^*_j}(\ga) .
\end{align*}
If $\ga\neq \zeta$, then this is equal to
$1 - |\Gamma|c_\ga$. If $\ga=\zeta$, then it is equal to
$1 - |\Gamma|c_\zeta - n |\Gamma|$.
Moreover,
$$  \sum_j \chi_j (\delta_j - \Tr|_{N^*_j}(\zeta))
= 2\sum_{\{j\in Q\,\mid\, j \text{ is a source } \}} \chi_j \delta_j = 2\psi. $$

We also have
$ \sum_j \chi_j \delta_j= 1- \nu = 1- k|\Gamma|/2. $
\end{proof}

\section{Dunkl representation}
\subsection{Dunkl operators.}\label{thth}
Recall the principal $\TT$-bundle
$p: X=\lreg\to\X=\lreg/\TT,$ and other notation introduced in
\S\ref{bundles}. We are going to
define a certain   representation  $\eta$ of the Lie algebra
$\ft=\Lie\TT$ which is normalized by the natural $\GG$-action on
$\ft$. Thus,  there is an action of $\GG$
on $\D(\X,p,\eta)$, 
the corresponding algebra of twisted differential operators.

Our goal is to define certain elements in
$\D(\X,p,\eta)\rtimes\GG$ which may be thought of 
as  $\G$-analogues of {\em Dunkl operators}.
The construction of these operators will be given
in five steps.

\step{1.} Write $\LL$ for the preimage of
$\P_L\sminus\P^\text{sing}_L$ under
the projection $L^\times\onto \P_L$.
Thus $\LL\sset L$ is an
open dense subset  formed by
the points $v\in L$ such
that,
 for any
$\gamma \in \G\sminus\{1,\zeta\},$
we have $\gamma(v)\not\in\C v$.

For any $\gamma\in\G$, we have a quadratic
function $\om^\ga: L\to\C, u\mto \om(u,\gamma
u).$ This function does not vanish on $\LL$,
thus, we have $1/\om^\ga\in\C[\LL].$
Given $v\in L,$ we also have the linear function
$v^\vee: u\mto \om(v,u),$ on $L$.

Recall the coefficients $c_\gamma\in\C$ given by \eqref {c}.
To each $v\in L$ we assign the following element \beq{rk1end}
D^v:=2|\Gamma|^{-1}\frac{\partial}{\partial v}+\sum_{\gamma\ne
1,\zeta} c_{\gamma}\frac{(\gamma v+v)^\vee}{\om^\ga}\gamma \in
\D(\LL)\rtimes\G.
\end{equation}

\step{2.} Let $F=\C^2$ be a 2-dimensional vector space with fixed
basis $(f^+,f^-)$, and identify $\End F$ with the algebra of
$2\times 2$-matrices. We set $\D(\LL,F)=\D(\LL)\ot\End F$,
and form the smash product $\D(\LL,F)\rtimes\G=
(\D(\LL)\otimes\End F)\rtimes\G,$ where
$\G$ acts trivially on $F$ and on $\End F$.

 For each $v\in L$, we introduce
the following element written as a matrix with entries in
${\D(\LL)\rtimes\G}$: \beq{DE} D^v_F:=\begin{pmatrix} 0& - v^\vee\\
D^v&0\end{pmatrix}= D^v\cd\begin{pmatrix} 0& 0\\  1&0\end{pmatrix}
- v^\vee\cd\begin{pmatrix} 0& 1\\  0&0\end{pmatrix}
\in\TD(\LL,F)\rtimes\G. \eeq

\step{3.} For any affine
variety ${\mathcal Y}$ and $n\geq 1,$ one has the standard algebra isomorphism
$\D({\mathcal Y}^n)\cong \D({\mathcal Y})^{\ot n}$.
Since
$\End(\FF)\cong(\End F)^{\ot n}$,
we deduce an
 algebra isomorphism
$\D({\mathcal Y}^n,\FF)$
$\cong \D({\mathcal Y},F)^{\ot n}$.
The symmetric group $S_n$ acts naturally on
${\mathcal Y}^n$ and on $(\End F)^{\ot n}$, hence,
also on the tensor product $\D({\mathcal Y}^n,\FF)\cong \D({\mathcal Y},F)^{\ot n}$.

We take ${\mathcal Y}=\LL$ and put $X:=\lreg,$
cf. \S\ref{d_fn}.
Thus, $X$ is a $\GG$-stable affine open dense subset of $(\LL)^n$,
and we have a chain of natural
 algebra imbeddings
$$\big(\D(\LL,F)\rtimes\G\big)^{\ot n}=\D((\LL)^n,\FF)\rtimes\G^n
\into \D(X,\FF)\rtimes\G^n\into
\D(X,\FF)\rtimes\GG,$$
where the middle imbedding is given by
 restriction from $(\LL)^n$ to $X$.

For any $v\in L$ and $l=1,\ldots,n,$
let $D^{v,n}_{F,l}\in \D(\lreg,\FF)\rtimes\GG$ denote the image of
the element $1^{\ot (l-1)}\ot D^v_F\ot 1^{\ot (n-l)}$,
cf. \eqref{DE}, under this imbedding.

\step{4.} 
For any $l=1,\ldots,n,$ and  $\gamma\in\G$, let $\gamma_l\in \G^n$
denote a copy of the  element $\gamma$
placed in the $l$-th factor of $\G^n$.
In particular, given any  $1\le m\ne l \le n$,
and $\gamma\in\G$, we have the corresponding transposition
$s_{ml}=\{m\leftrightarrow l\}\in S_n$
and the element $s_{ml}
\gamma_m\gamma_l^{-1}\in\GG$.
Given $\gamma\in \G,\,v\in L$ and
 any  $1\le m\ne l \le n$, we also have 
regular functions $v_l^\vee$ and
$1/\omg;m,l)$ on $\lreg,$
see \eqref{notation}. 

With this notation, for any $v\in L$ and $1\le m\ne l \le n$, we
will now define an element \beq{R} R^v_{ml}\in (\C[X]\ot
\End \FF)\rtimes\GG=
\Hom\big(\FF,(\C[X]\ot\FF)\rtimes\GG\big).
\eeq 

To this end, write
\beq{basis}
\{\overset{{\,}_\to}f=f_1\ot\ldots\ot f_n\;\mid\;f_i=f^\pm,\, i=1,\ldots,n\}.
\eeq
for  the standard basis of the vector  space $\FF$.
Given such
 a basis element $\overset{{\,}_\to}f=f_1\ot\ldots\ot f_n$
and $1\leq m\leq n$, let
$\overset{{\,}_\to}f^{_{_-}}_m:=f_1\ot\ldots\ot f_{m-1}\ot
f^-\ot f_{m+1}\ot\ldots\ot f_n$. 
Now,  view each $R^v_{ml}$ in \eqref{R} as a linear map
$\FF\to(\C[X]\ot\FF)\rtimes\GG$,
which we  define as follows

\begin{align}\label{rml}
R^v_{ml}(\overset{{\,}_\to}f)
&=\frac{1}{2}\sum\nolimits_{\gamma\in\Gamma} \big(\frac{(\gamma
v_l)^\vee}{\omg;m,l)}\ot\overset{{\,}_\to}f^{_{_-}}_m\big)
 \cd
s_{ml}
\gamma_m\gamma_l^{-1},
\en\text{ if }\en f_m=f^+, f_l=f^+,\nonumber\\
R^v_{ml}(\overset{{\,}_\to}f)
&=\frac{1}{2}\sum\nolimits_{\gamma\in\Gamma}
\big(\frac{v_m^\vee}{\omg;m,l)}\ot\overset{{\,}_\to}f^{_{_-}}_m\big)
 \cd s_{ml}\gamma_m\gamma_l^{-1},
 \en\text{ if }\en f_m=f^+, f_l=f^-,\nonumber\\
R^v_{ml}(\overset{{\,}_\to}f ) &=0, \en \text{ if }\en f_m=f^-.
\end{align}

We identify the algebra $\C[X]\ot\End\FF$ with
the subalgebra of $\D(X,\FF)$ formed by zero order
differential operators. Therefore, we may view the elements
$R^v_{ml}$, in \eqref{rml}, as being elements of
 $\D(X,\FF)\rtimes\GG$, which have
zero order as  differential operators. Given $k\in\C$ and
$v\in L$, we  define  first  order
differential operators 
\beq{dunop} \opp{Dunkl}^v_l:=D^{v,n}_{F,l}
+k\sum\nolimits_{l\ne m} R_{lm}^v\in \D(X,\FF)\rtimes\GG,
\quad\forall 1\le l \le n. 
\eeq

\step{5.} 
Let  $\mu\in\C$
be the constant defined in \eqref{mu}.
We introduce a representation of
 the 1-dimensional Lie algebra
$\C$  on the  vector space $F$.
Specifically, we let
the generator $\bone\in\C$ act, in the basis $\{f^+,f^-\}$, via the
diagonal matrix
$\diag(\mu+\mbox{$\frac{1}{2}$},\mu-\mbox{$\frac{1}{2}$})$.
The $n$-th
cartesian power of this 2-dimensional representation  gives a
Lie algebra representation $\eta: \ft=\C^n\to\End \FF$.

We consider the Lie algebra
homomorphism
$\dis\bxi-\eta:\ft\to \D(X,\FF),
h\mto\bxi_h\ot\id_\FF-\id_\D\ot\eta(h).$
The group $\GG$ acts naturally both on the
Lie algebra $\ft$ and on $\D(X,\FF)$,
and it is clear that the map
$\bxi-\eta$ is  $\GG$-equivariant. It follows in particular
that $\D(X,\FF)^{\ad\ft}$,
the centralizer of the image of the map $\bxi-\eta$ in
$\D(X,\FF)$,
 is a $\GG$-stable subalgebra.

Now, we apply the general construction of algebras of twisted
differential operators given in \S\ref{tdo} to the torus
$\TT=(\C^\times)^n$ acting on $X=\lreg$ and to the representation
$\eta$ defined above. This way, with the  notation of
\S\ref{bundles}, we  get an algebra $\D(\X,p,\eta)$.
By construction, the algebra $\D(\X,p,\eta)$
 is a quotient of $\D(X,\FF)^{\ad\ft}$, and this quotient
inherits
 a natural $\GG$-action.
Thus, we may form the
corresponding algebra $\D(\X,p,\eta)\rtimes\GG$.

It is straightforward
to verify, counting homogeneity degrees of the coefficients, that
for any $v\in L$ the operator in \eqref{dunop} is
$\ad\ft$-invariant. That is, for each $l=1,\ldots,n$, we have
$\opp{Dunkl}^v_l\in \D(X,\FF)^{\ad\ft}\rtimes\GG.$
Therefore, the element $\opp{Dunkl}^v_l$ has a well
defined image in
$\D(\X,p,\eta) \rtimes\GG,$ to be denoted by the same
symbol $\opp{Dunkl}^v_l$ and to be  called the $l$-th {\em Dunkl
operator} associated with $v\in L$.

\subsection{Equalizer construction.}\label{equal}
Recall from \S\ref{d_fn}, the group
$S=\{1,\zeta\}\cong\mathbb{Z}/2$.
Thus, $S=\C^\times\cap\G\sset\GL(L)$ may be (and will be)  viewed
either as a subgroup of $\C^\times$ or as a subgroup of $\G$.
We put $\ss:=S^n$. This group  comes equipped
with a natural group
 imbedding
$\beps:\ss\into \G^n\sset \GG,$ such that
the image of $\ss$ is a
{\em normal} subgroup in $\GG$,
and also with a natural imbedding $\ss\into\TT$.

In general, let
$A$ be an associative algebra
equipped with a
$\GG$-action $\GG\ni g: a\mto a^g$,
by algebra automorphisms,
and also with a
homomorphism $\at: \ss\to A$,
that is, with a map such
that $\at(1)=1,$ and such
that 
$\at(ss')=\at(s)\cd\at(s'),\,\forall
s,s'\in\ss.$
Assume in addition that the following  identities hold
(the one on the left says that $\at$ is a
$\GG$-{\em equivariant} map):
\beq{agat}
\at(s)^g=\at(gsg\inv),
\quad\text{and}\quad
\at(s)\cd a\cd\at(s\inv)=a^{\beps(s)},\quad\forall
s\in\ss,\,g\in\GG,\,a\in A.
\eeq

We form the  smash product $A\rt\GG$
and introduce the following
 two homomorphisms
$$\Upsilon_1, \Upsilon_2:\
 \ss\to A\rtimes\GG,\quad
\text{where}\quad\Upsilon_1:
s\mto \at(s)\rtimes 1,\quad
\Upsilon_2: s\mto 1\rtimes \beps(s).$$

It is straightforward to verify that
equations
\eqref{agat} imply that
the left ideal in the algebra
$A\rtimes\GG$ generated by
the set $\{\Upsilon_1(s)-\Upsilon_2(s),\, s\in\ss\}$
is in effect a two-sided ideal.
We let $A\ttimes\GG$ be the
quotient of  $A\rtimes\GG$
by that  two-sided ideal,
to be called
 {\em equalizer smash product} algebra.

\subsection{The homomorphism $\th^{\opp{Dunkl}}$.}\label{xxx}
Let $$\mathbb{I} = \{(\epsilon_1, \ldots, \epsilon_n) \mid
\epsilon_m = 0 \mbox{ or } 1 \mbox{ for all }m\} \subset \ZZ^n.$$
Let $\bbe=(\epsilon_1, \ldots, \epsilon_n) \in \mathbb I$,
and write $F_\bbe$ for the 
one dimensional subspace of $F^{\ot n}$
spanned by $f_1\ot \cdots\ot f_n$, where
$f_m = f^+$ if $\epsilon_m=0$, and $f_m=f^-$ if $\epsilon_m=1$
($m=1, \ldots, n$).
The representation $\eta$ of $\mathfrak t$ on
the vector space $F^{\ot n}$ induces an adjoint action of
$\mathfrak t$ on $\End(F^{\ot n})$.
We have a decomposition
$$
\End(F^{\ot n}) = \bigoplus\nolimits_{\bbe, \bbe'\in \mathbb I}\,
F_{\bbe'} \ot (F_{\bbe})^*,
$$
where each component in the direct sum is stable under
the action of $\mathfrak t$.
Moreover, $\mathfrak t$ 
acts on $F_{\bbe'} \ot (F_{\bbe})^*$ by the character
$\bbe'-\bbe$. 
Therefore, the $\mathfrak t$-action on
$F_{\bbe'} \ot (F_{\bbe})^*$ exponentiates to
a $\mathbb T$-action. 
Taking the direct sum over various pairs
$(\bbe,\bbe')$,
 we obtain a well defined {\em adjoint} $\TT$-action
on $\End F^{\ot n}=F^{\ot n}\ot (F^{\ot n})^*.$
Thus, for any $\btt\in\TT$,  the
adjoint action of $\btt$ gives an algebra automorphism
 $\Ad_{_F}(\btt): \End F^{\ot n}\to \End F^{\ot n}.$

The torus $\TT$ also acts naturally on $X=\lreg$.
The tensor product of the induced $\TT$-action
on $\D(X)$ with the  adjoint $\TT$-action
on $\End F^{\ot n}$ gives
a $\TT$ action on $\D(X,F^{\ot n})=\D(X)\ot\End F^{\ot n}$, to be
called the {\em adjoint} action 
$\Ad_{_{\D\ot F}}: \TT\to\opp{Aut}(\D(X,F^{\ot n})).$ 
The map $\Ad_{_{\D\ot F}}$ is clearly $\GG$-equivariant.
It is also clear from the construction that
the differential of the $\Ad_{_{\D\ot F}}$-action of $\TT$
is nothing but the adjoint action of the Lie algebra
$\ft$.
In particular, we have $\D(X,F^{\ot n})^{\ad\ft}=
\D(X,F^{\ot n})^{\Ad_{_{\D\ot F}}\TT}.$

Next, we are going to apply the general construction of
\S\ref{equal}
in the following special case.
Let
$S\to \D(\LL,F)=\D(\LL)\ot\End F$ 
be a homomorphism given by the assignment
$$1\mto 1_\D\ot\id_F, \quad
\zeta \mto 1_\D\ot \begin{pmatrix} 1& 0\\
0& -1\end{pmatrix}.
$$
We  define a homomorphism
$\at_{_F}:\ss\to \D(X,\FF)$ to
be the composite of the $n$-th cartesian power of the above homomorphism
$S\to \D(\LL,F)$,
followed by the natural imbedding $\D((\LL)^n,\FF)\into \D(X,\FF)$.
This homomorphism is clearly $\GG$-equivariant
and the image of  $\at_{_F}$ is
contained in
$\D(X,\FF)^{\ad\ft}.$

 Write $a\mto a^g$ for the action of an element
$g\in\GG$ on $a\in\D(X,\FF).$
One checks by direct computation
that the map $\at_{_F}$ is related
to the two natural imbeddings $\bept: \ss\into\TT$
and $\beps: \ss\into\GG$
 via the following
identity 
\beq{kkk}
\big(\Ad_{_{\D\ot F}}\ccirc\bept(s)\big)(a)=
\at_{_F}(s\inv)\cd a^{\beps(s)}\cd\at_{_F}(s),\quad\forall s\in \ss,\, a\in\D(X,\FF).
\eeq
It follows from \eqref{kkk} that, for any $a\in\D(X,\FF)^{\ad\ft}=
\D(X,\FF)^{\Ad_{_{\D\ot F}}\TT}$ and $s\in\ss$, one has
$\dis a^{\beps(s)}=\at_{_F}(s)\cd a\cd\at_{_F}(s\inv)$.
We conclude
that both conditions in \eqref{agat} hold for
the thus obtained  homomorphism  $\at:\ss\to A:=\D(X,\FF)^{\ad\ft}$.

Further, we have the algebra projection
 $\pr:\D(X,\FF)^{\ad\ft}\onto\D(\X,p,\eta)$,
and we let
$\pr\ccirc\at: \ss\to\D(\X,p,\eta)$
be the composite homomorphism.
The $\Ad_{_{\D\ot F}}$-action of  $\TT$ on $\D(X,\FF)$ clearly descends to
an $\Ad_{_{\D\ot F}}$-action on $\D(\X,p,\eta)$.
It follows that  conditions \eqref{agat} hold for
the map $\pr\ccirc\at$ as well.
 Thus, we are in a position to form
$\D(\X,p,\eta)\ttimes\GG$,
the corresponding equalizer smash product.
We let  ${\opp{Dunkl}}(v,l)$ denote the image in $\D(\X,p,\eta)\ttimes\GG$
of
the element $\opp{Dunkl}^v_l\in \D(\X,p,\eta)\rtimes\GG$.

The main result of this section reads

\begin{theorem}\label{mess} Put $t=1/|\Gamma|$. The assignment,
given on the generators $g\in\G,\,v_l\in L_l,\,
l=1,\ldots,n$ (where $L_l$ stands for the $l$-th direct summand in
$L^n$),
of the algebra $\hh_{t,k,c}(\GG)$
 by the formulas below
extends uniquely to a well defined and injective
algebra homomorphism
$$
\th^{\opp{Dunkl}}:\ \hh_{t,k,c}(\GG)\too
\TD(\X,p,\eta)\ttimes\GG, \quad g\mto g,\;
v_l\mto{\opp{Dunkl}}(v,l).$$
\end{theorem}

The injectivity statement in the theorem is not difficult;
it follows easily from the PBW theorem for the algebra
 $\hh_{t,k,c}(\GG)$, by considering principal symbols of differential
operators. The difficult part is to verify that the assignment
of the theorem does define an algebra homomorphism.
The proof of  this is quite long and involves
a lot of explicit computations. That
 proof  will be given later, in  \S\ref{pf_mess}.
In the  special
 case $n=1$,  the proof is less technical and is presented below.

\subsection{Proof of Theorem \ref{mess} in the special
case: $n=1$.}\label{ddd}
Let $u_1,u_2$ denote the coordinates in the
symplectic vector space $(L,\om)$.

For $n=1,$ the assignment of Theorem \ref{mess} reduces to the
 map  $L \to \D(\LL,F)\ttimes\Gamma$ that reads
$$
\th^{\opp{Dunkl}}:
v\mapsto\begin{pmatrix} 0& -v^\vee\\  D^v&0\end{pmatrix},
\en\text{where}\en
D^v=\frac{2}{|\Gamma|}\frac{\partial}{\partial v}+\sum_{\gamma\ne
1,\zeta} c_{\gamma}\frac{(\gamma v+v)^\vee}{\om^\ga}\gamma.
$$

For any $v,w\in L$ we are going to  compute all 4 entries of the
$2\times 2$-matrix representing
the operator $[\theta^{\text{Dunkl}}(v),\theta^{\text{Dunkl}}(w)]\in
\D(\LL,F)\rtimes\G$.
First of all, it is easy to see that
$\dis[\theta^{\text{Dunkl}}(v),\theta^{\text{Dunkl}}(w)]_{12}$
$\dis=
[\theta^{\text{Dunkl}}(v),\theta^{\text{Dunkl}}(w)]_{21}=0.$

Next, write
 $\eu=u_1\frac{\partial}{\partial
u_1}+u_2\frac{\partial}{\partial u_2}$ for the Euler operator.
We compute
\begin{align*}
& [\theta^{\text{Dunkl}}(v),\theta^{\text{Dunkl}}(w)]_{11}
=w^\vee D^v- v^\vee D^w  \\
=& \frac{2}{|\Gamma|}\left(w^\vee
\frac{\partial}{\partial v}- v^\vee \frac{\partial}{\partial
w}\right)
+ \sum_{\gamma\ne 1,\zeta}\frac{c_\gamma}{\om^\ga}
\big(w^\vee (\gamma v+v)^\vee - v^\vee (\gamma
w+w)^\vee\big) \gamma \\
=& -\frac{2\omega(v,w)}{|\Gamma|}\eu+
\sum_{\gamma\ne 1,\zeta}\frac{c_\gamma}{\om^\ga}
(w^\vee (\gamma v)^\vee-(\gamma w)^\vee
v^\vee))\gamma=
\omega(v,w)\big(-\frac{2}{|\Gamma|}\eu
+\sum_{\gamma\ne1,\zeta}c_\gamma
\gamma\big).
\end{align*}

One proves similarly that
$
[\theta^{\text{Dunkl}}(v),\theta^{\text{Dunkl}}(w)]_{22}=\omega(v,w)
(\frac{2}{|\Gamma|}(\eu+2)+\sum_{\gamma\ne
1,\zeta}c_\gamma \gamma).
$ 
Thus, we find
$$
[\theta^{\text{Dunkl}}(v),\theta^{\text{Dunkl}}(w)]=
\frac{2\omega(v,w)}{|\Gamma|}\left(\mbox{\small$\left(\begin{matrix} -1& 0\\
0&1\end{matrix}\right)$}\eu
+\mbox{\small$\left(\begin{matrix} 0& 0\\
0&2\end{matrix}\right)$}\right)+\omega(v,w)\sum_{\gamma\ne
1,\zeta}c_\gamma \gamma.
$$

Now, in the 1-dimensional Lie algebra
$\ft=\C$, we have the generator $\bone$
which acts in $F$ via the matrix
$\diag(\mu+\mbox{$\frac{1}{2}$},\mu-\mbox{$\frac{1}{2}$}).$
By  definition of twisted differential operators,
in the algebra $\D(\X,p,\eta)$, we have
$\eu=\bone=\diag(\mu+\mbox{$\frac{1}{2}$},\mu-\mbox{$\frac{1}{2}$}).$ Therefore,
in  the algebra $\D(\X,p,\eta)$, we get
$$
[\theta^{\text{Dunkl}}(v),\theta^{\text{Dunkl}}(w)]=
\frac{2\omega(v,w)}{|\Gamma|}\left(
 \mbox{\small$\left(\begin{matrix} -\mu-\half & 0\\
0&\mu-\half \end{matrix}\right)$}
+
\mbox{\small$\left(\begin{matrix} 0& 0\\
0&2\end{matrix}\right)$}\right)+\omega(v,w)\sum_{\gamma\ne
1,\zeta}c_\gamma \gamma.
$$
We have
$$
 \frac{2}{|\Gamma|}\left( {\mbox{\small$\left(\begin{matrix} -\mu-\half & 0\\
0&\mu-\half \end{matrix}\right)$}} 
+
{\mbox{\small$\left(\begin{matrix} 0& 0\\
0&2\end{matrix}\right)$}}   \right) 
= \frac{2}{|\Gamma|} \left(
{\mbox{\small$\left(\begin{matrix} \half & 0\\
0&\half \end{matrix}\right)$}} 
  - (\mu + 1) {\mbox{\small$\left(\begin{matrix} 1& 0\\
0& -1\end{matrix}\right)$}}   \right)
= \frac{1}{|\Gamma|}{ \mbox{\small$\left(\begin{matrix} 1& 0\\
0&1\end{matrix}\right)$}} 
+  c_\zeta {\mbox{\small$\left(\begin{matrix} 1& 0\\
0& -1\end{matrix}\right)$}}
$$
where in the last equality we have used the definition of
$\mu$ from  \eqref{mu}. We find
\beq{qq}
[\theta^{\text{Dunkl}}(v),\theta^{\text{Dunkl}}(w)]=
\omega(v,w)\left( |\Gamma|^{-1} +
c_\zeta\begin{pmatrix} 1& 0\\
0& -1\end{pmatrix}+\sum_{\gamma\ne
1,\zeta}c_\gamma \gamma\right).
\eeq

In these formulas, the matrix $\text{diag}( 1, -1 )\in \End F$
is viewed as an element of $\D(\LL,F)$.
The image of this element in $\D(\X,p,\eta)\ttimes\Gamma$,
 the equalizer smash product
algebra, equals  $1\rtimes\beps(\zeta)$, by definition.
Hence, from \eqref{qq} we deduce
$$
[\theta^{\text{Dunkl}}(v),\theta^{\text{Dunkl}}(w)]=
 \omega(v,w)(1/|\Gamma|+\sum\nolimits_{\gamma\ne 1}\,c_\gamma \gamma).
$$

This completes the proof of Theorem \ref{mess} in the
 special case $n=1$.\qed

\subsection{The map $\Th^{\text{Dunkl}}$.}
In \S\ref{xxx},
for any $\bbe$ and $\bbe'$ and $\btt\in\TT$,  we have defined 
the {\em adjoint} action 
$\Ad_{_F}(\btt): \Hom(F_{\bbe},
F_{\bbe'})\to\Hom(F_{\bbe},
F_{\bbe'}).$
This $\Ad_{_F}$-action of $\TT$ 
 descends to an action on  $\D(\X,p,\eta, F_{\bbe}\to
F_{\bbe'})$, the corresponding quotient space.
The reader should be alerted that
the resulting $\Ad_{_F}$-action on  $\D(\X,p,\eta, F_{\bbe}\to
F_{\bbe'})$ that we are considering at the moment is
{\em different} from the $\Ad_{_{\D\ot F}}$-action of $\TT$
considered in \S\ref{xxx}: the latter action,
comes from the action of $\TT$ on {\em both}
factors of the tensor product
$\D(X,\FF)=\D(X)\ot\End\FF$, while  the former
comes from the action of $\TT$ on the {\em second}
tensor factor, $\End\FF,$ only.

\begin{lemma} \label{ett}
Let $i,j\in I$. For the algebra $\D(\X,p,\eta)\rtimes_\ss \GG$, as defined
in \S\ref{xxx}, one has:
\begin{gather*}
e_{i,n-1}( \D(\X,p,\eta)\rtimes_\ss \GG )e_{j,n-1} 
= \D(\X,p,\varrho, \N_j\to \N_i)^\GG,\\
e_{i,n-1}( \D(\X,p,\eta)\rtimes_\ss \GG )\e 
= \D(\X,p,\varrho, \N_s\to \N_i)^\GG,\\
\e( \D(\X,p,\eta)\rtimes_\ss \GG )e_{j,n-1} 
= \D(\X,p,\varrho, \N_j\to \N_s)^\GG,\\
\e( \D(\X,p,\eta)\rtimes_\ss \GG ) \e 
= \D(\X,p,\varrho, \N_s\to\N_s)^\GG.
\end{gather*}
\end{lemma}

\begin{proof}
We prove the first equality; the rest are similar.

Note that $\C\GG e_{j,n-1} = \bigoplus_{l=1}^n
N_o^{\ot (l-1)} \ot N_j \ot N_o^{\ot (n-l)} = \N^*_j$, so
$$ ( \D(\X,p,\eta)\rtimes_\ss \GG )e_{j,n-1} 
= \bigoplus_{1\leq l\le n} \D(\X,p,\eta)\otimes_\ss (
N_o^{\ot (l-1)} \ot N_j \ot N_o^{\ot (n-l)}), $$
where on the right hand side, $s\in \ss$ acts on $\D(\X,p,\eta)$ 
by right multiplication by $\at(s)$.

For any $\bbe=(\epsilon_1,\ldots,\epsilon_n)\in \mathbb I$,
we write $s^\bbe$ for the character of $\ss$ whose
value at $\zeta_{( l)}$ is $(-1)^{\epsilon_ l}$.

Suppose $\ss$ acts on $N_o^{\ot (l-1)} \ot N_j \ot N_o^{\ot (n-l)}$
by the character $s^{\bbe(l)}$, where $\bbe(l)\in \mathbb I$.
If $j$ is a sink in $Q$, then
$\bbe(l) = (0, ..., 0)$, while if $j$ is a source in $Q$, then
$\bbe(l)=(0,..., 1, ...,0)$ (where the $1$
is in the $l$-th position). 
Under the above right action of $\ss$ on $\D(\X,p,\eta)$
via $\at$, 
the $s^{\bbe(l)}$-isotypic component of $ \D(\X,p,\eta)$
is $\oplus_{\bbe\in\mathbb I}
\D(\X,p,\eta, F_{\bbe(l)} \to F_{\bbe})$, so
$$
\D(\X,p,\eta)\otimes_\ss (
N_o^{\ot (l-1)} \ot N_j \ot N_o^{\ot (n-l)})
= \bigoplus\nolimits_{\bbe\in\mathbb I}\,
\D(\X,p,\eta, F_{\bbe(l)} \to F_{\bbe}) \ot N_j.
$$

Now, the space $e_{i,n-1} (\D(\X,p,\eta)\rtimes_\ss \GG) e_{j,n-1}$
can be written as
$$
\left(\bigoplus_{1\leq m,l\leq n}
(N_o^{\ot (m-1)} \ot N^*_i \ot N_o^{\ot (n-m)})
\bigotimes \left( \oplus_{\bbe\in\mathbb I}
\D(\X,p,\eta, F_{\bbe(l)} \to F_{\bbe}) \ot N_j \right)\right)^\GG.
$$
The subgroup $\ss\subset \GG$ acts on 
$N_o^{\ot (m-1)} \ot N^*_i \ot N_o^{\ot (n-m)}$ by the 
character $s^{\bbe'(m)}$, where
$\bbe'(m)$ is $(0,...,0)$ if $i$ is a sink, 
or $(0,..,1,..,0)$ (where the $1$ is in $m$-th position)
if $i$ is a source.

Recall the $\Ad_{_F}$-action of $\TT$ on 
$\D(\X,p,\eta, F_{\bbe(l)} \to F_{\bbe})$ 
described before this Lemma.
We have the group imbedding $\bept: \ss \into \TT$, see \S\S\ref{equal},\ref{xxx}. 
Equation \eqref{kkk} implies that
the restriction, via $\bept$, of  the
$\Ad_{_F}$-action of $\TT$ to the subgroup $\ss$ coincides with  the
restriction,
via $\beps$, of
the action of
$\GG$ on $\D(\X,p,\eta, F_{\bbe(l)} \to F_{\bbe})$
to its subgroup $\ss$. We have
\begin{align*}
\left((N_o^{\ot (m-1)} \ot N^*_i \ot N_o^{\ot (n-m)})\right.
&\left.\bigotimes \left( \oplus_{\bbe\in\mathbb I}
\D(\X,p,\eta, F_{\bbe(l)} \to F_{\bbe}) \ot N_j \right)\right)^\ss \\
&=
 N^*_i \ot \D(\X,p,\eta, F_{\bbe(l)}\to 
F_{\bbe'(m)}) \ot N_j 
\end{align*}
We conclude that
\begin{align*}
 e_{i,n-1} (\D(\X,p,\eta)\rtimes_\ss \GG) e_{j,n-1} 
=& \left(\bigoplus_{1\leq m,l\leq n}
N_i^*\ot\D(\X,p,\eta, F_{\bbe(l)}\to 
F_{\bbe'(m)}) \ot N_j\right)^\GG.
\end{align*}
The last expression is equal to $\D(\X, p, \varrho, \N_j\to \N_i)^\GG.$
\end{proof}

Recall the homomorphism $\theta^{\mathrm{Dunkl}}$ of
Theorem \ref{mess}.
For any $i,j\in \Ic$, recall the subspace
$\mathsf B_{i,j}$ of $\hh$ defined in (\ref{MM}).
Using Lemma \ref{ett}, we obtain by restricting 
$\theta^{\mathrm{Dunkl}}$ to $\mathsf B_{i,j}$,
a homomorphism
$$\Theta^{\mathrm{Dunkl}}: \mathsf{B}_{i,j}
\too \D(\X,p,\varrho, \N_j\to \N_i)^\GG 
\subset \D(\X,p,\varrho)^\GG.$$
We define the following algebra homomorphism
$$\Theta^{\mathrm{Dunkl}}:\
 \mathsf{B} \to \D(\X,p,\varrho)^\GG,
\quad\sum\nolimits_{i,j}\, u_{i,j}\mto
 \sum\nolimits_{i,j}\, \Theta^{\mathrm{Dunkl}} (u_{i,j}),\quad\forall u_{i,j}
\in\mathsf B_{i,j},\,i,j\in \Ic.
$$

\subsection{Computation of
$\mbox{$\Th^{{\text{Dunkl}}}\ccirc\th^{\text{Quiver}}$}$.}\label{yyy}
Let  $a$ be an edge
of $\overline{Q_{\text{CM}}}$, viewed as an
element of the algebra $\Pi'$.
We would like to 
 compute  $D^a:=\Th^{\text{Dunkl}}\ccirc\th^{\text{Quiver}}(a)\in
\D(\X,p,\varrho)^\GG$,  the image of that element under
the composite map $\Th^{\text{Dunkl}}\ccirc\th^{\text{Quiver}}$.

 We will freely use the notation from \S\ref{formulas}.

\begin{proposition}\label{Da_rk1}
If $a\in Q$ then
$
D^a=-\phi^n_a,
$ and $D^{a^*}$ is an $n\times n$-matrix with the  entries
\begin{gather*}
(D^{a^*})_{mm}=2|\Gamma|^{-1}\frac{\dd}{(\dd\phi_{a^*})_{m}}+\sum_{\gamma\ne
1,\zeta}c_\gamma\frac{(\phi_{a^*}\ccirc(\gamma^{-1}+\id))_m}{\omg;m,m)}
\gamma^{-1},\quad\text{and}\\
(D^{a^*})_{ml}=-\frac{k}{2}\sum_{\gamma\in\Gamma}\frac{(\phi_{a^*}\ccirc\gamma)_l
}{\omg;m,l)}\gamma,\quad\text{for}\en
l\ne m.\end{gather*}
\end{proposition}

\begin{proof}[Proof of Proposition \ref{Da_rk1} for $n=1$.]
In this special case, we have $N_i=\CC\Gamma e_i$ and
$N_i^*=e_i\CC\Gamma$, with the pairing defined by
$$ (e_i\gamma,\gamma' e_i)=e_i\gamma\gamma' e_i\in\CC.$$

Thus, for $n=1$ the formulas of Proposition \ref{Da_rk1} read
\begin{equation*}
D^a= -\phi_a,\quad
D^{a^*}=2|\Gamma|^{-1}\frac{\dd}{\dd \phi_{a^*}}+\sum_{\gamma\ne
1,\zeta}c_\ga\frac{\phi_{a^*}\ccirc(\ga^{-1}+\id)}{\om^\ga}\ga^{-1},\quad\forall a\in Q.
\end{equation*}

To verify these formulas,
 we  write the Dunkl map in the form
$\Theta^{\text{Dunkl}}(v)=\sum_{\gamma\in\Gamma}d_\gamma(v)\gamma,$
where $v\in L$ and $d_\gamma(v)\in\D(L_{\text{reg}},F)$. We  prove
the formula for $D^{a^*}$ because the formula for $D^a$ is easier
to prove.

We recall the construction of the map $\theta^{\text{Quiver}}$.
Let $\tilde{\phi}_{a^*}\in\Hom_\Gamma(N_{h(a^*)},N_{t(a^*)}\otimes L)$
be the element 
corresponding to $\phi_{a^*}\in
\Hom_\Gamma(N^*_{t(a^*)},N^*_{h(a^*)}\otimes L)$. Then the map
$$\Hom_\Gamma(N_{h(a^*)},N_{t(a^*)}\otimes L)\simeq
\Hom_\Gamma(N^*_{t(a^*)},N^*_{h(a^*)}\otimes L)\to
e_{h(a^*)}{\mathbb C}[\Gamma]\otimes L e_{t(a^*)}$$
 from the
construction of $\Theta^{\text{Quiver}}$ is defined by
$\dis\tilde{\phi}_{a^*}\mapsto \tilde{\phi}_{a^*}(1\cdot e_{h(a^*)}).$

Choose a basis $v_1,v_2$ in $L$.
 Then
$\tilde{\phi}_{a^*}=\sum_{s=1}^2 \tilde{\phi}_{a^*}^s\otimes v_s$
where $\tilde\phi_{a^*}^s\in \Hom(N_{h(a^*)},N_{t(a^*)})$. We
consider $D^{a^*}$ as element of $\D(\LL, N_{h(a^*)}\to
N_{t(a^*)})^\Gamma$. From the construction of
$\theta^{\text{Quiver}}$ we find
\begin{align*}
D^{a^*}=&
\sum_{s=1}^2\Theta^{\text{Dunkl}}(v_s)\tilde{\phi}^s_{a^*} =
\sum_{s=1}^2\sum_{\gamma\in
\Gamma}d_\gamma(v_s)(\gamma\ccirc\tilde{\phi}^s_{a^*}) \\
=& 2|\Gamma|^{-1}\frac{\partial}{\partial \tilde{\phi}_{a^*}}+
\sum_{s=1}^2\sum_{\gamma\ne 1,\zeta}c_\gamma\frac{(\gamma
v_s+v_s)^\vee}{\om^\ga}(\gamma\ccirc\tilde{\phi}^s_{a^*}) \\
=& 2|\Gamma|^{-1}\frac{\partial}{\partial \tilde{\phi}_{a^*}}+
\sum_{\gamma\ne
1,\zeta}c_\gamma\frac{\gamma\ccirc\tilde{\phi}_{a^*}\ccirc(\gamma^{-1}+\id)}
{\om^\ga}.
\end{align*}
We have natural isomorphisms
$\Hom(N_{h(a^*)}, N_{t(a^*)})\cong N_{h(a^*)}^*\ot N_{t(a^*)}\cong
\Hom(N^*_{t(a^*)},N^*_{h(a^*)}).$ We deduce
$$\D(\LL,
N_{h(a^*)}\to N_{t(a^*)})^\Gamma\simeq \D(L_{\text{reg}},
N^*_{t(a^*)}\to N^*_{h(a^*)})^\Gamma.$$
Under this isomorphism, the element
$\gamma\tilde{\phi}_{a^*}\ccirc(\gamma^{-1}+\id)$ corresponds to
${\phi}_{a^*}\ccirc(\gamma^{-1}+\id)\ccirc\gamma^{-1}$ and
$\frac{\partial}{\partial \tilde{\phi}_{a^*}}$ corresponds to
$\frac{\partial}{\partial \phi_{a^*}}$.
This
 completes the proof.
\end{proof}

We omit 
the proof of Proposition \ref{Da_rk1} for $n>1$;
 it is  similar to the above computation in the case $n=1$.

It is easy to see that for the edge $b$: $s\to o$, we have
$D^b:=\Th^{{\text{Dunkl}}}\ccirc\th^{\text{Quiver}}(b)=-(1,\dots,1)^t$
and
$D^{b^*}:=\Th^{{\text{Dunkl}}}\ccirc\th^{\text{Quiver}}(b^*)=\nu\cdot(1,\dots,1)$.
Thus, for all $a\in\QQc$,  we have computed the  operators
$D^a:=\Th^{{\text{Dunkl}}}\ccirc\th^{\text{Quiver}}(a)$ where
$D^a\in \D(\X,p,\varrho)$.


\begin{theorem}   \label{diagramcommutes}
For all values of $c,k$,  we have
$\dis\en
\Theta^{\opp{Radial}}\ccirc
\th^{\opp{Holland}}=\Theta^{{\opp{Dunkl}}}\ccirc\th^{\opp{Quiver}}.
$\qed
\end{theorem}

In the special case $n=1$, for any edge $a\in Q$, we have
$$\th^{\text{Radial}}\theta^{\text{Holland}}(a^*)=
2|\Gamma|^{-1}\frac{\partial}{\partial\phi_{a^*}} + \sum_{\ga \neq
1,\zeta} c_\ga \frac{\phi_{a^*}\ccirc(\ga^{-1}+\id)}{\om^\ga}
\ga^{-1} - |\Gamma|^{-1} \sum_{\ga \neq 1,\zeta}
\frac{\phi_{a^*}\ccirc\ga^{-1}}{\om^\ga}.
$$
Therefore, replacing here the map $\th^{\text{Radial}}$ by 
$\Th^{\text{Radial}},$ we find
$$
\Theta^{\text{Radial}} \theta^{\text{Holland}}(a^*) =
2|\Gamma|^{-1}\frac{\partial}{\partial\phi_{a^*}} + \sum_{\ga \neq
1,\zeta} c_\ga \frac{\phi_{a^*}\ccirc(\ga^{-1}+\id)} {\om^\ga}\ga^{-1}  .$$
When $n>1$, it is completely similar.


\section{Harish-Chandra homomorphism}

Recall that we assume $\lambda\cdot\delta=1$,
i.e. $t=|\Gamma|^{-1}$.
We shall write $\hh_{k,c}$ for $\hh_{t,k,c}(\GG)$.

\subsection{Modified Holland's map.}
In this subsection, we define a  map
$\Theta^{\text{Holland}}: e_s\Pi' e_s \to \A_{\chi'}$,
cf. \eqref{achi}.
To this end, assume for the moment that $\nu$ is a formal variable,
and the algebras $\Pi$, $\Pi'$, $\T_\chi$, $\A_{\chi'}$
are all defined over $\C[\nu]$.

\begin{lemma} \label{radialinjective}
The map $\gr(\theta^{\opp{Radial}}): \gr\A_{\chi'} \too
\gr(\D( \X,p,\rho_s )^\GG \ot \C[\nu])$ is injective.
\end{lemma}
\begin{proof}
This follows from Proposition \ref{ppp},
Theorem \ref{idealreduced} and Proposition \ref{gris}.
\end{proof}

The preceding lemma implies that $\gr\A_{\chi'}$ and
$\A_{\chi'}$ are free $\C[\nu]$-modules.

We define a homomorphism
$$\tilde\Theta^{\text{Holland}}:
\C\QQc \ot \C[\nu] \too \big( \D(\Qc,\al) \ot (\End\N)
\ot \C[\nu, \nu^{-1}] \big)^{\GL(\al)}$$
by
\begin{align*}
\tilde\Theta^{\text{Holland}}(e_j) =& \tilde\th^{\text{Holland}}(e_j)
\quad \text{ for all vertices } j,\\
\tilde\Theta^{\text{Holland}}(a) =&
\tilde\th^{\text{Holland}}(a) \quad \text{ for any edge } a\neq b^*,\\
\tilde\Theta^{\text{Holland}}(b^*) =&
\nu^{-1} \tilde\th^{\text{Holland}}(b^*).
\end{align*}
It is easy to see that since $\tilde\th^{\text{Holland}}$ descends
to a homomorphism $\th^{\text{Holland}}: \Pi \to \T_\chi$,
the homomorphism $\tilde\Theta^{\text{Holland}}$ descends to a homomorphism
$$\Theta^{\text{Holland}}: \Pi' \to \T_\chi [\nu^{-1}],
\quad\text{such that}\quad
\Theta^{\text{Holland}}: e_s\Pi' e_s \to\A_{\chi'}[\nu^{-1}].
$$

Suppose that $X\in\C\QQ$. By (\ref{radialbb}),
$\th^{\opp{Radial}}\th^{\text{Holland}}(b^*Xb)$ vanishes if we set $\nu=0$,
so by Lemma \ref{radialinjective},
$$ \tilde\th^{\text{Holland}}(b^*Xb) \in
( \D(\Qc,\al)(\bxi - (\la-\partial-n e_s))
(\gl(\al)) )^{\GL(\al)}. $$
Since $\chi'= (\la -\partial- n e_s) - \nu e_o + n\nu e_s$, we have
that $\tilde\th^{\text{Holland}}(b^*Xb)$ belongs to
$$
 \D(\Qc,\al)(\bxi - \chi')(\gl(\al))
 + \nu \cd\D(\Qc,\al) (e_o -n e_s)(\gl(\al)).
$$
Since $\GL(\al)$ is a reductive group, we have
a projection map
$$\pr:\D(\Qc,\al) \to \D(\Qc,\al)^{\GL(\al)}$$
such that
$$\pr\big(\D(\Qc,\al)(\bxi - \chi')(\gl(\al))\big)
\sset ( \D(\Qc,\al)(\bxi - \chi')(\gl(\al)) )^{\GL(\al)}.$$
Thus, the element
$\tilde\th^{\text{Holland}}(b^*Xb) =
\pr (\tilde\th^{\text{Holland}}(b^*Xb))$
belongs to
$$
\big( \D(\Qc,\al)(\bxi - \chi') (\gl(\al)) \big)^{\GL(\al)}
+ \nu\cd \D(\Qc,\al)^{\GL(\al)}.
$$
Therefore, $\Th^{\text{Holland}}(e_s\Pi' e_s) \subseteq \A_{\chi'}$.
Thus, for any $\nu\in\C$, we have a homomorphism
$$\Th^{\text{Holland}}: e_s\Pi' e_s \too \A_{\chi'}.$$

\begin{theorem} \label{diagramcommutes2}
The following diagram commutes:
$$
\xymatrix{
e_s\Pi'e_s
\ar[rr]^{\Theta^{\opp{Quiver}}}
\ar[d]_{\Theta^{\opp{Holland}}} && \e\hh_{k,c}\e
\ar[d]^{\Theta^{{\opp{Dunkl}}}} \\
\A_{\chi'}
\ar[rr]^{\Theta^{\opp{Radial}} }
&& \D(\X,p,\varrho_s)^\GG
}
$$
\end{theorem}
\begin{proof}
This follows from Theorem \ref{diagramcommutes}.
\end{proof}

\begin{proposition}  \label{surjective}
The map $\Theta^{\opp{Holland}}: e_s\Pi'e_s \too \A_{\chi'}$
is surjective.
\end{proposition}
\begin{proof}
The algebra $\Pi'$ has a filtration with
$\deg(a)=1$ for edges $a\neq b,b^*$, and
$\deg(b)=\deg(b^*)=0$.
It suffices to show that the associated graded map
$$\gr(\Theta^{\text{Holland}}): e_s\gr(\Pi')e_s \too
\gr(\A_{\chi'})$$ is surjective.

Now $e_s\gr(\Pi')e_s$ is generated by
$e_s$ and  $b^*(e_o\Pi(Q)e_o)b$
where $\Pi(Q)$ is the preprojective algebra of $Q$.
By \cite[Theorem 3.4]{CB}, we have
$$ \C[\Rep_{n\delta}(\Pi(Q))]^{GL(n\delta)}
= \left( (\C[\Rep_{\delta}(\Pi(Q))]^{\GL(\delta)})^{\ot n}
\right)^{S_n}. $$
Denote by $$\Theta^{\text{CBH}}: \Pi(Q) \too
\C[\Rep_{\delta}(\Pi(Q))\ot
\End(\oplus_{i\in I}\C^{\delta_i})]^{\GL(\delta)}$$
the natural morphism defined in \cite{CBH}, which gives a morphism
$$\Theta^{\text{CBH}}: e_o\Pi(Q)e_o \too
\C[\Rep_{\delta}(\Pi(Q))]^{GL(\delta)}.$$
This latter morphism is an isomorphism by \cite[Theorem 8.10]{CBH}.

Now, given an element $X\in e_o\Pi(Q)e_o$, we claim that
\begin{equation} \label{grholl}
\gr(\Th^{\text{Holland}})(b^*Xb) =
 - \sum\nolimits_{p=1}^n\, 1^{\ot (p-1)}\ot\Th^{\text{CBH}}(X)\ot 1^{\ot (n-p)}.
\end{equation}
Indeed, we have
\begin{align*}
 \gr(\th^{\text{Radial}})&\gr(\Th^{\text{Holland}})(b^*Xb) =
 \gr(\Th^{{\text{Dunkl}}})\gr(\Th^{\text{Quiver}})(b^*Xb) \qquad \text{ (by
Theorem \ref{diagramcommutes2}) } \\
=& \gr(\th^{\text{Radial}}) \left(- \sum_{p=1}^n 1^{\ot
(p-1)}\ot\Th^{\text{CBH}}(X)\ot 1^{\ot (n-p)} \right)
\qquad \quad\text{ (by (\ref{grg}) below). }
\end{align*}
Hence, (\ref{grholl}) follows from injectivity of $\gr(\th^{\text{Radial}})$.
It follows from (\ref{grholl}) and
 Lemma \ref{symmetricinvariants}
below that
$\gr(\Th^{\text{Holland}})$ is surjective.
\end{proof}

\subsection{}\label{hc_sec}
It follows from Theorem \ref{diagramcommutes2},
Proposition \ref{surjective}, and the injectivity
of $\Th^{{\text{Dunkl}}}$ that we have a homomorphism
$$
(\Th^{{\text{Dunkl}}})^{-1}\ccirc\Th^{\text{Radial}} :
\A_{\chi'} \too \e\hh_{k,c}\e.
$$
Since $\Rep_\al(\Qc) = \Rep_{n\delta}(Q) \oplus \C^n$,
we have an obvious embedding
$\varpi:
\D(Q,n\delta) \rightarrow \D(\Qc,\al)$.
\begin{definition}\label{hc_def}
The Harish-Chandra homomorphism $\Phi_{k,c}$ is defined to be the
composition
\begin{equation}  \label{hc}
\xymatrix{ \D(Q,n\delta)^{\GL(n\delta)}
 \ar[r]^{\qquad\quad \varpi} & \A_{\chi'}
\ar[rrr]^{(\Th^{{\opp{Dunkl}}})^{-1}\ccirc\Th^{\opp{Radial}}} &&&
\e\hh_{k,c}\e .}
\end{equation}
\end{definition}

Following \cite{EG}, we define a $1$-parameter space
of representations $V_d$ of $\gl_n$ as follows.
As a vector space, $V_d$ is spanned by expressions
$(x_1\cdots x_n)^d\cdot P$, where $P$ is a Laurent polynomial in
$x_1,\ldots,x_n$ of total degree $0$. The Lie algebra
$\gl_n$ has an action on $V_d$ by formal differentiation,
where $e_{p,q}$ acts by $x_p\frac{\partial}{\partial x_q}$.
We restrict this to an $\sln$ action.
The desired $\gl_n$ action on $V_d$ is
obtained by pulling back the $\sln$ action via the natural Lie
algebra projection $\gl_n \to \sln$, so that the center
of $\gl_n$ acts trivially.

Let $Fun(...)$ denote the vector space
of functions on a formal neighborhood
of a point of the slice $\lreg$.
Recall that
$\partial = -n \sum_{i\in I}\Tr|_{N_i}(\zeta) e_i$.
We have
\begin{equation}  \label{qqqq}
\left( Fun(\Rep_\al(\Qc)) \otimes \C_{-\chi'}
\right)^{\mathfrak{gl}(\al)}
= \left( Fun(\Rep_{n\delta}(Q))
\widehat\otimes V_{\nu-1} \otimes \C_{-\lambda+e_o+\partial}
\right)^{\mathfrak{gl}(n\delta)}.
\end{equation}

The $\GL(n\delta)$ action on $\Rep_{n\delta}(Q)$
induces a Lie algebra map
$\ad: \gl(n\delta) \to \D(Q,n\delta)$.
Let $\ad: U\gl(n\delta) \to \D(Q,n\delta)$
be the induced map on the universal enveloping algebra
of $\gl(n\delta)$.
Define the left ideal
$$ J_{k,c} := \D(Q,n\delta)\cd \ad(
\Ann(V_{\nu-1} \otimes \C_{-\lambda+e_o+\partial}) ) \sset\D(Q,n\delta).$$
By (\ref{qqqq}),
the ideal $J_{k,c}^{\GL(n\delta)}$ is in the kernel of the
map $\varpi$ in (\ref{hc}).

\begin{theorem} \label{maintheorem}
The Harish-Chandra homomorphism induces  an algebra isomorphism
$$
 \Phi_{k,c}:  \D(Q,n\delta)^{\GL(n\delta)}
\big/J_{k,c}^{\GL(n\delta)}\iso\e\hh_{k,c}\e.
$$
\end{theorem}
\begin{proof}
By Theorem \ref{idealreduced} and \cite[Theorem 1.3]{EG},
the associated graded map, $\gr\Phi_{k,c},$
 is the
isomorphism in (\ref{gradedmapred}), hence $\Phi_{k,c}$ is itself an
isomorphism.
\end{proof}

\subsection{Proof of Corollary \ref{maincorollary}.} \label{fisc}
Given any $C=\sum_{\ga\in\Gamma} C_\ga \ga\in \C[\Gamma]$,
we let $\overline C = \sum_{\ga\in\Gamma} C_\ga \ga^{-1}$.
Correspondingly, if $\la= \sum_{i\in I}\Tr|_{N_i}(C)e_i$,
then let $\overline\la= \sum_{i\in I}\Tr|_{N_i}(\overline C)e_i$.
We have an anti-isomorphism
\begin{equation}  \label{antiisomorphism}
\hh_{k,c} \iso \hh_{k,\overline c},\quad
 g\mapsto g^{-1}, \quad u \mapsto \sqrt{-1}u, \qquad
\forall g\in \GG,\ u\in L^n.\end{equation}
We also have an isomorphism
\begin{equation} \label{isomorphism}
\hh_{k,c} \iso \hh_{-k,c},\quad
 \sigma \mapsto (-1)^\sigma\sigma, \quad g\mapsto g,\quad
u \mapsto u, \qquad\forall \sigma\in S_n,\
g\in\Gamma^n,\ u\in L^n. \end{equation}
The isomorphism in (\ref{isomorphism}) sends $\e$ to $ \e_-.$

Now, for any $i\in I$ set $\la_i^\dag := \Tr|_{N_i} (t\cdot 1 +
c^\dag).$
We put
$$
c^\dag := -c + 2|\Gamma|^{-1}\sum\nolimits_{\ga\neq 1}\,\ga,
\quad\text{and}\quad
\la^\dag := \sum\nolimits_{i\in I}\, \la_i^\dag e_i
= -\la + 2e_o.
$$

The group $\GL(n\delta)$ acts on $\det(\Rep_{n\delta}(Q)^*)$
by the character $2\partial$.
We have
\begin{align*}
( V_{\nu-1} \ot \C_{-\la+e_o+\partial} )^*
\ot \det(\Rep_{n\delta}(Q)^*)
\simeq & V_{-\nu} \ot \C_{\la-e_o-\partial}
\ot \C_{2\partial} \\
= & V_{-\nu} \ot \C_{-\la^\dag+e_o+\partial}.
\end{align*}

Let $\mathsf i: \D(Q,n\delta) \to \D(Q,n\delta)$
be the anti-isomorphism sending a differential operator to
its adjoint.
Then for any $\GL(n\delta)$-module $V$,
$$\mathsf i(\ad (\Ann(V)))
= \ad( (\Ann( V^* \ot \det(\Rep_{n\delta}(Q)^*) )) ).$$

The proof of the first isomorphism in 
Corollary \ref{maincorollary}
is now completed by the following isomorphisms
\begin{align*}
\e\hh_{k,c}\e \simeq & (\e\hh_{ 2|\Gamma|^{-1}-k,c^\dag } \e)^{op}
\qquad  \text{ using  Theorem \ref{maintheorem} } \\
\simeq & \e\hh_{ 2|\Gamma|^{-1}-k,\overline{c^\dag} } \e
\qquad\quad\en  \text{ using  (\ref{antiisomorphism}) } \\
\simeq & \e_-\hh_{ k-2|\Gamma|^{-1},\overline{c^\dag} } \e_-
\qquad  \text{ using  (\ref{isomorphism}). }
\end{align*}

We will prove the second isomorphism in 
Corollary \ref{maincorollary} later in \S\ref{shiftsubsec}.

\section{Reflection isomorphisms}\label{refl_sec}

Except for \S\ref{shiftsubsec},
this section is independent of the earlier sections.

\subsection{}\label{refl1}
Let $Q$ be an arbitrary quiver (not necessarily of
type $\widetilde A$, $\widetilde D$, or $\widetilde E$).
Denote by $I$ the set of vertices of $Q$.
Let $R= \bigoplus_{i\in I}\C$, and $E$ the vector space
with basis formed by the set of edges $\{a\in \QQ\}$.
Thus, $E$ is naturally a $R$-bimodule.
The path algebra of $\QQ$
is $\C\QQ := T_R E = \bigoplus_{n\geq 0} T^n_R E$, where
$T^n_R E = E\ot_R \cdots \ot_R E$ is the $n$-fold
tensor product.
The trivial path for the vertex $i$ is denoted by
$e_i$, an idempotent in $R$.

Fix a positive integer $n$. Let $\sfR := R^{\ot n}$.
For any $\ell\in [1,n]$, define the $\sfE$-bimodules
$$ \sfE_{\ell}:= R^{\ot (\ell-1)}\ot E\ot R^{\ot (n-\ell)}
\qquad \mathrm{and}\qquad \sfE :=
\bigoplus\nolimits_{1\leq\ell\leq n}\, \sfE_\ell\,.$$
The natural inclusion $\sfE_\ell \hookrightarrow
R^{\ot (\ell-1)}\ot T_R E\ot R^{\ot (n-\ell)} $
induces a canonical identification
$T_\sfR\sfE_\ell = R^{\ot (\ell-1)}\ot T_R E\ot R^{\ot (n-\ell)}$.
Given two elements $\varepsilon\in\sfE_\ell$
and $\varepsilon'\in\sfE_m$ of the form
\begin{equation} \label{eqnve}
\varepsilon= e_{i_1}\ot e_{i_2}\otot a\otot h(b) \otot e_{i_n}\,,
\end{equation}
\begin{equation} \label{eqnve'}
\varepsilon'= e_{i_1}\ot e_{i_2}\otot t(a)\otot b\otot e_{i_n}\,,
\end{equation}
where $\ell\neq m$, $a,b\in \QQ$ and $i_1, \ldots, i_n \in I$,
we define
\begin{eqnarray}
\lfloor\varepsilon,\varepsilon'\rfloor &:=&
(e_{i_1}\otot a\otot h(b) \otot e_{i_n})
(e_{i_1}\otot t(a)\otot b \otot e_{i_n}) \nonumber\\ & &
-(e_{i_1}\otot h(a)\otot b \otot e_{i_n})
(e_{i_1}\otot a\otot t(b) \otot e_{i_n}) . \nonumber
\end{eqnarray}
Note that $\lfloor\varepsilon,\varepsilon'\rfloor$
is an element in $T^2_\sfR \sfE$.

\begin{definition}[\cite{GG}, Def. 1.2.3] \label{dea}
For any $\la=\sum_{i\in I} \la_i e_i$
where $\la_i\in \C$,
and $\nu\in \C$, define the algebra
$\sfA_{n,\la,\nu}(Q)$ to be the quotient of
$T_{\sfR}\sfE\rtimes\C[S_n]$
by the following relations.
\begin{itemize}
\item[\vi]
For any $i_1,\ldots, i_n\in I$ and $\ell\in [1,n]$:
\begin{gather*}
e_{i_1}\otot \left(
\sum_{\{a\in Q\,|\, h(a)=i_\ell\}}
a\cdot a^*- \sum_{\{a\in Q\,|\, t(a)=i_\ell\}} a^*\cdot a
-\la_{i_\ell}e_{i_\ell} \right)
\otot e_{i_n}  \\
= \nu \sum_{\{ j\neq\ell \,|\, i_j=i_\ell\}}
(e_{i_1}\otot e_{i_\ell} \otot e_{i_n})s_{j\ell}.
\end{gather*}
\item[\vii]
For any $\varepsilon,\varepsilon'$
of the form (\ref{eqnve})--(\ref{eqnve'}):
\[  \lfloor \varepsilon,\varepsilon'\rfloor
= \left\{ \begin{array}{ll}
\nu\cdot(e_{i_1}\otot h(a)\otot t(a)\otot e_{i_n})s_{\ell m}
& \textrm{if $b\in Q$,\ $a = b^*$} ,\\
- \nu\cdot(e_{i_1}\otot h(a)\otot t(a)\otot e_{i_n})s_{\ell m}
& \textrm{if $a\in Q$,\ $b = a^*$} ,\\
0 & \textrm{else}\,.  \end{array} \right.  \]
\end{itemize}
\end{definition}

When $n=1$, there is no parameter $\nu$, and $\sfA_{n,\la,\nu}(Q)$
is the deformed preprojective algebra $\Pi_{\la}(Q)$ defined
in \cite{CBH}.

\subsection{Quiver functors.}  \label{qf}
The goal of this section is to put the
construction of the functor $M\to \wt M$ exploited in \S\ref{quivermap}
into an appropriate, more general, context.

Let $T$ be a nonempty subset of $I$, and let
$e_T := \sum_{i\in T} e_i$. In particular, $e_I = 1$.
Let $Q_T$ be a quiver obtained from $Q$ by adjoining a vertex $s$,
and arrows $b_i : s \to i$ for $i\in T$.
We call $s$ the \emph{special} vertex.
We shall define a functor $G$ from $\sfA_{n,\la,\nu}(Q)$-modules
to $\Pi_{\la-\nu e_T +n\nu e_s}(Q_T)$-modules.

Let $M$ be an $\sfA_{n,\la,\nu}(Q)$-module.
We want to define a $\Pi_{\la-\nu e_T +n\nu e_s}(Q_T)$-module
$G(M)$. For each $i\in I$, let
$G(M)_i := e_{n-1}(e_i\ot e_T^{\ot (n-1)})M$.
Also, let $G(M)_s := e_n e_T^{\ot n}  M$.
If $a$ is an edge in $\QQ$, then define
$a: G(M)_{t(a)} \too G(M)_{h(a)}$
to be the map given by the element
$a\ot e_T^{\ot (n-1)} \in \sfA_{n,\la,\nu}(Q)$.
We have an inclusion $G(M)_s \subset
e_{n-1} e_T^{\ot n} M = \bigoplus_{j\in T} G(M)_j$.
For $i\in T$, we have a projection map
$\mathsf{pr}_i: \bigoplus_{j\in T} G(M)_j \too G(M)_i$.
Define $b_i: G(M)_s \too G(M)_i$ to be the restriction
of $\mathsf{pr}_i$ to $G(M)_s$.
Define $b^*_i: G(M)_i \too G(M)_s$ to be
$-\nu\cdot(1+s_{12}+\cdots + s_{1n})$.

The following lemma is a generalization of Lemma \ref{pimodule}.
\begin{lemma}
With the above actions,  $G(M)$ is a
$\Pi_{\la-\nu e_T +n\nu e_s}(Q_T)$-module.
\end{lemma}
\begin{proof}
It is clear that $(1+s_{12}+\cdots+s_{1n})e_{n-1}= ne_n$.

On $G(M)$, at the special vertex $s$,
we have
$\sum_{i\in T} b^*_ib_i = -n\nu.$

At a vertex $i\in I$, $i\notin T$, we have
$$ \sum_{a\in Q; h(a)=i} aa^* - \sum_{a\in Q; t(a)=i} a^*a
= \la_i $$
by the relation \vi in Definition \ref{dea}.
At a vertex $i\in T$, we have
\begin{align*}
\sum_{a\in Q; h(a)=i} aa^* - \sum_{a\in Q; t(a)=i} a^*a
=& \la_i + \nu\cd \mathsf{pr}_i(s_{12} + \cdots + s_{1n})
= \la_i - \nu -  b_ib_i^*,
\end{align*}
using again the relation \vi in Definition \ref{dea}.
\end{proof}

It is clear that the assigment $M\mapsto G(M)$ is functorial.
We have constructed a functor
\begin{equation}
G:\Mod{\sfA_{n,\la,\nu}(Q)}
\too \Mod{\Pi_{\la-\nu e_T +n\nu e_s}(Q_T)}.
\end{equation}

Recall the symmetrizer $e_n:=\frac{1}{n!}\sum_{s\in S_n} s.$

\begin{definition}
Let $\UU_{n,\la,\nu}(Q):=e_n\sfA_{n,\la,\nu}(Q)e_n$
be the {\em spherical subalgebra} in $\sfA_{n,\la,\nu}(Q).$
\end{definition}

The idempotents $e_n$ and
$\eel:=e_T^{\ot n}$  commute.
For  $M := \sfA_{n,\la,\nu}(Q) e_n\eel,$
we get
$$G(M)_s =  \eel\UU_{n,\la,\nu}\eel $$
In this case, $G(M)_s$ is an algebra, and the action of
$e_s \Pi_{\la-\nu e_T +n\nu e_s}(Q_T) e_s$ on $G(M)_s$
commutes with right multiplication by the elements of $G(M)_s$.
Thus,  our construction yields an algebra homomorphism
\begin{equation}
\widehat G: e_s \Pi_{\la-\nu e_T +n\nu e_s}(Q_T) e_s  \too
 \eel\UU_{n,\la,\nu}\eel .
\end{equation}

\subsection{Modified version.}
The map $\widehat G$ is $0$ on nonconstant paths when $\nu=0$.
For this reason, we shall need a slight modification of the
constructions in the previous subsection.

Define $\Pi'_{\la-\nu e_T +n\nu e_s}(Q_T)$ to be the quotient of
the path algebra $\C\overline{Q_T}$ by the following relations:
$$
\sum_{a\in Q}[a,a^*]+\nu\sum_{i\in T} b_ib_i^*=\la-\nu e_T,
\qquad  \sum_{i\in T} b_i^*b_i = -ne_s.
$$
We have an algebra morphism
$\Pi_{\la-\nu e_T +n\nu e_s}(Q_T) \too
\Pi'_{\la-\nu e_T +n\nu e_s}(Q_T)$
defined on the edges by
$$ a \mapsto a \text{ for } a \neq b_i^*,
\qquad b^*_i \mapsto \nu b_i^*. $$
This is an isomorphism only when $\nu\neq 0$.

Given a $\sfA_{n,\la,\nu}(Q)$-module $M$, we construct a
$\Pi'_{\la-\nu e_T +n\nu e_s}(Q_T)$-module $G'(M)$ analogous to
$G(M)$ in the previous subsection, the only difference is
that now, we let
$b_i^*: G'(M)_i \too G'(M)_s$ be $-(1+s_{12}+\cdots+s_{1n})$.
Hence, as above, we obtain a functor
\begin{equation}
G':\Lmod{\sfA_{n,\la,\nu}(Q)}
\too \Lmod{\Pi'_{\la-\nu e_T +n\nu e_s}(Q_T)}
\end{equation}
as well as a morphism
\begin{equation} \label{mqm}
\widehat G': e_s \Pi'_{\la-\nu e_T +n\nu e_s}(Q_T) e_s  \too
 \eel\UU_{n,\la,\nu}\eel .
\end{equation}

The algebra $\Pi'_{\la-\nu e_T +n\nu e_s}(Q_T)$
has a filtration with $\deg(a)=1$ for $a \neq b_i, b^*_i$,
and $\deg(b_i)=\deg(b^*_i)=0$ for $i\in T$. Also,
$ e_s \gr(\Pi'_{\la-\nu e_T +n\nu e_s}(Q_T)) e_s $ is
generated by $e_s$ and
$\left( \bigoplus_{i,j\in T} b^*_j
\Pi_0 (Q) b_i \right)$.

\emph{We shall assume that $Q$ is a connected quiver without edge-loops,
and $Q$ is not a finite Dynkin quiver.}
Then, by \cite[Theorem 2.2.1]{GG} and \cite[Remark 2.2.6]{GG},
$\gr\sfA_{n,\la,\nu}(Q) = \Pi_0(Q)^{\otimes n}\rtimes \C[S_n]$.
Thus,
$$ \gr\big(\eel\UU_{n,\la,\nu}\eel\big)
= (( e_T \Pi_0(Q) e_T )^{\ot n})^{S_n}. $$
Now given an element $X\in e_j\Pi_0 (Q) e_i$ where $i,j\in T$,
we have
\begin{equation} \label{grg}
\gr(\widehat G')(b^*_jXb_i) = -\sum_{p=1}^n e_T^{\otimes (p-1)}
\otimes X \otimes e_T^{\otimes (n-p)}.
\end{equation}

\begin{lemma}  \label{symmetricinvariants}
Let $A$ be any associative algebra with unit $1\in A$.
Then $(A^{\ot n})^{S_n}$ is generated as an algebra by elements of the form
$$ \sum\nolimits_{p=1}^n\, 1^{\otimes (p-1)}
\otimes X \otimes 1^{\otimes (n-p)}, \qquad X\in A. $$
\end{lemma}
\begin{proof}
Since $(A^{\ot n})^{S_n}$ is spanned by elements of the form
$a^{\ot n}$ where $a\in A$, it suffices to show that the lemma is true
for $A=\C[a]$, but this follows from the
main theorem on symmetric functions.
\end{proof}

\begin{proposition} \label{gsurjective}
The map $\widehat G'$ in (\ref{mqm}) is surjective.
\end{proposition}
\begin{proof}
It suffices to show that $\gr(\widehat G')$ is surjective. This
follows from (\ref{grg}) and the preceding lemma.
\end{proof}

\subsection{Reflection functors.} Recall the setting of
reflection functors as in \eqref{gan}. In particular,
we have
the Weyl group $W$  generated by the
simple reflections $r_i$ for $i\in I$.
We also have a nonempty subset $T\sset I$ and
we fix a vertex $i\notin T$.

Let us apply the reflection functor
$F_i$ to the $\sfA_{n,\la,\nu}(Q)$-module
$\sfA_{n,\la,\nu}(Q) \eel $.
By construction, we have
$$\eel  F_i(\sfA_{n,\la,\nu}(Q) \eel )
=\eel  \sfA_{n,\la,\nu}(Q) \eel $$
and the left action of
$\eel  \sfA_{n,r_i(\la),\nu}(Q) \eel $
on $\eel  F_i(\sfA_{n,\la,\nu}(Q) \eel )$
commutes with the right multiplication by
$\eel  \sfA_{n,\la,\nu}(Q) \eel $.
Hence, for $i\notin T$, we obtain a homomorphism
\begin{equation}   \label{rfisom}
\widehat F_i : \eel  \sfA_{n,r_i(\la),\nu}(Q)
\eel  \too
\eel  \sfA_{n,\la,\nu}(Q) \eel .
\end{equation}
Note that $\widehat F_i(\eel\UU_{n,r_i(\la),\nu}(Q) \eel)
\subset \eel\UU_{n,\la,\nu}\eel $.

In the special case when $n=1$, reflection functors were
constructed in \cite{CBH}; let us recall their definition.
Since $\Pi_\la(Q)$ does not depend on the orientation of $Q$,
we may assume
without loss of generality that $i$ is a sink in $Q$.
Let $M$ be a $\Pi_\la(Q)$-module, and $M_j=e_jM$ for each $j\in I$.
For each edge $a\in Q$ such that $h(a)=i$, write
$$ \pi_a :\bigoplus_{\xi\in Q; h(\xi)=i} M_{t(\xi)} \to M_{t(a)},
\quad \mu_a : M_{t(a)} \to \bigoplus_{\xi\in Q; h(\xi)=i} M_{t(\xi)} $$
for the projection map and inclusion map, respectively.
Define
$$\pi: \bigoplus_{a\in Q; h(a)=i} M_{t(a)} \too M_i,
\quad \pi := \sum_{a\in Q; h(a)=i} a \ccirc \pi_a,$$
and
$$\mu: M_i \too \bigoplus_{a\in Q; h(a)=i} M_{t(a)},
\quad \mu := \sum_{a\in Q; h(a)=i} \mu_a \ccirc a^*.$$
Observe that $\pi\mu=\la_i$.
Let $(F_i(M))_j:= M_j$ if $j\neq i$, and let $(F_i(M))_i :=
\mathrm{Ker}(\pi)$.
If $a\in \overline Q$ and $h(a),t(a)\neq i$, then
let $a: (F_i(M))_{t(a)} \to (F_i(M))_{h(a)}$ be the same map as
$a: M_{t(a)} \to M_{h(a)}$. If $a\in Q$ and $h(a)=i$, then let
$a: (F_i(M))_{t(a)} \to (F_i(M))_i$ be the map
$(-\la_i+\mu\pi)\mu_a$, and let
$a^*: (F_i(M))_i \to (F_i(M))_{t(a)}$ be the map $\pi_a$ restricted
to $(F_i(M))_i$.
Letting $F_i(M)=\oplus_{j\in I}(F_i(M))_j$, we have defined the functor
$$ F_i : \Lmod{\Pi_\la(Q)} \too \Lmod{\Pi_{r_i(\la)}(Q)}
\quad \text{ for any } i\in I.$$
In particular, for the quiver $Q_T$, and for
$i\in I$ but $i\notin T$, we have
\begin{equation} \label{fiqt}
F_i : \Lmod{\Pi_{\la-\nu e_T+ n\nu e_s}(Q_T)} \too
\Lmod{\Pi_{r_i(\la)-\nu e_T+n\nu e_s}(Q_T)}
\end{equation}

Let $i\in I$ but $i\notin T$. We define a functor
\begin{equation*}
F'_i : \Lmod{\Pi'_{\la-\nu e_T+ n\nu e_s}(Q_T)} \too
\Lmod{\Pi'_{r_i(\la)-\nu e_T+n\nu e_s}(Q_T)}
\end{equation*}
in exactly the same way as $F_i$ in (\ref{fiqt}).
It is easy to see from definitions that the diagram
(\ref{diag_commutes}) commutes.

\subsection{Relations in rank 1.}
In this subsection, the rank $n$ is equal to $1$.
Let $\mathfrak C = (\mathfrak C_{ij})$ be the generalized Cartan matrix
of $Q$.

\begin{proposition}  \label{relations1}
For all $\la\in R$, we have the following.

\vi The map $$\widehat F_i: (1-e_i)\Pi_\la(Q)(1-e_i) \to
(1-e_i)\Pi_{r_i(\la)}(Q)(1-e_i)$$ is an isomorphism, and
$\widehat F_i^2 =\id$.

\vii If $\mathfrak C_{ij}=0$, then
$$\widehat F_i \ccirc \widehat F_j = \widehat F_j \ccirc \widehat F_i:
(1-e_i-e_j)\Pi_\la(Q)(1-e_i-e_j) \to
(1-e_i-e_j)\Pi_{r_i r_j(\la)}(Q)(1-e_i-e_j).$$

\viii If $\mathfrak C_{ij}=-1$, then
$$\widehat F_i \ccirc \widehat F_j \ccirc\widehat F_i =
\widehat F_j \ccirc\widehat F_i \ccirc\widehat F_j :
(1-e_i-e_j)\Pi_\la(Q)(1-e_i-e_j) \to
(1-e_i-e_j)\Pi_{r_i r_j r_i(\la)}(Q)(1-e_i-e_j).$$
\end{proposition}

\begin{proof}
\vi The algebra $(1-e_i)\Pi_\la(Q)(1-e_i)$ is generated by
edges $a\in\QQ$ with $h(a),t(a)\neq i$, and paths of length two:
$a_2a_1$ with $h(a_2), t(a_1)\neq i$ and $t(a_2)=h(a_1)=i$.

If $a\in \QQ$ and $h(a),t(a)\neq i$, then $\widehat F_i(a)=a$.

Now let $a_2a_1$ be a path with $h(a_2), t(a_1)\neq i$ and $t(a_2)=h(a_1)=i$.
If $a_2\neq a_1^*$ or $a_1\neq a_2^*$, then $\widehat F_i(a_2a_1) = a_2a_1$.
If $a_2=a_1^*$, then $\widehat F_i(a_2a_1) = -\la_i e_{t(a_1)} + a_2a_1$,
and so
$\widehat F_i(\widehat F_i(a_2a_1)) = -\la_i e_{t(a_1)}
+ \la_i e_{t(a_1)} + a_2a_1 =a_2a_1$.

\vii When $\mathfrak C_{ij}=0$, there is no edge joining $i$ and $j$.
In this case, it is clear that $F_i F_j=F_j F_i$,
so $\widehat F_i \widehat F_j = \widehat F_j \widehat F_i$.

\viii When $\mathfrak C_{ij}=-1$, there is precisely one edge in $Q$
joining $i$ and $j$, say $a:i\to j$.
The algebra $(1-e_i-e_j)\Pi_\la(Q)(1-e_i-e_j)$ is generated by:
\begin{itemize}
\item edges $a_1\in \QQ$ with $h(a_1),t(a_1)\neq i,j$;

\item paths $a_2a_1$ with $t(a_2)=h(a_1)=i$ and $h(a_2),t(a_1)\neq i,j$;

\item paths $a_2a_1$ with $t(a_2)=h(a_1)=j$ and $h(a_2),t(a_1)\neq i,j$;

\item paths $a_3 a_2 a_1$ with $a_2=a$, $t(a_3)=j$, $h(a_1)=i$ and
$h(a_3), t(a_1)\neq i,j$;

\item paths $a_3 a_2 a_1$ with $a_2=a^*$, $t(a_3)=i$, $h(a_1)=j$ and
$h(a_3), t(a_1)\neq i,j$.
\end{itemize}

In the first case above, we have
$\widehat F_i \widehat F_j \widehat F_i(a_1) = a_1 =
\widehat F_j \widehat F_i \widehat F_j(a_1)$.

In the second case above, when $a_2\neq a_1^*$ or $a_1\neq a_2^*$, we have
$\widehat F_i \widehat F_j \widehat F_i(a_2a_1) = a_2a_1 =
\widehat F_j \widehat F_i \widehat F_j(a_2a_1)$.
When $a_2=a_1^*$, we have
\begin{align*}
\widehat F_i \widehat F_j \widehat F_i(a_2a_1) =&
\widehat F_i \widehat F_j( - \la_i + a_2a_1 ) =
\widehat F_i( -\la_i + a_2a_1 ) \\ =&
-\la_i -\la_j + a_2a_1.
\end{align*}
since $r_j(r_i(\la)) e_i = \la_j$; and on the other hand,
since $r_j(\la) e_i= \la_i + \la_j$, we find
$$
\widehat F_j \widehat F_i \widehat F_j(a_2a_1) =
\widehat F_j \widehat F_i (a_2a_1) =
\widehat F_j (-\la_i-\la_j + a_2a_1) =
-\la_i-\la_j + a_2a_1.
$$

The third case above is similar to the second case.

In the fourth and fifth cases above,
note that no two of the edges $a_1, a_2, a_3$ are reverse of
the other, so we have
$\widehat F_i \widehat F_j \widehat F_i (a_3a_2a_1)
= a_3a_2a_1 = \widehat F_j \widehat F_i \widehat F_j (a_3a_2a_1)$.
\end{proof}

\begin{lemma}   \label{eqc}
\vi If $\la_i\neq 0$, then $$\Pi_\la(Q)= \Pi_\la(Q) (1-e_i) \Pi_\la(Q).$$

\vii If $\mathfrak C_{ij}=-1$, and
$\la_i\neq 0$, $\la_j\neq 0$, $\la_i+\la_j\neq 0$, then
$$\Pi_\la(Q)= \Pi_\la(Q) (1-e_i-e_j) \Pi_\la(Q).$$
\end{lemma}
\begin{proof}
\vi As a $\Pi_\la(Q)$-module,
$\frac{\Pi_\la(Q)}{\Pi_\la(Q) (1-e_i) \Pi_\la(Q)}$ is zero at all
vertices not equal to $i$, so all edges of $\QQ$ must act by $0$.
But then it must also be zero at the vertex $i$ since
$\la_ie_i = \sum_{a\in Q; h(a)=i} aa^* - \sum_{a\in Q; t(a)=i} a^*a$.

\vii There is only one edge in $Q$ joining $i$ and $j$, say $a:i\to j$.
Let $V$ be the $\Pi_\la(Q)$-module
$\frac{\Pi_\la(Q)}{\Pi_\la(Q) (1-e_i-e_j) \Pi_\la(Q)}$.
Now $V$ is zero at all vertices not equal to $i$ or $j$, so $V=V_i\oplus V_j$.
Suppose $V\neq 0$, say $V_j\neq 0$. Then $aa^*=\la_je_j$ on $V_j$
implies that $a, a^*$ are nonzero maps, and $a$ has a right inverse
$\la_j^{-1} a^*$. But then $a^*a= - \la_ie_i$ on $V_i$ implies that
$a$ has a left inverse $-\la_i^{-1}a^*$. Hence, $\la_j=-\la_i$, a contradiction.
\end{proof}

Using  Proposition \ref{relations1}\vi, $\Pi_{r_i(\la)}(Q) (1-e_i)$ is a
right $(1-e_i)\Pi_{\la}(Q) (1-e_i)$-module, and
$\Pi_{r_i(\la)}(Q) (1-e_i-e_j)$ is a
right $(1-e_i-e_j)\Pi_{\la}(Q) (1-e_i-e_j)$-module.

\begin{corollary}   \label{eqccor}
\vi If $\la_i\neq 0$, then
$$F_i(M) = \Pi_{r_i(\la)}(Q) (1-e_i) \ot_{(1-e_i)\Pi_{\la}(Q) (1-e_i)}
(1-e_i)M$$ for any $M\in \Lmod{\Pi_{\la}(Q)}$.

\vii If $\mathfrak C_{ij}=-1$, and
$\la_i\neq 0$, $\la_j\neq 0$, $\la_i+\la_j\neq 0$, then
$$F_i(M) = \Pi_{r_i(\la)}(Q) (1-e_i-e_j)
\ot_{(1-e_i-e_j)\Pi_{\la}(Q)(1-e_i-e_j)}
(1-e_i-e_j) M$$
for any $M\in \Mod{\Pi_{\la}(Q)}$.
\end{corollary}

\begin{proof}
\vi Let $M\in \Pi_{\la}(Q)-\mod$. By Lemma \ref{eqc}\vi,
\begin{align*}
F_i(M)=& \Pi_{r_i(\la)}(Q) (1-e_i)
\ot_{(1-e_i)\Pi_{r_i(\la)}(Q) (1-e_i)} (1-e_i) F_i(M) \\ =&
\Pi_{r_i(\la)}(Q) (1-e_i)
\ot_{(1-e_i)\Pi_{\la}(Q) (1-e_i)} (1-e_i) M.
\end{align*}

The proof of \vii is similar, using Lemma \ref{eqc}\vii.
\end{proof}

\begin{corollary}  \label{frelations1}
\vi If $\la_i\neq 0$, then $F_i^2 =\id$.

\vii If $\mathfrak C_{ij}=0$, then $F_iF_j=F_jF_i$.

\viii If $\mathfrak C_{ij}=-1$ and $\la_i\neq 0$, $\la_j\neq 0$,
$\la_i+\la_j\neq 0$, then $F_iF_jF_i=F_jF_iF_j$.
\end{corollary}
\begin{proof}
\vii is trivial, while \vi and \viii are immediate from
Proposition \ref{relations1} and Corollary \ref{eqccor}.
\end{proof}

Our proof of Corollary \ref{frelations1} appears to be simpler
than earlier proofs, 
see \cite[Theorem 5.1]{CBH} (for ({\sf i})),
\cite[Remark 3.20]{Na}, \cite[Theorem 3.4]{Na},
\cite{Lu2}, and \cite{Maf}.

\subsection{Relations in higher rank.}\label{high}
In this subsection, $n$ is an integer greater than $1$.
We shall show that the reflection functors $F_i$ of (\ref{gan})
satisfy the Weyl group relations when the parameters are
generic.

%
%
%
%

We omit the proof of the following proposition since
it is completely similar to the proof of Proposition 
\ref{relations1}.

\begin{proposition}  \label{relations2}
Let $i,j\in I$. The homomorphisms $\widehat F_i$
of (\ref{rfisom}) satisfy the following for any $\la\in R$ and
$\nu \in \C$:

\vi Let $T= I\setminus\{i\}$. 
Then the map $$\widehat F_i:
\eel  \sfA_{n,r_i(\la),\nu}(Q)\eel  \to
\eel \sfA_{n,\la,\nu}(Q)\eel $$
is an isomorphism, and $\widehat F_i\ccirc \widehat F_i = \id$.

\vii Let $T=I\setminus\{i,j\}$. 
If $\mathfrak C_{ij}=0$, then
$$ \widehat F_i\ccirc \widehat F_j = \widehat F_j\ccirc\widehat F_i:
\eel \sfA_{n,r_ir_j(\la),\nu}(Q)\eel  \to
\eel \sfA_{n,\la,\nu}(Q)\eel .$$

\viii Let $T=I\setminus\{i,j\}$.
If $\mathfrak C_{ij}=-1$, then
$$ \widehat F_i\ccirc \widehat F_j \ccirc\widehat F_i=
\widehat F_j\ccirc\widehat F_i \ccirc\widehat F_j:
\eel \sfA_{n,r_ir_jr_i(\la),\nu}(Q)\eel  \to
\eel \sfA_{n,\la,\nu}(Q)\eel .$$  \qed
\end{proposition}

Next, we have the following generalization of Lemma \ref{eqc}.

\begin{lemma}
\vi Let $T=I\setminus\{i\}$.
If $\la_i\pm p\nu\neq 0$ for $p=0,1,...,n-1$, then
$ \sfA_{n,\la,\nu} =
 \sfA_{n,\la,\nu} \eel  \sfA_{n,\la,\nu}. $

\vii Let $T=I\setminus\{i,j\}$ and suppose 
$\mathfrak{C}_{ij}=0$.
If $\la_i\pm p\nu\neq 0$ and $\la_j\pm p\nu\neq 0$ 
for $p=0,1,...,n-1$, then
$ \sfA_{n,\la,\nu} =
\sfA_{n,\la,\nu} \eel  \sfA_{n,\la,\nu}. $

\vii Let $T=I\setminus\{i,j\}$ and suppose
$\mathfrak{C}_{ij}=-1$.
If $\la_i\pm p\nu\neq 0$, $\la_j\pm p\nu\neq 0$ and
$\la_i+\la_j\pm p\nu\neq 0$ for $p=0,1,...,n-1$, then
$ \sfA_{n,\la,\nu} =
\sfA_{n,\la,\nu} \eel  \sfA_{n,\la,\nu}. $
\end{lemma}

\begin{proof}
The proof is similar to the proof of Lemma \ref{eqc}.

To prove \vi, let $V$ be the $\sfA_{n,\la,\nu}$-module
$\frac{\sfA_{n,\la,\nu}}
{\sfA_{n,\la,\nu} \eel  \sfA_{n,\la,\nu}}$
where $T=I\setminus\{i\}$.
For any $n$-tuple of vertices $i_1,\ldots,i_n$, we let
$V_{i_1, \ldots, i_n} := (e_{i_1}\ot\cdots\ot e_{i_n}) V$,
so $V=\bigoplus_{i_1,...,i_n\in I}  V_{i_1, \ldots, i_n}$.
Since $\eel V=0$, we have
$V_{i_1, \ldots, i_n}=0$ when none of $i_1,...i_n$ is $i$.
Suppose now that
$i$ appears $m$ times in $i_1,...i_n$.
We shall prove by induction on $m$ that 
$V_{i_1, \ldots, i_n}=0$, so we assume that the statement
is true whenever $i$ appears less than $m$ times.
Without loss of generality, say $i_1=\cdots=i_m=i$.
Then by the relation \vi in Definition \ref{dea} and the 
induction hypothesis, we have
$$(\la_i + \nu \sum_{\ell=2}^m s_{1\ell})V_{i_1, \ldots, i_n}=0.$$
By \cite[Prop. 5.12]{Ga}, the element 
$\la_i + \nu \sum_{\ell=2}^m s_{1\ell}$ is invertible
in the group algebra $\C[S_n]$. Hence, 
$V_{i_1, \ldots, i_n}=0$, and \vi follows by induction.

The proofs of \vii and \viii are similar, using induction.
\end{proof}

As in the previous subsection, we obtain

\begin{corollary}\label{cox}
\vi Let $T=I\setminus\{i\}$. 
If $\la_i\pm p\nu\neq 0$ for $p=0,1,...,n-1$, then
$$ F_i(M) = \sfA_{n,r_i(\la),\nu}\eel 
\ot_{_{\eel  \sfA_{n,\la,\nu} \eel }} 
\eel  M,\quad \forall M\in  \Lmod{\sfA_{n,\la,\nu}}. $$

\vii Let $T=I\setminus\{i,j\}$ and suppose
$\mathfrak{C}_{ij}=0$.
If $\la_i\pm p\nu\neq 0$ and $\la_j\pm p\nu\neq 0$
for $p=0,1,...,n-1$, then
$$ F_i(M) = \sfA_{n,r_i(\la),\nu}\eel 
\ot_{_{\eel  \sfA_{n,\la,\nu} \eel }}
\eel  M,\quad \forall M\in  \Lmod{\sfA_{n,\la,\nu}}.$$

\viii Let $T=I\setminus\{i,j\}$ and suppose
$\mathfrak{C}_{ij}=-1$.
If $\la_i\pm p\nu\neq 0$, $\la_j\pm p\nu\neq 0$ and
$\la_i+\la_j\pm p\nu\neq 0$ for $p=0,1,...,n-1$, then
$$ F_i(M) = \sfA_{n,r_i(\la),\nu}\eel 
\ot_{_{\eel  \sfA_{n,\la,\nu} \eel }}
\eel  M,\quad \forall M\in  \Lmod{\sfA_{n,\la,\nu}}.\qquad\Box $$
\end{corollary}

Proposition \ref{weylrel} is immediate from
Proposition \ref{relations2} and Corollary \ref{cox}.


\subsection{Shift functors.} \label{shiftsubsec}
In this subsection, we return to the case when $Q$ is the
affine Dynkin quiver associated to $\Gamma$.

Let $\mathfrak C=(\mathfrak C_{ij})$ be the
generalized Cartan matrix of $Q$.
The affine Weyl group $\widetilde W$ is generated by the
simple reflections $r_i$ for $i\in I$. It acts on $\C^I$
by $r_i: \C^I\to\C^I$, where
$r_i(\la)= \la-\sum_{j\in I} \mathfrak C_{ij}\la_ie_j$.

Let $Q'$ be the finite Dynkin quiver obtained from $Q$
by deleting the vertex $o$. The Weyl group $W$
of $Q'$ is the subgroup
of $\widetilde W$ generated by the $r_i$ for $i\neq o$.
Let $\mathfrak C'= (\mathfrak C'_{ij})$ be the
Cartan matrix  of $Q'$. Then $W$ acts on
$\oplus_{i\neq o}\C e_i$ by $r_i(\la)= \la - \sum_{j\neq o}
\mathfrak C'_{ij} \la_i e_j$.
Denote by $w_0\in W$ the longest element of $W$.

If $i\in I$, then let $i^*\in I$ be the vertex such that
$N_{i^*} = N_i^*$.
Recall that if $\la=\sum_{i\in I}\la_i e_i$, then
$\overline \la = \sum_{i\in I}\la_{i^*} e_i$.

\begin{lemma} \label{longest}
For any $\la\in\C^I$ with $\la\cdot\delta=1$,
we have $w_0(\la)= -\overline \la +2e_0$.
\end{lemma}
\begin{proof}
The projection
$\C^I \too \oplus_{i\neq o}\C e_i$
is $W$-equivariant with kernel $\C e_o$.
We write $\la= \la_o e_o + \la'$ where $\la'\in
\oplus_{i\neq o}\C e_i$.
Now $w_0(\la)-w_0(\la')\in \C e_0$ and
 $w_0(\la')= -\overline{\la'}$. It follows that
$w_0(\la) = -\overline\la + 2(\la\cdot\delta)e_0$.
\end{proof}

We will now prove the second isomorphism in Corollary
\ref{maincorollary}.
For each vertex $i\neq o$, we have, from (\ref{rfisom}),
the homomorphism
$$\widehat F_i: \e_-\sfA_{n,r_i(\la),\nu-1}(Q)\e_-
\too \e_-\sfA_{n,\la,\nu-1}(Q)\e_-$$
which is an isomorphism by Proposition \ref{relations2}\vi.
By writing $w_0$ as a product of simple reflections, we get an
isomorphism
\begin{equation} \label{fwo}
\widehat F_{w_0} : \e_-\sfA_{n,w_0(\la),\nu-1}(Q)
\e_-\stackrel{\sim}{\too} 
\e_-\sfA_{n,\la,\nu-1}(Q)\e_-. 
\end{equation}
Proposition \ref{relations2} implies that this isomorphism
does not depend on the choice of presentation of $w_0$ as a
product of simple reflections.

Write $\hh_{t,k,c}=\hh_{t,k,c}(\GG).$
By \cite[(3.5.2)]{GG}, there is an
isomorphism $f^{\ot n}\hh_{t,k,c}f^{\ot n}
= \mathsf A_{n,\la,\nu}(Q)$ where $f=\sum_{i\in I}e_i$.
In particular, $\e_-\hh_{t,k-2t,c'}\e_- =
\e_- \mathsf A_{n,\overline{\la^\dag},\nu-1}(Q)\e_-$, and
$\e_-\hh_{t,k-2t,c}\e_- = \e_-\sfA_{n,\la, \nu-1}\e_-$.
By Lemma \ref{longest}, $\overline{\la^\dag}=w_0(\la)$,
so by (\ref{fwo}) we have the isomorphism
$$ \widehat F_{w_0} : \e_-\hh_{t,k-2t,c'}\e_-
 \stackrel{\sim}{\too} \e_-\hh_{t,k-2t,c}\e_-.$$
This completes the proof of Corollary \ref{maincorollary}.


Using the isomorphism
$\e\hh_{t,k,c}\e \simeq \e_-\hh_{t,k-2t, c}\e_-$ of
Corollary \ref{maincorollary}, we can consider
$\hh_{t,k-2t,c}\e_-$ as a
$(\hh_{t,k-2t,c}, \e\hh_{t,k,c}\e)$-bimodule.

\begin{definition}
The shift functor is defined to be
$$ \mathbb S : \Mod{\hh_{t,k,c}} \too \Mod{\hh_{t,k-2t, c}},
\quad V \mapsto \hh_{t,k-2t,c}\e_-\ot_{\e\hh_{t,k,c}\e} \e V. $$
\end{definition}

\section{Extended Dynkin quiver}

\subsection{$\Gamma$-analogue of commuting variety.}
In this subsection, we will prove a generalization of
\cite[Theorem 12.1]{EG}.

Let $R(\Gamma,n)$ be the space of extensions of the representation
$\C\Gamma \ot \C^n$ of $\Gamma$ to a representation
of $T(L) \rtimes \C\Gamma$, i.e.,
 $$ R(\Gamma,n) := \Hom_\Gamma(L, \End_\C (\C\Gamma \ot \C^n) ). $$
Let $\Z=\Z(\Gamma,n)$ be the (not necessarily reduced) subscheme
of $R(\Gamma,n)$ consisting of those representations $\rho$
such that $\rho([X, Y])=0$ for all $X,Y\in L$.

Now $\C\Gamma = \oplus_{i\in I} \End(N_i)$. Let
$p_i\in \C\Gamma$ be the idempotent element corresponding
to the identity element of $\End(N_i)$.
Define $\Z_1=\Z_1(\Gamma,n)$ to be the (not necessarily reduced) subscheme
of $R(\Gamma,n)$ consisting of those representations $\rho$
such that, for all $X,Y\in L$, we have
$\rho([X, Y]p_i)=0$ for $i\neq o$,
and $\wedge^2 \rho([X,Y]p_o)=0$.
We remark that $\rho([X,Y]p_o)$ is a $n\times n$-matrix
and $\wedge^2 \rho([X,Y]p_o)=0$ means that all $2\times 2$ minors
of this matrix vanish (so its rank is at most 1).

We shall denote by $J$ and $J_1$ the defining ideals of
$\Z(\Gamma,n)$ and $\Z_1(\Gamma,n)$, respectively.
Thus, $\Z(\Gamma,n) =\mathrm{Spec} \C[R(\Gamma,n)] / J$
and $\Z_1(\Gamma,n) =\mathrm{Spec} \C[R(\Gamma,n)] / J_1$.
Let $G := \mathtt{Aut}_\Gamma (\C\Gamma \ot \C^n).$
Observe that the group $G$ acts on $R(\Gamma,n)$, $\Z(\Gamma,n)$,
and $\Z_1(\Gamma,n)$.

\begin{theorem}  \label{idealrank1}
One has: $J^G=J_1^G$.
\end{theorem}

It is clear that $J\supset J_1$, so $J^G\supset J_1^G$.
To prove Theorem \ref{idealrank1}, we have to show that
$J^G \subset J_1^G$.
We need the following lemmas.
First, let us fix a basis $X,Y$ for $L$.

\begin{lemma}
The ideal $J^G$ is generated in $\C[R(\Gamma,n)]^G$ by
functions of the form $\rho\mapsto \Tr(\rho(Q[X,Y]))$,
where $Q\in T(L)\rtimes \C\Gamma$.
\end{lemma}
\begin{proof}
This follows from Weyl's fundamental theorem of invariant theory.
\end{proof}

Therefore, it suffices to show that $\Tr(\rho(Q[X,Y]p_i)) = 0$
(mod $J_1$) for all $\rho\in R(\Gamma,n)$, $Q\in T(L)\rtimes \C\Gamma$,
and $i\in I$.
This is obvious for $i\neq o$ from the definition of $J_1$.
For $i=o$, we shall prove it by induction on the degree of $Q$.
The case $\deg Q=0$ is clear, so let $d>0$ and assume that
$\Tr(\rho(Q[X,Y]p_o)) = 0$ (mod $J_1$) whenever $\deg Q<d$.

\begin{lemma}
Let $\deg Q =d$.
If $Q= Q_1[X,Y]Q_2$ for some $Q_1,Q_2 \in T(L)\rtimes \C\Gamma$,
then $\Tr(\rho(Q[X,Y]p_o))=0$ (mod $J_1$).
\end{lemma}
\begin{proof}
We may replace $Q, Q_1, Q_2$ by $p_oQp_o, p_oQ_1p_o, p_oQ_2p_o$
respectively.
Modulo $J_1$, and writing in terms of
matrix elements, we have
\begin{align*}
& \Tr(\rho(Q[X,Y]p_o)) = \Tr(\rho(Q_1[X,Y]p_o Q_2[X,Y]p_o)) \\
&= \sum \rho(Q_1)_{lm} \rho([X,Y]p_o)_{mq} \rho(Q_2)_{qr}
\rho([X,Y]p_o)_{rl} \\
&= \sum \rho(Q_1)_{lm} \rho([X,Y]p_o)_{ml} \rho(Q_2)_{qr}
\rho([X,Y]p_o)_{rq}
\quad\text{ (since $\wedge^2 \rho([X,Y]p_o) = 0$) }\\
&= \Tr(\rho(Q_1[X,Y]p_o)) \Tr(\rho(Q_2[X,Y]p_o)).
\end{align*}
This is equal to zero by induction hypothesis.
\end{proof}

Let $\varphi: T(L)\rtimes \C\Gamma \to S(L)\rtimes \C\Gamma$
be the quotient map, where $S(L)$ denotes the symmetric
algebra on $L$.
The preceding lemma implies the following corollary.

\begin{corollary}
If $\deg Q\leq d$ and $Q\in \Ker \varphi$, then
$\Tr(\rho(Q[X,Y]p_o)) = 0$ (mod $J_1$).  \qed
\end{corollary}

Note that elements of the form
$(aX+bY)^m$ (where $a,b\in \C$) span
a set of representatives of $S(L)$ in $T(L)$.
Thus, it remains to show that
$\Tr(\rho((aX+bY)^m[X,Y]p_o)) = 0$ (mod $J_1$) for any $a,b\in\C$
and $m\leq d$.
This is equivalent to showing that
$\Tr(\rho((aX+bY)^m[X,Y])) = 0$ (mod $J_1$).
But we have
$$ \Tr(\rho((aX+bY)^m[X,Y])) =
\frac{1}{a(m+1)}\Tr(\rho([(aX+bY)^{m+1},Y])) = 0. $$
This completes the proof of Theorem \ref{idealrank1}.

\subsection{}
Let $\mu_{\mathsf{CM}}: \Rep_{\al}(\QQc) \to 
\mathfrak{gl}(\al)$ be the moment map,
and $\Zc = \mu_{\mathsf{CM}}^{-1}(0)$ the scheme theoretic
inverse image of the point $0$.
It was proved in \cite[Theorem 1.3.1]{GG2}
that $\Zc$ is a reduced scheme.
Now, there are natural algebra morphisms
\begin{equation} \label{pstar}
\C[\Z]^G \stackrel{f}{\longleftarrow} \C[\Z_1]^G 
\stackrel{g}{\too} \C[\Zc]^G .
\end{equation}
By Theorem \ref{idealrank1}, $f$ is an isomorphism.
The following proposition and its proof is a straightforward
generalization of \cite[Proposition 2.8.2]{GG2}, given our
Theorem \ref{idealrank1}.

\begin{proposition} \label{ppp}
The morphism $g$ in (\ref{pstar}) is an isomorphism.  \qed
\end{proposition}

From Proposition \ref{ppp}
and \cite[Theorem 1.3.1]{GG2}, we have the following 
generalization of \cite[Theorem 1.2.1]{GG2}.

\begin{theorem} \label{idealreduced}
One has: $J^G= \sqrt{J}^G$.  \qed
\end{theorem}

Let $\Z^{red} := \mathrm{Spec} \C[\Rep_{n\delta}(\QQ)]/\sqrt{J}$,
a closed subvariety of $\Rep_{n\delta}(\QQ)$.
Define an embedding
$\jmath: \lreg\too \Rep_{n\delta}(\QQ)$ by
$\jmath(u_1,\ldots,u_n)_a = (\phi_a(u_1), \ldots, \phi_a(u_n))$
for any $a\in \QQ$.
Using formulas \eqref{rpp} from \S\ref{lllpf} below,
we deduce that the image of $\jmath$ lies in $\Z^{red}$.
Pullback of functions gives a morphism
\begin{equation} \label{gradedmapred}
\jmath^*:  \C[\Z^{red}]^G \too \C[\lreg]^{\GG}.
\end{equation}
By \cite[Theorem 3.4]{CB} and \cite[Corollary 3.2]{Kr}, we have
the following proposition.
\begin{proposition} \label{gris}
The map $\jmath^*$ in (\ref{gradedmapred}) is an isomorphism. \qed
\end{proposition}

\section{Proof of Proposition \ref{radialpart}}\label{pf_radialpart}
\subsection{} The formula of Proposition \ref{radialpart}(i)  is clear.
Next, we have
$$ \th^{\text{Holland}}(a^*) =
\sum_{p,q} e^a_{q,p} \ot \frac{\partial}{\partial t^a_{p,q}}. $$

To compute the restriction of
$e^a_{q,p} \ot \frac{\partial}{\partial t^a_{p,q}}$
to $\jmath(\lreg)$ at a point
$\mathbf u=(u_1,\ldots,u_n)\in
\lreg$, let $g_{p,q}(\varepsilon)=
\id+\varepsilon\mathbf B_{p,q}$ be an element of $\GL(\alpha)$
such that
\begin{equation} \label{gpq}
g_{p,q}(\varepsilon)\cdot (\jmath(\mathbf u) +
\varepsilon e^a_{p,q})
= \jmath(\mathbf u)+ \varepsilon\jmath(\mathbf w),
\quad \mathbf w\in \lreg,
\end{equation}
where we omit terms of higher order in $\varepsilon$.
Then for a function $f\in \mathcal O(\chi,N)$, we have
\begin{align}  \label{rest}
e^a_{q,p}\ot \frac{\partial}{\partial \varepsilon}
f(\jmath(\mathbf u) + \varepsilon e^a_{p,q})
=& e^a_{q,p}\ot\frac{\partial}{\partial \varepsilon}
f( g_{p,q}(\varepsilon)^{-1}\cdot (
\jmath(\mathbf u)+ \varepsilon\jmath(\mathbf w) ) \nonumber \\
=&
e^a_{q,p}\ot \frac{\partial}{\partial \varepsilon}
\chi(g_{p,q}(\varepsilon)) g_{p,q}(\varepsilon)^{-1}
f(\jmath(\mathbf u)+ \varepsilon\jmath(\mathbf w)) \\
=&
e^a_{q,p}\ot
\frac{\partial}{\partial \mathbf w} f(\jmath(\mathbf u))
+ e^a_{q,p}\ot\left(  \sum\nolimits_{j\in I}\,
\chi_j \Tr(\mathbf B_{p,q}^{(j)}) \id - \mathbf B_{p,q}
\right) f(\jmath(\mathbf u))\nonumber
\end{align}
where $\mathbf B_{p,q}^{(j)}$ is the component of $\mathbf B_{p,q}$
in $\gl(\al_j)$.
We shall write $\mathbf B_{p,q}^{(j)}$ for $j\in I$
as a $n\times n$ block matrix $\oplus_{1\leq \ell,m \leq n}
\mathbf B_{p,q}^{(j)}(\ell,m)$ where
$\mathbf B_{p,q}^{(j)}(\ell,m)\in \gl(\delta_j)$ is the
$(\ell, m)$-th block.
Similarly, we write $e^a_{p,q}$ as $\oplus_{1\leq\ell,m\leq n}
e^a_{p,q}(\ell,m)$.
By (\ref{gpq}), we need to solve the equations:
\begin{gather*}
\mathbf B_{p,q}^{(h(a))}(\ell,m) \phi_a(u_m) -
\phi_a(u_\ell) \mathbf B_{p,q}^{(t(a))}(\ell,m) + e^a_{p,q}(\ell,m)
=  \left\{ \begin{array}{ll}
0 & \text{ if } \ell\neq m   \\
\phi_a(w_\ell) & \text{ if } \ell=m
\end{array} \right.   \\
\sum_{m=1}^n \mathbf B_{p,q}^{(o)}(\ell,m) - \mathbf B_{p,q}^{(s)}
=0.
\end{gather*}
where $1\leq \ell,m \leq n$.
We shall set $\mathbf B^{(s)}_{p,q}=0$.

Suppose $(\ell-1)\delta_{h(a)}< p \leq \ell\delta_{h(a)}$
and $(m-1)\delta_{t(a)}< q \leq m\delta_{t(a)}$
where $\ell, m \in [1,n]$.
If $\ell\neq m$, then we set $\mathbf B_{p,q}^{(j)}(\ell',m')=0$
whenever $\ell'\neq m'$ and $(\ell',m')\neq (\ell,m)$.
If $\ell=m$, then we set $\mathbf B_{p,q}^{(j)}(\ell',m')=0$
whenever $\ell'\neq m'$.

\subsection{Proof of (\ref{radialml}).}\label{lllpf}
First of all, 
it is  immediate from (\ref{rp}) that
\begin{align}\label{rpp}
&\sum_{a\in Q; t(a)=i} \phi_{a^*}(u) \phi_a(w) =
\delta_i \om (w,u) \id_{N_i^*},\quad\text{for each source $i$ in $Q$};\\
& \sum_{a\in Q; h(a)=j} \phi_a(w) \phi_{a^*}(u) =
\delta_j \om (w,u) \id_{N_j^*},\quad\text{for each  sink $j$ in $Q$.}\nonumber
\end{align}

Next, we find a collection
of operators $\beta_i \in \End(N_i^*)$ such that
\begin{equation}  \label{offdiagproblem}
\phi_a(u) \beta_i - \beta_j \phi_a(w) = f_a
\end{equation}
where $i=t(a)$, $j=h(a)$, and
$f_a : N_i^* \to N_j^*$ are given operators.
We write the collection $\beta_i$ as an element
$\sum \beta_\ga \ga$ of $\C[\Gamma]$.
Since $\ga \phi_a(w) = \phi_a(\ga w)\ga$, we get
$$ \sum \beta_\ga \phi_a(u-\ga w)\ga = f_a ,\qquad
 \sum \beta_\ga \ga \phi_a(\ga^{-1}u-w) = f_a .$$
We  multiply the first equation above  by
$\phi_{a^*}(u)$ on the left and add over all edges going out from $i$.
Similarly, let us multiply the second equation  above  by
$\phi_{a^*}(w)$ on the right and add over all edges going into $j$.

Using  formulas \eqref{rpp}, we obtain:
\begin{align*}
& \delta_i \sum \beta_\ga \omega(u,\ga w) \ga|_{N_i^*}
  = \sum_{a\in Q; t(a)=i} \phi_{a^*}(u) f_a,\quad\text{for sources $i$;}\\
&\delta_j \sum \beta_\ga \omega(u,\ga w) \ga|_{N_j^*}
  = \sum_{a\in Q; h(a)=j} f_a \phi_{a^*}(w),\quad\text{for sinks $j$.}
\end{align*}
This implies that
\beq{offdiagsoln}
\beta_\ga = \omega(u,\ga w)^{-1}|\Gamma|^{-1}
\sum_{a \in Q}
\left( \Tr|_{N^*_{h(a)}}( f_a\phi_{a^*}(w) \ga^{-1} ) +
\Tr|_{N^*_{t(a)}}( \phi_{a^*}(u) f_a\ga^{-1} )
\right).
\eeq
Hence, if $\ell\neq m$, for $\oplus_{j\in I} \mathbf
B_{p,q}^{(j)}(\ell,m)$
we get the expression
\begin{align} \label{mllm}
&\sum_{\ga\in\Gamma}\omega(u_\ell,\ga u_m)^{-1}|\Gamma|^{-1}
\left( \Tr|_{N^*_{h(a)}}( e^a_{p,q}\phi_{a^*}(u_m) \ga^{-1} ) +
\Tr|_{N^*_{t(a)}}( \phi_{a^*}(u_\ell) e^a_{p,q}\ga^{-1} )
\right) \ga.
\end{align}
and so,  for $\mathbf B_{p,q}^{(o)} (\ell,\ell)
= -\mathbf B_{p,q}^{(o)} (\ell,m)$ we obtain the expression
$$
 - \sum_{\ga\in\Gamma}\omega(u_\ell,\ga u_m)^{-1}|\Gamma|^{-1}
\left( \Tr|_{N^*_{h(a)}}( e^a_{p,q}\phi_{a^*}(u_m) \ga^{-1} ) +
\Tr|_{N^*_{t(a)}}( \phi_{a^*}(u_\ell) e^a_{p,q}\ga^{-1} )
\right).
$$
Thus, for all $j\in I$, for $\mathbf B_{p,q}^{(j)} (\ell,\ell)$ we
 obtain the expression
$$
- \sum_{\ga\in\Gamma}\omega(u_\ell,\ga u_m)^{-1}|\Gamma|^{-1}
\left( \Tr|_{N^*_{h(a)}}( e^a_{p,q}\phi_{a^*}(u_m) \ga^{-1} ) +
\Tr|_{N^*_{t(a)}}( \phi_{a^*}(u_\ell) e^a_{p,q}\ga^{-1} )
\right)  \id_{\delta_i\times \delta_i}.
$$

It follows from the last formula and from (\ref{rest})
that for $\ell\neq m$, the $(m,\ell)$-entry of the
radial part of $\th^{\text{Holland}}(a^*)$ is
\begin{gather*}
\sum_{p=(\ell-1)\delta_{h(a)}+1}^{\ell\delta_{h(a)}}
\sum_{q=(m-1)\delta_{t(a)}+1}^{m\delta_{t(a)}}  e^a_{q,p}(m,\ell)
\left( \sum_{j\in I} \chi_j \Tr(\mathbf B_{p,q}^{(j)} (\ell,\ell))
- \sum_{j\in I}\mathbf B_{p,q}^{(j)} (\ell,\ell)
\right) \\
= |\Gamma|^{-1} \sum_\ga
\frac{(\phi_{a^*})_m\cdot\ga^{-1} +
\ga^{-1}(\phi_{a^*})_\ell}
{\omega(\ga;\ell,m)} (1 - \sum_j \chi_j\delta_j)
= \frac{k}{2} \sum_\ga
\frac{(\phi_{a^*}\ccirc(\id+ \ga^{-1}))_{m,\ell}}
{\omega(\ga;\ell,m)}\ga^{-1}.
\end{gather*}
Note that, since $\zeta$ acts by $-1$ on $L$
and by $1$ on $N^*_{h(a)}$,
$$ \frac{(\phi_{a^*})_m}{\omega(\ga \zeta;\ell,m)}
(\ga\zeta)^{-1}
= - \frac{(\phi_{a^*})_m}
{\om(\ga;\ell,m)}\ga^{-1} $$
and so
$$ \sum_\ga \frac{(\phi_{a^*})_m}
{\om(\ga;\ell,m)}\ga^{-1} = 0. $$
Hence, the $(m,\ell)$-entry of the
radial part of $\th^{\text{Holland}}(a^*)$ is equal to
$$ -\frac{k}{2} \sum_\ga \frac{(\phi_{a^*}\ccirc\ga)_\ell}
{\omega(\ga;m,\ell)} \ga. $$

\begin{proof}[Proof of Lemma \ref{lll}.]
We set $f_a=0$ in (\ref{offdiagproblem}).
Then from (\ref{offdiagsoln}), we have
$\omega(u, \ga w)\beta_\ga =0$ for all $\ga\in \Gamma$.
Since not all $\beta_\ga$ are zero, we must have
$\omega(u, \ga w)=0$ for some $\gamma$.\end{proof}

\subsection{Proof of (\ref{radialmm}).}
We  need to solve
$\dis \phi_a(u) \beta_i - \beta_j \phi_a(u) =
f_a - \phi_a(w) .$

As above, for $\ga \neq 1, \zeta$, we obtain
$$ \beta_\ga = \omega(u,\ga u)^{-1}|\Gamma|^{-1}
\sum_{a \in Q}
\left(
\Tr|_{N^*_{h(a)}}( f_a\phi_{a^*}(u) \ga^{-1} )
+  \Tr|_{N^*_{t(a)}}( \phi_{a^*}(u)f_a \ga^{-1} )
\right). $$

Moreover, multiplying on the right by $\phi_{a^*}(v)$
and summing over all incoming edges $a\in Q$
at the vertex $j$, we get
$$ \delta_j \sum \beta_\ga \omega(\ga^{-1}u- u,v)  \ga|_{N^*_j}
= \sum_{a\in Q: h(a)=j} f_a\phi_{a^*}(v) - \delta_j \omega(w,v). $$
Take the trace of both sides this equation
and sum up over all sinks $j$. We have
$\oplus_{\text{ sink } j} (N_j^*)^{\oplus \delta_j} =
\C[\Gamma/S],$ where $S=\{1,\zeta\}$. It follows that
the trace of $\ga$ in the last sum vanishes if $\ga\neq 1,\zeta$.
Let $\beta_\zeta=0$. Then
\begin{equation} \label{diffw}
w = 2|\Gamma|^{-1}\sum_{a\in Q} \Tr(f_a \phi_{a^*}).
\end{equation}

Hence, for $\ell=m$, we get
\begin{align} \label{Bll}
& \oplus_{j\in I} \mathbf B_{p,q}^{(j)}(m,m) \\ =&
\sum_{\ga\neq 1,\zeta}
\omega(u_m,\ga u_m)^{-1}|\Gamma|^{-1}
\left(
\Tr|_{N^*_{h(a)}}( e^a_{p,q}(\phi_{a^*})_m \ga^{-1} )
+  \Tr|_{N^*_{t(a)}}( (\phi_{a^*})_me^a_{p,q} \ga^{-1} )
\right) (\ga-1).\nonumber
\end{align}

It follows from (\ref{rest}), (\ref{diffw}), (\ref{Bll}) and
(\ref{mllm})
that the $(m,m)$-entry of
the radial part of $\th^{\text{Holland}}(a^*)$ is
\begin{gather}
\frac{2}{|\Gamma|} \frac{\partial}{\partial (\phi_{a^*})_m} +
\sum_{p=(m-1)\delta_{h(a)}+1}^{m\delta_{h(a)}}
\sum_{q=(m-1)\delta_{t(a)}+1}^{m\delta_{t(a)}}
e^a_{q,p}(m,m)
\left(
\sum_{j\in I} \chi_j \Tr(\mathbf B_{p,q}^{(j)}(m,m))
- \sum_{j\in I} \mathbf B_{p,q}^{(j)}(m,m)
\right)  \nonumber \\ 
- \sum_{\ell\neq m}
\sum_{p=(\ell-1)\delta_{h(a)}+1}^{\ell\delta_{h(a)}}
\sum_{q=(m-1)\delta_{t(a)}+1}^{m\delta_{t(a)}}
e^a_{q,p}(m,\ell)
\sum_{j\in I} \mathbf B_{p,q}^{(j)}(\ell,m)=\nonumber
\end{gather}
\begin{gather} \label{last}
= \frac{2}{|\Gamma|} \frac{\partial}{\partial (\phi_{a^*})_m} +
\frac{1}{|\Gamma|} \sum_{\ga \neq 1,\zeta}
\frac{\ga^{-1}\cdot(\phi_{a^*})_m+(\phi_{a^*})_m\cdot\ga^{-1}}
{\om(\ga;m,m)}
\left( -\ga +1 - \sum_j \chi_j (\delta_j-\Tr|_{N^*_j}(\ga))
\right)  \nonumber \\
- \frac{1}{|\Gamma|}\sum_{\ell\neq m} \sum_\ga
\frac{(\phi_{a^*}\ccirc(\id+\ga^{-1}))_{m,\ell}}{\om(\ga;\ell,m)}
\\
= \frac{2}{|\Gamma|} \frac{\partial}{\partial (\phi_{a^*})_m} +
\frac{1}{|\Gamma|} \sum_{\ga \neq 1,\zeta}
\frac{(\phi_{a^*}\ccirc(\ga^{-1}+\id))_{mm}}{\om(\ga;m,m)}
\left( -1 + |\Gamma|c_\ga \ga^{-1} \right)
+ \frac{1}{|\Gamma|}\sum_{\ell\neq m}  \sum_\ga
\frac{(\phi_{a^*}\ccirc\ga)_\ell}{\omega(\ga;m,\ell)} . \nonumber
\end{gather}
The last term in (\ref{last}) comes from (\ref{mllm}).

It is even easier to compute the radial part
for the edge $b: s\to o$. We omit this computation.
This completes the proof of  Proposition \ref{radialpart}.

\section{Proof of  Theorem \ref{mess}}\label{pf_mess}
\subsection{} It easy to check that the operators $R_{ml}^v$ have 
the following
$\GG$-equivariance properties:
\begin{gather*}
\gamma_m R_{ml}^v=R_{ml}^{\gamma(v)}\gamma_m,\quad \gamma_l
R_{ml}^v=R_{ml}^v\gamma_l,\\
s_{ml}R^v_{ml}=R^v_{lm}s_{ml},\quad
s_{mj}R_{ml}^v=R^v_{jl}s_{mj},\quad s_{lj}R^v_{ml}=R^v_{mj}s_{lj},
\end{gather*}
where $j\ne m,l$. It implies that
$\dis g \Theta^{{\text{Dunkl}}}(v)=\Theta^{{\text{Dunkl}}}{g(v)},$
for any $g\in\GG$ and $v\in \lreg$.

Next, we prove that
\begin{equation}
[\Theta^{{\text{Dunkl}}}(w_i),\Theta^{{\text{Dunkl}}}(v_i)]=\omega(w,v)
\left( t\cdot 1+
\frac{k}{2} \sum_{j\neq i}\sum_{\ga\in\Gamma}
s_{ij}\ga_{i}\ga_{j}^{-1} + \sum_{\ga\in\Gamma\smallsetminus\{1\}}
c_{\ga}\ga_{i}\right) , \qquad 1\le i\le n.
\label{relation1Theta}\end{equation}

First we prove
\begin{equation}\label{RD1}
[R_{ij}^w,\Theta_i^{{\text{Dunkl}}}(v)]+[\Theta_i^{{\text{Dunkl}}}(w),R_{ij}^v]=\frac12
\omega(w,v) \sum_{j\neq i}\sum_{\ga\in\Gamma} s_{i
j}\ga_{i}\ga_{j}^{-1},\quad 1\le i\ne j\le n.
\end{equation}

We  prove that (\ref{RD1}) holds if we apply both sides to
$\overset{{\,}_\to}f$, a basis vector in
$\FF$ such that  $f_i,f_j=f^+$,
cf. \eqref{basis}. Indeed, in that case
we have 
\begin{align*}
&[R_{ij}^w,\Theta_i^{{\text{Dunkl}}}(v)](\overset{{\,}_\to}f)=\left(R_{ij}^w (D^{v})_i
+\Theta_i^{{\text{Dunkl}}}(v)\frac12\sum_{\gamma\in\Gamma}
\frac{(\gamma w)^\vee_j}{\om({1};i,j)}\right) (\overset{{\,}_\to}f)=A \overset{{\,}_\to}f\\
&[\Theta^{{\text{Dunkl}}}_i(w),R_{ij}^v](\overset{{\,}_\to}f)=\left(-\frac12
\Theta^{{\text{Dunkl}}}_i(w)\sum_{\gamma\in\Gamma}
\frac{(\gamma v)^\vee_j }{\om({\ga\inv};i,j)}-R_{ij}^w(D^w)_i\right)
(\overset{{\,}_\to}f)=B \overset{{\,}_\to}f,
\end{align*}
where
$$
A=-\frac12
\sum_{\gamma\in\Gamma}\frac{  v_i^\vee  (\gamma w)^\vee_j}{\om({\ga\inv};i,j)}
 s_{ij}\gamma_i\gamma_j^{-1} ,\quad\text{resp.},
\quad
B=
\frac12\sum_{\gamma\in\Gamma}\frac{w^\vee _i (\gamma
 v)_j^\vee }{\om({\ga\inv};i,j)}s_{ij}\gamma_i\gamma_j^{-1}.
$$
These formulas yield
\begin{multline*}
\left([R_{ij}^w,\Theta_i^{{\text{Dunkl}}}(v)]+
[\Theta_i^{{\text{Dunkl}}}(w),R_{ij}^v]\right)\overset{{\,}_\to}f=
-\frac12\sum_{\gamma\in\Gamma}\frac{ v_i^\vee  (\gamma
w)^\vee_j-w^\vee _i(\gamma v)^\vee _j}{\om({\ga\inv};i,j)}
s_{ij}\gamma_i\gamma_j^{-1} \overset{{\,}_\to}f=\\
\frac12 \omega(w,v)\sum_{\gamma\in\Gamma}
s_{ij}\gamma_i\gamma_j^{-1} \overset{{\,}_\to}f.
\end{multline*}

We  consider the case $f_i=f^-$,
$f_j=f^+$. Then we have:
\begin{equation*}
[R_{ij}^w,\Theta^{{\text{Dunkl}}}_i(v)]\overset{{\,}_\to}f=-R_{ij}^w v_i^\vee
\overset{{\,}_\to}f=\frac12\sum_{\gamma\in\Gamma}\frac{(\gamma
w)_j^\vee }{\om({\ga\inv};i,j)} v_i^\vee
s_{ij}\gamma_i\gamma_j^{-1} \overset{{\,}_\to}f=0,
\end{equation*}
because in the sum the terms corresponding $\gamma$ and
$\zeta\gamma$ mutually cancel. Analogously,
\begin{equation*}
[\Theta^{{\text{Dunkl}}}_i(w),R_{ij}^w]\overset{{\,}_\to}f=R_{ij}^v
w_i^\vee \overset{{\,}_\to}f=-\frac12\sum_{\gamma\in\Gamma}\frac{(\gamma
v)_j^\vee }{\om({\ga\inv};i,j)}s_{ij}\gamma_i\gamma_j^{-1} w_i^\vee
\overset{{\,}_\to}f=0.
\end{equation*}
For the same reason we also have
$\sum_{\gamma\in\Gamma} s_{ij}\gamma_i\gamma_j^{-1}\overset{{\,}_\to}f=0.$

We  consider the case $f_i=f^+,f_j=f^-.$ A similar argument yields
\begin{equation*}
[R_{ij}^w,\Theta^{{\text{Dunkl}}}_i(v)]\overset{{\,}_\to}f=
[\Theta^{{\text{Dunkl}}}_i(w),R_{ij}^w]\overset{{\,}_\to}f=
\sum_{\gamma\in\Gamma} s_{ij}\gamma_i\gamma_j^{-1}\overset{{\,}_\to}f=0.
\end{equation*}

 The case $f_i,f_j=f^-$ is analogous to the first case and we
have
\begin{gather*}
[R_{ij}^w,\Theta_i^{{\text{Dunkl}}}(v)]\overset{{\,}_\to}f=\frac12\sum_{\gamma\in\Gamma}
\frac{w_i^\vee (\gamma  v)_j^\vee }{\om({\ga\inv};i,j)}\overset{{\,}_\to}f,\\
[\Theta_i^{{\text{Dunkl}}}(w),R_{ij}^v]\overset{{\,}_\to}f=-\frac12
\sum_{\gamma\in\Gamma}\frac{v_i^\vee (\gamma w)_j^\vee }{
\om({\ga\inv};i,j)}
s_{ij}\gamma_i\gamma_j^{-1}\overset{{\,}_\to}f,
\end{gather*}
and
\begin{multline*}
\left([R_{ij}^w,\Theta_i^{{\text{Dunkl}}}(v)]+
[\Theta_i^{{\text{Dunkl}}}(w),R_{ij}^v]\right)\overset{{\,}_\to}f=
-\frac12\sum_{\gamma\in\Gamma}\frac{w^\vee _i (\gamma
 v)_j^\vee - v_i^\vee  (\gamma w)_j^\vee }{\om({\ga\inv};i,j)}
s_{ij}\gamma_i\gamma_j^{-1} \overset{{\,}_\to}f=\\
\frac12 \omega(w,v)\sum_{\gamma\in\Gamma}
s_{ij}\gamma_i\gamma_j^{-1} \overset{{\,}_\to}f.
\end{multline*}

We  remark that for any $1\le j\ne i\ne k\le n$ and $w,v\in L$ we
have
$\dis R^w_{ij}R^v_{ik}=R^v_{ik}R^w_{ij}=0.$
Now, (\ref{relation1Theta}) follows from this equation,
the $n=1$ case of
Theorem \ref{mess}, and (\ref{RD1}).

\subsection{}
Next we prove:
\begin{equation}
[\Theta^{\text{Dunkl}}(w_i),\Theta^{\text{Dunkl}}(v_m)]= -\frac{k}{2}
\sum_{\ga\in\Gamma} \omega_{L}(\ga u,v) s_{i
m}\ga_{i}\ga_{m}^{-1}, \quad  1\le i\ne m\le
n.\label{relation2Theta}
\end{equation}
To this end, we rewrite the RHS of (\ref{relation1Theta}) as follows
\begin{multline*}
[\Theta^{\text{Dunkl}}(w_i),\Theta^{\text{Dunkl}}(v_m)]=k([\Theta^{\text{Dunkl}}_i(w),R_{mi}^v]+
[R^w_{im}\Theta^{\text{Dunkl}}_m(v)])+\\k^2([R_{im}^w,R_{mi}^v]+\sum_{j\ne
m,i}[R_{im}^w,R_{mj}^v]+[R_{ij}^w,R_{mj}^v]+[R_{ij}^w,R_{mi}^v]).
\end{multline*}
We  first prove that
\begin{equation*}
[R_{im}^w,R_{mj}^v]+[R_{ij}^w,R_{mj}^v]+[R_{ij}^w,R_{mi}^v]=0,\quad
j\ne m,i.
\end{equation*}
For that it is enough to show that
\beq{BraidDunkl1}
R_{im}^wR_{mj}^v-R_{mj}^vR_{ij}^w+R^w_{ij}R_{mi}^v=0,
\quad\text{and}\quad
-R_{mj}^vR_{im}^w+R_{ij}^wR_{mj}^v-R_{mi}^v
R_{ij}^w=0.
\eeq

We  prove the first equation, the second is proved similarly.
Let $\overset{{\,}_\to}f$ be a basis vector in $\FF$ such that
 $f_i,f_j,f_m=f^+$, cf. \eqref{basis}.
We compute
\begin{align*}
4(R_{im}^wR_{mj}^v-R_{mj}^vR_{ij}^w+R^w_{ij}R_{mi}^v)(\overset{{\,}_\to}f) &=
\left(
\sum_{\beta ,\gamma\in\Gamma}\frac{(w)_i^\vee (\gamma
v)_j^\vee }{\om({\beta \inv};i,m)\om({(\gamma\beta )\inv};i,j)}
s_{im}s_{mj}\beta _i\gamma_m
(\gamma\beta )^{-1}_j \right. \\ & -
\left. \sum_{\gamma,\beta \in\Gamma}
\frac{(\beta v)_j^\vee (\beta ^{-1}\gamma w)_m^\vee }{
\om({\beta \inv};m,j)\om({\beta };i,m)}
s_{mj}s_{ij}(\beta ^{-1}\gamma)_i\beta _m\gamma^{-1}_j
 \right. \\ & +
\left. \sum_{\gamma,\beta \in\Gamma}\frac{(\beta 
w)_j^\vee (\beta \gamma v)_j^\vee }{
\om({\beta \inv};i,j)\om({(\beta \ga)\inv};m,j)}
s_{ij}s_{mi}\gamma_i^{-1}(\beta \gamma)_m\beta _j^{-1}\right)\overset{{\,}_\to}f.
\end{align*}

We  change summation indices at the first term as
$\gamma\to\gamma\beta ,\beta \to \gamma$ and at the
third term as $\gamma\to\beta ^{-1},\beta \to
\gamma\beta $. We get
\begin{multline*}
\sum_{\gamma,\beta \in\Gamma}\left(\frac{ w_i^\vee
(\gamma v)_j^\vee }{
\om({\beta \inv};i,m)
\om({\beta \inv};i,j)}-\frac{(\gamma
v)_j^\vee (\beta  w)^\vee _m}
{\om({\ga\inv};m,j)\om({\beta \inv};i,m)}\right.\\
+\left.
\frac{(\gamma\beta  w)_j^\vee \gamma
v_j^\vee }{\om({\beta \inv};i,j)
\om({\gamma\inv};m,j)}\right)
\cd s_{im}s_{mj}\beta _i\gamma_m(\beta ^{-1}\gamma^{-1})_j\overset{{\,}_\to}f
\end{multline*}
$$
=4
\left(\sum_{\gamma,\beta \in\Gamma}\frac{\gamma
 v_j^\vee }{\om({\beta \inv};i,m)
\om({\beta \inv};i,j)
\om({\ga\inv};m,j)}\cd Y_{\beta,\ga,m,j}\right)\overset{{\,}_\to}f=0,
$$
where $Y_{\beta,\ga,m,j}$ is given by the following expression
$$Y_{\beta,\ga,m,j}=\left(w_i^\vee \om({\ga\inv};m,j)
-(\beta w)_m^\vee \om({\beta \inv};i,j)
+(\gamma\beta w)_j^\vee
\om({\beta \inv};i,m)\right)\cd
s_{im}s_{mj}\beta _i
\gamma_m(\beta ^{-1}\gamma^{-1})_j.
$$

We  consider the case $f_i,f_m=f^+$,
$f_j=f^-$:
\begin{align*}
 4(R_{im}^wR_{mj}^v-R_{mj}^vR_{ij}^w+R^w_{ij}R_{mi}^v) (\overset{{\,}_\to}f)& =
\left(
\sum_{\beta ,\gamma\in\Gamma}\frac{(w)_i^\vee \beta ^{-1}
v_i^\vee }{\om({\beta\inv};i,m)\om({(\ga\beta)\inv};i,j)}
s_{im}s_{mj}\beta _i\gamma_m
(\gamma\beta )^{-1}_j \right. \\ & -
\sum_{\gamma,\beta \in\Gamma}\frac{ v_m^\vee
w_i^\vee }{\om({\beta\inv};m,j)\om({\ga\inv\beta};i,m)}
s_{mj}s_{ij}(\beta ^{-1}\gamma)_i\beta _m\gamma^{-1}_j \\
& +
\left. \sum_{\gamma,\beta \in\Gamma}\frac{ w_i^\vee
\beta \gamma v_j^\vee }{\om({\beta\inv};i,j)
\om({(\beta\ga)\inv};m,j)}
s_{ij}s_{mi}\gamma_i^{-1}(\beta \gamma)_m\beta _j^{-1}\right)\overset{{\,}_\to}f.
\end{align*}

Performing the  same change of summation indices as in the
previous paragraph, we deduce that the following sum vanishes
$$
\sum_{\gamma,\beta \in\Gamma}\frac{
w^\vee _i(\beta ^{-1}v_i^\vee 
\om({\ga\inv};m,j)
-v_m^\vee \om({(\gamma\beta)\inv};i,j) +\gamma v_j^\vee
\om({\beta\inv};i,m))}{\om({\beta\inv};i,m)
\om({(\gamma\beta)\inv};i,j)\om({\gamma\inv};m,j)}
s_{im}s_{mj}\beta _i
\gamma_m(\beta ^{-1}\gamma^{-1})_j\overset{{\,}_\to}f.
$$

It easy to see that if $f_i=f^-$ or
$f_m=f^-$ then
$R_{im}^wR_{mj}^v\overset{{\,}_\to}f=R_{mj}^vR_{ij}^w
\overset{{\,}_\to}f=R^w_{ij}R_{mi}^v\overset{{\,}_\to}f=0$. Thus, we
have
proved (\ref{BraidDunkl1}).

 Now we prove that
\begin{equation*}
[R_{im}^w,R_{mi}^v]=0.
\end{equation*}
It is easy to see that $[R_{im}^w,R_{mi}^v]\overset{{\,}_\to}f=0$ when
$f_i=f^-$ or $f_m=f^-$. If
$f_i,f_m=f^+$ then $[R_{im}^w,R_{mi}^v]\overset{{\,}_\to}f=0$ is
equivalent to:
\begin{equation}\label{badterm1}
\sum_{\gamma,\beta \in\Gamma}\frac{w_i^\vee \beta \gamma
v_m^\vee -v_m^\vee \beta ^{-1}\gamma^{-1}
w_i^\vee }{\om({\be\inv};i,m)\om({(\beta\gamma\beta)\inv};i,m)}
(\beta \gamma)_m(\gamma\beta )_i^{-1}\overset{{\,}_\to}f=0.
\end{equation}

Fix some $\gamma,\beta \in\Gamma$. Let
$\alpha:=\beta \gamma$, $\beta :=\gamma\beta $. It is
easy to see that the coefficient in front of
$\alpha_m\beta _i^{-1}\overset{{\,}_\to}f$ in (\ref{badterm1}) is equal to
\begin{equation}\label{badterm2}
(w_i^\vee \alpha v_m^\vee -v_m^\vee \beta ^{-1}w_i^\vee )\sum_{\delta\in
Z_\beta }\frac1{\om({\delta\beta};i,m)
\om({\beta\inv\delta\beta};i,m)}
\end{equation}
where $Z_\beta $ is the notation for the centaralizer of
$\beta \in\Gamma$. We  notice that $\alpha$ is conjugate to
$\beta $, hence if $\beta =1$ or $\beta =\zeta$ then
$\alpha=\beta =\beta ^{-1}$ and (\ref{badterm2}) is zero.

Thus we can assume that $\beta \ne 1,\zeta$. In this case $Z_\beta $
is a cyclic group. We  denote the order of $Z_\beta $ by $l$ and
let $\rho$ be a generator  of $Z_\beta $. Then we can assume that
$\beta =\rho^q$ for some $q,0<q<l$.

  Let $a,b\in L$ be the basis in $L$ such that
$\rho a=\epsilon a$, $\rho b=\epsilon^{-1} b$ where $\epsilon$ is
some primitive $l$th root of $1$. Let $x=a^*$ and $y=b^*$. We
make change of variables $\beta u_m\to u_m$ in
(\ref{badterm2}). Let  $z_{mi}=x_iy_m/(x_m y_i)$. Then,
 (\ref{badterm2}) is proportional to
\begin{multline*}
\sum_{p=0}^{l-1}\frac1{\omega(u_i,\delta^pu_m)\omega(\delta^q
u_i,\delta^pu_m)}=\frac1{(x_my_i)^2}\sum_{p=0}^{l-1}\frac1{(\epsilon^{-p}-\epsilon^p
z_{mi})(\epsilon^{-p+q}-\epsilon^{p-q}z_{mi})}=\\
\frac1{(\epsilon^q-\epsilon^{-q})(x_mx_i)^2}\sum_{p=0}^{l-1}
\frac1{z_{mi}-\epsilon^{2q-2p}}-\frac1{z_{mi}-\epsilon^{-2p}}=0.
\end{multline*}

 Finally we show that
\begin{equation*}
[\Theta^{\text{Dunkl}}_i(w),R_{mi}^v]+ [R^w_{im}\Theta^{\text{Dunkl}}_m(v)]=
-\frac{1}{2} \sum_{\ga\in\Gamma} \omega_{L}(\ga w,v) s_{i
m}\ga_{i}\ga_{m}^{-1}
\end{equation*}

If $f_i,f_m=f^{\pm}$ then
the terms in the LHS of the sum below corresponding to
$\gamma$ and $\gamma\zeta$  mutually cancel out, and we deduce
$$-\frac{1}{2} (\sum_{\ga\in\Gamma} \omega_{L}(\ga w,v) s_{i
m}\ga_{i}\ga_{m}^{-1})\overset{{\,}_\to}f=0.$$

In the case $f_i,f_m=f^-$ we know that
$[\Theta^{\text{Dunkl}}_i(w),R_{mi}^v]\overset{{\,}_\to}f=[R_{im}^w,\Theta^{\text{Dunkl}}_m(v)]\overset{{\,}_\to}f=0$.
In the case $f_i=f_m=f^+$ we have
\begin{gather*}
[\Theta^{\text{Dunkl}}_i(w),R_{mi}^v]\overset{{\,}_\to}f=
\frac12\left(\sum_{\gamma\in\Gamma}-(D^w)_i\frac{\gamma
v_i^\vee }{\om({\ga\inv};m,i)}s_{mi}\gamma_m\gamma_i^{-1}+
\frac{v_m^\vee }{\om({\ga\inv};m,i)}
(D^{\gamma^{-1}w})_ms_{mi}\gamma_m\gamma_i^{-1}\right)\overset{{\,}_\to}f\label{11}\\
[R_{im}^w,\Theta^{\text{Dunkl}}_m(v)]\overset{{\,}_\to}f=
\frac12\left(\sum_{\gamma\in\Gamma}-
\frac{w_i^\vee }{\om({\ga\inv};i,m)}
(D^{\gamma^{-1}v})_is_{im}\gamma_i\gamma_m^{-1}+
(D^v)_m\frac{\gamma w_m^\vee }{\om({\ga\inv};i,m)}
s_{im}\gamma_i\gamma_m^{-1}\right)\overset{{\,}_\to}f\label{22}.
\end{gather*}

Fix  a conjugacy class ${\rm C}\subset\Gamma$. 
Then the
coefficient in front of $-\frac12 c_{\beta }$,
$\beta \in{\rm C}$ at  $([\Theta^{\text{Dunkl}}_i(w),R_{mi}^v]+
[R^w_{im}\Theta^{\text{Dunkl}}_m(v)])\overset{{\,}_\to}f$ is equal to
\begin{multline*}\label{33}
\left(\sum_{\gamma\in\Gamma,\beta \in{\rm
C}}\frac{(\beta  w_i^\vee +w_i^\vee )\beta \gamma
v_i^\vee }{\om({\be};i,i)\om({\ga\inv\be\inv};m,i)}
s_{mi}(\beta {\gamma})_m(\gamma^{-1})_i-
\frac{v_m^\vee (\beta \gamma^{-1} w_m^\vee +\gamma^{-1}
w_m^\vee )}{\om({\ga\inv};m,i)\om({\be};m,m)}
s_{mi}\gamma_m(\beta \gamma^{-1})_i\right.\\
-\left.\frac{(\beta  v_m^\vee +v_m^\vee )\beta \gamma
w_m^\vee }{\om({\be};m,m)\om({\ga\inv\be\inv};i,m)}
s_{im}(\beta {\gamma})_i(\gamma^{-1})_m+
\frac{w_i^\vee (\beta \gamma^{-1} v_i^\vee +\gamma^{-1}
v_i^\vee )}{\om({\ga\inv};i,m)\om({\be};i,i)}
s_{im}\gamma_i(\beta \gamma^{-1})_m\right)\overset{{\,}_\to}f.
\end{multline*}

Fix  $\gamma\in\Gamma,\beta \in{\rm C}$. Then the
coefficient in front of
$s_{im}(\beta \gamma)_m\gamma_i^{-1}$ is equal to
\begin{multline*}
F=\frac{(\beta w_i^\vee +w_i^\vee )\beta \gamma
v_i^\vee }{\omega(\be;i,i)\om(\beta \gamma;m,i)}
-\frac{v_m^\vee (\gamma^{-1}w_m^\vee +\gamma^{-1}\beta ^{-1}
w_m^\vee )}{\omega(\beta \gamma;m,i)\omega(\gamma^{-1}\beta \gamma;m,m)}-\\
\frac{(\gamma^{-1}\beta \gamma
v_m^\vee +v_m^\vee )\gamma^{-1}w_m^\vee }{\om({\gamma^{-1}\beta \gamma};m,m)
\om({\ga};i,m)}+ \frac{w_i^\vee (\beta \gamma
v_i^\vee +\gamma
v_i^\vee )}{\om({\gamma^{-1}};i,m)\om({\be};i,i)}.
\end{multline*}

We see that $F=F(u_i,u_m)$ is a homogeneous function in two variables, of
bidegree $(-1,-1)$, that is $F\in H^0({\mathbb{P}} \times
{\mathbb{P}} ,\mathcal{O}(-1)\boxtimes\mathcal{O}(-1))$. It could
have simple poles along the divisors $u_i\sim \beta u_i$,
$u_m\sim\gamma^{-1}\beta \gamma u_m$,
$u_i\sim\beta \gamma u_m$, $u_i\sim\gamma u_m$ where $\sim$
stands for being proportional. But is easy to check that the
residues  actually vanish. We deduce that $F=0$.

 The coefficient in front of $s_{mi}\gamma_m\gamma_i^{-1}\overset{{\,}_\to}f$ in
the part of $([\Theta^{\text{Dunkl}}_i(w),R_{mi}^v]+
[R^w_{im}\Theta^{\text{Dunkl}}_m(v)])\overset{{\,}_\to}f$ that does not contain
coefficients $c_{\gamma}$, $\gamma\in\Gamma$, equals
$$
-\frac{\partial}{\partial w_i}\frac{\gamma v_i^\vee }{\om({\ga\inv};m,i)}
+\frac{v_m^\vee }{\om({\ga\inv};m,i)}
\frac{\partial}{\partial(\gamma^{-1}w)_m}-\frac{w_i^\vee 
}{\om({\gamma^{-1}};i,m)}\frac{\partial}{\partial(\gamma
v)_i}+\frac{\partial}{\partial v_m}\frac{\gamma^{-1}
w_m^\vee }{\om({\gamma^{-1}};i,m)}.
$$
It is easy to show that this expression vanishes since we have
$$w_i^\vee \frac{\partial}{\partial(\gamma v)_i^\vee }-\gamma
v_i^\vee \frac{\partial}{\partial w_i}
=\gamma^{-1}w_m^\vee \frac{\partial}{\partial (\gamma v)_m^\vee }-
v_m^\vee \frac{\partial}{\partial (\gamma^{-1}w)_m}=\omega(w,v)
\eu.$$

We consider the case $f_i=f^-$,
$f_m=f^+$. Then we have $[R_{im}^w,\Theta_m(v)]\overset{{\,}_\to}f=0$
and
\begin{multline*}
[\Theta^{\text{Dunkl}}_i(w),R_{mi}^v]\overset{{\,}_\to}f=\frac12
\sum_{\gamma\in\Gamma}\frac{ w_i^\vee  v_m^\vee - (\gamma
v)_i^\vee (\gamma^{-1}
w)_m^\vee}{\om({\gamma\inv};m,i)}
s_{mi}\gamma_m\gamma_i^{-1}\overset{{\,}_\to}f=
-\sum_{\gamma\in \Gamma}\omega(\gamma w,v)s_{mi}\gamma_i\gamma_m^{-1}
\overset{{\,}_\to}f.
\end{multline*}

The analysis of the case $f_i=f^+$, $f_m=f^-$
is similar.

{\footnotesize{

\smallskip

\noindent
Department of Mathematics,  Massachusetts Institute
of Technology, Cambridge, MA 02139, USA\\
{\tt etingof@math.mit.edu}
\smallskip

\noindent
Department of Mathematics,  Massachusetts Institute
of Technology, Cambridge, MA 02139, USA\\
{\tt wlgan@math.mit.edu}
\smallskip

\noindent
Department of Mathematics, University of Chicago,
Chicago, IL 60637, USA\\
{\tt ginzburg@math.uchicago.edu}
\smallskip

\noindent 
School of Mathematics,  
Institute for Advanced Study, Princeton, NJ 08540, USA\\
{\tt oblomkov@math.ias.edu}

}}

\end{document}